\documentclass{article}    
\usepackage{xy}
\usepackage{array}
\usepackage{amsmath,amsfonts,amssymb}
\newenvironment{Proof}{
{\setlength{\parindent}{0.0in}

{\em Proof.} }
}{
\rule[-2.2mm]{1.5mm}{3.7mm}  
\vspace{0.13in} }
\def\veestrut{\rule{0mm}{3mm}}
\def\mm{\,\o{\veestrut\wedge}\,}\def\jj{\,\o{\veestrut\vee}\,}%
\newtheorem{theorem}{Theorem}
\newtheorem{lemma}[theorem]{Lemma}
\newtheorem{corollary}[theorem]{Corollary}

\newcommand{\Dash}{\rule[.9mm]{1.5cm}{.1mm}\hspace{2mm}}
\newcommand{\EQ}{\mbox{$\;\; = \;\;$}}
\newcommand{\Wavy}{\;\approx\;}
\newcommand{\WAVY}{\;\approx\;}

\newcommand{\INV}[1]{#1^{-1}}
\newcommand{\Compos}{\!\circ\! }
\newcommand{\FROM}{\!:\!}
\newcommand{\TO}{\longrightarrow}
\newcommand{\GOESTO}{\longmapsto}
\newcommand{\Mapnamed}[1]{\stackrel{#1}{\longrightarrow}}
\newcommand{\ISO}{\cong}
\newcommand{\la}{\langle}
\newcommand{\ra}{\rangle}
\newcommand{\Reals}{{\mathbb R}}
\newcommand{\Complex}{{\mathbb C}}
\newcommand{\Integers}{{\mathbb Z}}

\newcommand{\ENDPROOF}{\hspace{0.03in}
\rule[-2.2mm]{1.5mm}{3.7mm}  
\vspace{0.03in}\\}\def\o{\overline}
\def\downstrut{\relax}

\def\wavy{\approx}
\def\Models{\;\models\;}

\def\o{\overline}

\newcommand{\Tri}[3]{\mbox{$\bigtriangleup#1#2#3$}} 
\def\join{\vee}     
\def\meet{\wedge}   

\begin{document}
\begin{center}
    {\bf Approximate satisfaction of identities}\\[0.1in]
                   by\\[0.1in]
      {\bf Walter {
      {Taylor}}} (Boulder, CO)
\end{center}

\begin{center}{Version \today.}\\[0.1in]
\end{center} \vspace{0.1in}

{\small {\bf Abstract.} For a metric space $(A,d)$, and a set
$\Sigma$ of equations, a quantity is introduced that measures how
far continuous operations must deviate from satisfying $\Sigma$ on
$(A,d)$. } \vspace{0.05in}

\tableofcontents 

 \setcounter{section}{-1}

\section{Introduction.}                  \label{sec-intro}

This paper is part of a continuing investigation---see the
author's papers \cite{wtaylor-cots} (1986), \cite{wtaylor-sae}
(2000), and \cite{wtaylor-eri} (2006)---into the {\em
compatibility} relation (see (\ref{eq-space-models-sigma}) below)
between a topological space $A$ and a set $\Sigma$ of equations,
which we will briefly review in \S\ref{sub-compat-backg}.

\subsection{Compatibility---context and background.}
            \label{sub-compat-backg}
In this context, $\Sigma$ typically denotes a set (finite or
infinite) of equations\footnote{%
A (formal) equation is an ordered pair of terms $(\sigma,\tau)$,
more frequently written $\sigma\wavy\tau$. As such it makes no
assertion, but merely presents two terms for consideration. The
actual mathematical assertion is made by the satisfaction relation
$\models$.}, %
which are understood as universally quantified. We usually expect
that $\Sigma$ has a specified similarity type. This means that we
are given a set $T$ and whole numbers $n_t\geq0$ ($t\in T$), such that
for each $t\in T$ there is an operation symbol\footnote{%
In our examples we give the operation symbols familiar names like
$+$ or $\wedge$. These may be thought of as colloquial expressions
for the more formal $F_t$.} %
$F_t$ of arity $n(t)$, and such that the operation symbols of $\Sigma$
are included among these $F_t$.

Given a set $A$ and for each $t\in T$ a function $\o{F_t}\FROM
A^{n(t)} \TO A$ (called an operation), we say that the operations
$\o{F_t}$ {\em satisfy} $\Sigma$ and write
\begin{equation}             \label{eq-opening-set}
        (A,\o{F}_t)_{t\in T} \Models \Sigma,
\end{equation}
if for each equation $\sigma\wavy\tau$ in $\Sigma$, both $\sigma$
and $\tau$ evaluate to the same function when the operations
$\o{F_t}$ are substituted for the symbols $F_t$ appearing in
$\sigma$ and $\tau$. Given a {\em topological space} $A$ and a set
of equations $\Sigma$, we write
\begin{equation}             \label{eq-space-models-sigma}
        A \Models \Sigma,
\end{equation}
and say that $A$ and $\Sigma$ are {\em compatible}, iff there
exist {\em continuous} operations $\o{F_t}$ on $A$ satisfying
$\Sigma$.

While the definitions are simple, the relation
(\ref{eq-space-models-sigma}) remains mysterious. The algebraic
topologists long knew that the $n$-dimensional sphere $S^n$ is
compatible with H-space theory (see \S\ref{sub-sub-S2-epsilon}) if
and only if $n=1,3$ or $7$. For $A=\Reals$, the relation
(\ref{eq-space-models-sigma}) is algorithmically undecidable for
finte $\Sigma$ \cite{wtaylor-eri}; i.e.\ there is no algorithm that
inputs an arbitrary finite $\Sigma$ and outputs the truth value of
(\ref{eq-space-models-sigma}) for $A=\Reals$. In any case,
(\ref{eq-space-models-sigma}) appears to hold only sporadically,
and with no readily discernable pattern.

The mathematical literature contains many scattered examples of
the truth or falsity of specific instances of
(\ref{eq-space-models-sigma}). The author's earlier papers
\cite{wtaylor-cots}, \cite{wtaylor-sae}, \cite{wtaylor-eri},
collectively refer to most of what is known. We will therefore not
attempt to write a list of examples for this introduction. In any
case, many of the known results are recapitulated later in the
paper, as we endeavor to find specific metrical versions of
(\ref{eq-space-models-sigma}) or its negation. The reader is
invited to peruse the two Figures in \S\ref{sub-def-lamb-space} as
a starting point for this information.

\subsection{Metric approximation to compatibility.}
              \label{sub-metric-approx-compat}
In this paper we examine and refine the compatibility relation
(\ref{eq-space-models-sigma}) when the topology of $A$ is given by
a metric $d$. If (\ref{eq-space-models-sigma}) holds, that is if
$A\models\Sigma$, then there is little more to say in the context
of this paper. On the other hand, if (\ref{eq-space-models-sigma})
fails, that is, if $\Sigma$ cannot be modeled on $A$ with
continuous operations, then we ask whether it is possible to model
$\Sigma$ {\em approximately} with respect to the metric $d$. More
precisely, for real $\varepsilon>0$, we write
\begin{align}                 \label{eq-models-varep}
              A \models_{\varepsilon} \Sigma
\end{align}
iff there exist continuous operations $\o F_t$ on $A$ such that,
for each equation $\sigma\wavy\tau$ in $\Sigma$, the terms
$\sigma$ and $\tau$ evaluate to functions that are within
$\varepsilon$ of each other. We seek information about the number
$\lambda_A(\Sigma)$ that is defined as follows:
$\lambda_A(\Sigma)$ {\em is the smallest non-negative\/\footnote{%
If there is any real number satisfying this condition, then there
is a smallest one, by completeness. If there is no such real
number, then we let $\lambda_A(\Sigma)=\infty$.} %
real such that $A \models_{\varepsilon} \Sigma$ for every
$\varepsilon>\lambda_A(\Sigma)$.}

(The reader may check that $\lambda_A$ is characterized by the
validity for all real $\varepsilon$ of the following two
statements: (a) if $0<\varepsilon<\lambda_A(\Sigma)$, then
$A\not\models_{\varepsilon}\Sigma$; (b) if
$\lambda_A(\Sigma)<\varepsilon$, then
$A\models_{\varepsilon}\Sigma$. It therefore makes little
difference whether $\models_{\varepsilon}$ refers to distances
$<\varepsilon$ or $\leq\varepsilon$. Seldom in the paper do we
refer to $\models_{\varepsilon}$ as such, although see
\S\ref{sub-vision},  Problem 6 in \S\ref{sub-problems}, the proof
of Lemma \ref{lemma-comp-two-metrics} in \S\ref{sub-dep-metric},
and \S\ref{subsub-semil-homotop}.)

The paper contains some general facts about $\lambda_A(\Sigma)$.
For example, we have studied how $\lambda_A(\Sigma)$ depends on
deductions, on interpretation of operations by terms, on choice of
metric for A's topology, and so on. We have included a small
catalog of estimates for $\lambda_A(\Sigma)$, along with some
indication of various methods that are sometimes applicable. On
the theoretical side, near the end of the paper we show that if
the space $A$ has a finite triangulation, then simplicial maps can
be used for the approximate satisfaction of $\Sigma$. With this
tool we see that, for $\alpha$ a computable real, the collection
of finite sets $\Sigma$ with $\lambda_A(\Sigma)<\alpha$ is
recursively enumerable (Corollary \ref{cor-less-alpha-re}).

\subsection{The structure of the paper.}
In \S\ref{sec-def-lambda} we state again the definition of
$\lambda_A(\Sigma)$ and some further auxiliary definitions. At
this point, statements of the form $\lambda_A(\Sigma)<M$ or
$\lambda_A(\Sigma)>M$ will be intelligible to the reader; two
charts in \S\ref{sub-def-lamb-space} contain a number of such
results, selected from later in the paper (mostly from
\S\ref{sub-num-ex}).
\S\ref{sec-general} contains theoretical results about
$\lambda_A(\Sigma)$ and its calculation.

The long \S\ref{sub-num-ex} contains our proofs of estimates on
$\lambda_A(\Sigma)$ (both from above and from below).
\S\ref{sub-num-ex} is organized mostly by a rough typology of the
various $\Sigma$'s that can appear: inconsistent, group-theoretic,
lattice-theoretic, and so on. Right at the start of
\S\ref{sub-num-ex} we outline an alternate organization of the
material, in terms of the different methods of proof that are
available.

In \S\ref{sec-top-inv} we investigate the possibility of
eliminating any dependence on the choice of metric, for example by
taking the infimum of $\lambda_A(\Sigma)$-values over all metrics
that yield $A$'s topology and have diameter $1$. Like
\S\ref{sub-num-ex}, \S\ref{sec-top-inv} is mostly a collection of
small results and methods. In \S\ref{sec-prod-vars} we briefly
examine product varieties: how does
$\lambda_A(\Gamma\times\Delta)$ relate to $\lambda_A(\Gamma)$ and
$\lambda_A(\Delta)$? It was known that if $A$ is
product-indecomposable and $A\models \Gamma\times\Delta$, then
$A\models\Gamma$ or $A\models\Delta$. For certain special $A$, we
have an analogous result: if $\lambda_A(\Gamma\times\Delta)$ is
small, then $\lambda_A(\Gamma)$ is small or $\lambda_A(\Delta)$ is
small. \S\ref{sec-simplicial} deals with approximate satisfaction
by simplicial (piecewise linear) maps: if $A$ is a finite
simplicial complex, and if $\Sigma$ is finite, then in the
definition of $\lambda_A(\Sigma)$ we may restrict attention to
simplicial maps (Corollary \ref{cor-simp-lambda} on page
\pageref{cor-simp-lambda}). \S\ref{sec-alg}
contains the enumerability result mentioned above in
\S\ref{sub-metric-approx-compat}.

\S\ref{sec-filters} concerns a certain kind of filter that can be
defined using $\lambda_A$ on the lattice of interpretability types
of varieties. \S\ref{sub-approx-diff} is a very short excursion
into the question of approximate satisfaction by differentiable
operations.

After the brief introduction in \S\ref{sec-def-lambda}, the
remaining sections can be read in almost any order. Even
\S\ref{sec-general} may be skipped; one can return there when
necessary. \S\ref{sec-top-inv} may be read independently as a
collection of estimates, although most of its proofs depend on
\S\ref{sub-num-ex}. \S\ref{sec-filters} requires
\S\ref{sec-prod-vars} for a proof and also for an example, but
otherwise can be understood independently. \S\ref{sec-simplicial}
should probably precede \S\ref{sec-alg}. The other sections are
independent of one another. As for the long \S\ref{sub-num-ex},
its subsections (\S\ref{subsub-lim-zero},
\S\ref{sub-spheres-simple}, and so on) can each be read
independently.

\subsection{Vision and outlook}     \label{sub-vision}
We propose the incompatibility measure $\lambda_A$, along with
$\delta^-_A$ and $\delta^+_A$ of \S\ref{sec-top-inv}, as a new
tool for explaining, systematizing, and perhaps extending known
results on the truth or falsity of the compatibility relation
(\ref{eq-space-models-sigma}). In many cases that we have
investigated, it turns out that the approximate satisfiability
relation
(\ref{eq-models-varep}) is as tractable\footnote%
{For three notable exception to this assertion, see
\S\ref{subsub-group-cyl},
\S\ref{subsub-boolean} and \S\ref{subsub-square-not-exp2}.} %
as the apparently simpler relation (\ref{eq-space-models-sigma}).

It is our hope that some pattern may emerge from this finer
information that will begin to elucidate the overall idea of
compatibility. As one example we mention that the algorithmic or
recursion-theoretic character of the relation $A\models\Sigma$
(say, for finite $\Sigma$ and a finite complex $A$) is not
understood.\footnote{%
Although from \cite{wtaylor-eri} we know that the relation
$\Reals\models\Sigma$ is not algorithmic.} %
On the other hand we now know (Corollary
\ref{cor-less-alpha-re-pairs} of \S\ref{sub-re-ity}) that, for
fixed recursive\footnote{%
The recursive reals include all rational numbers, and almost
all well-known numbers such as $e$, $\pi$ and
so on.} %
$\varepsilon$, the relation
$A\models_{\varepsilon}\Sigma$, with $A$ and $\Sigma$ ranging over
finite complexes and finite theories, is recursively enumerable.
We also know (Corollary \ref{cor-pi-2} of
\S\ref{sub-ar-char-zero}) that for fixed finite $K$, the set of
$\Sigma$ with $\lambda_{|K|}(\Sigma)=0$ is (at worst) a
$\Pi_2$-set.

Finally, we remark that, for a fixed set $\Sigma$ of equations,
one may for example define a topological invariant
$H_{\Sigma}(A)=\delta^-_A(\Sigma)$ (from \S\ref{sec-top-inv}).
It may be that for certain $\Sigma$ the invariant $H_{\Sigma}$
could play some role in topology, for instance in the
classification of spaces $A$. Certainly some cases may be seen in
the paper where a given $H_{\Sigma}$ distinguishes two spaces;
each such non-homeomorphism was, however, already known.

\subsection{Preview of an alternate approach}
We will briefly mention a different metric refinement of the compatibility
relation $(A,d)\models\Sigma$, for a metric space $(A,d)$. It will
be discussed at length in the forthcoming paper \cite{wtaylor-dise}.
The basic idea is that we consider only operations $\o F_t$ satisfying
(\ref{eq-opening-set}), i.e. which model our equations $\Sigma$
exactly, but we allow discontinuities in the $\o F_t$. In
\cite{wtaylor-dise} we introduce a real quantity $\mu_n(A,\Sigma)$,
which measures how far at least one operation $\o F_t$ must deviate
from continuity, as we now very briefly describe.

For simplicity here, we
mention only the version that corresponds to uniform continuity.
Let $(A,d)$ and $(B,e)$ be metric spaces.
Given $\o F\FROM B\TO A$ (not necessarily continuous) and
$\delta,\varepsilon>0$, we say that $\o F$ is {\em
$(\delta,\varepsilon)$-constrained} it satisfies: for all $b,b'\in
B$, if $e(b,b')<\delta$, then $d(\o F(b),\o F(b'))<\varepsilon$.

For $0<\delta_0\leq\delta_n$, we say that $\o F$ is $n$-constrained
by $(\delta_0,\delta_n)$ iff
there exist $\delta_1,\ldots,\delta_{n-1}$ with $0<\delta_0
\leq\delta_1\allowbreak\leq\cdots\allowbreak\leq\delta_n$ such
that $\o F$ is $(\delta_0,\delta_1)$-constrained and
$(\delta_1,\delta_2)$-constrained, and so on, up to
$(\delta_{n-1},\delta_n)$-constrained.
(If $\o F$ is uniformly continuous, then for every $n$ and for
every $\varepsilon>0$
there exists $\delta>0$ so that $\o F$ is $n$-constrained by
$(\delta,\varepsilon)$.)

We now define
\begin{align*}
        A \Models{\!\!\!}_n^{\varepsilon}\; \Sigma
\end{align*}
to mean that {\em there exists an algebra\/ $\mathbf A=( A,\o
F_t)_{t\in T}$ modeling\/ $\Sigma$ and a real number\/
$\delta_0>0$ such that each\/ $\o F_t$ is $n$-constrained by\/
$(\delta_0,\varepsilon)$.} \vspace{0.1in}
We define
\begin{align*}
    \mu_n(A,\Sigma) \EQ
       \inf \,\{\varepsilon:A \!\Models{\!\!\!}_n^{\varepsilon}\;
       \Sigma\}.
\end{align*}
It is not hard to see that
\begin{align*}
    0\;\leq\;\mu_1(A,\Sigma) \;\leq\;\mu_2(A,\Sigma)
    \;\leq\; \cdots \;\leq\;\text{diam}(A).
\end{align*}

For several interesting equational theories $\Sigma$, the behavior
of $\mu_n(A,\Sigma)$, especially $\mu_1(A,\Sigma)$, is somewhat parallel
to---but not always the same as---that of the $\lambda_A(\Sigma)$ studied here.
We view the
utility of the concept of  $\mu_n(A,\Sigma)$ as similar to that proposed
for  $\lambda_A(\Sigma)$ in \S\ref{sub-vision}: to further our understanding
of compatibility in general.   The detailed
presentation will appear in \cite{wtaylor-dise}; there will be no further
mention of $\mu_n(A,\Sigma)$ in this paper.

\subsection{Problems}                    \label{sub-problems}
We have written few explicit problems into the body of this paper,
although obviously many things are not yet known. Here we collect
a few ideas for further study. In a sense the most important
problem is the vague one of elucidating the notion of
compatibility between a topological space and a set of equations.
However a list of some more focused problems may be of some use.
\vspace{0.1in}

\hspace{-\parindent}%
{\bf Problem 1.}
 Of course, every approximation given for some
$\lambda_A(\Sigma)$ calls for a better approximation, if not for
the exact value. Specific instances of this problem are implicit
throughout \S\ref{sub-num-ex}.\vspace{0.1in}

\hspace{-\parindent}%
{\bf Problem 2.} We mention an incompatibility result for Boolean
algebras in \S\ref{subsub-boolean}.  There we ask whether this
result can be
sharpened to a positive lower estimate on $\lambda_A(B)$ for $B$
some equations true in Boolean algebras, and $A$ some appropriate
space or spaces.\vspace{0.1in}

\hspace{-\parindent}%
{\bf Problem 3.} Are there theorems saying that, under certain
conditions, $\lambda_A(\Sigma)=0$ implies
$A\models\Sigma$?\vspace{0.1in}

\hspace{-\parindent}%
{\bf Problem 4.}  Every approximation that one discovers for some
$\delta^-(\Sigma)$ or some $\delta^+(\Sigma)$ calls for a better
approximation, if not for the exact value. Specific instances of
this problem are implicit throughout
\S\ref{sub-estim-delta-minus}.\vspace{0.1in}

\hspace{-\parindent}%
{\bf Problem 5.} Does $\delta^-$ take any values besides $0$,
$0.5$ and $1$? (See the various calculations of $\delta^-$ in
\S\ref{sub-estim-delta-minus}. Cf.\ Theorem \ref{th-delt-plus} of
\S\ref{sub-est-delt-plu} for the situation with
$\delta^+$.)\vspace{0.1in}

\hspace{-\parindent}%
{\bf Problem 6.} The recursive enumerations of \S\ref{sub-re-ity}
(mentioned in \S\ref{sub-vision}) essentially rely on blind luck:
try {\em all\/} possible piecewise linear operations (in a certain
appropriate finite collection), and check each one (using Tarski's
algorithm) to see if it yields $\models_{\varepsilon}$. Is there
an algorithm that proceeds more systematically toward the desired
results?\vspace{0.1in}

\hspace{-\parindent}%
{\bf Problem 7.} Corollary \ref{cor-pi-2} of
\S\ref{sub-ar-char-zero} says that for a fixed finite simplicial
complex $K$, the set of $\Sigma$ with $\lambda_{|K|}(\Sigma)=0$ is
a $\Pi_2$-set. Is it actually of a simpler arithmetic
character?\vspace{0.1in}

\hspace{-\parindent}%
{\bf Problem 8.} In \S\ref{sec-filters} we define $\mathcal L_Z$,
for a metric space $Z$, to be the class of deductively closed
theories $\Sigma^{\star}$ such that $\lambda_Z(\Sigma^{\star})>0$.
Does $\mathcal L_Z$ satisfy the finiteness condition of Mal'tsev
conditions \cite{wtaylor-cmc}? Namely if $\Sigma$ is in  $\mathcal
L_Z$, does $\Sigma$ have a finite subset that is also in $\mathcal
L_Z$? If this condition holds, and if $\mathcal L_Z$ is also
closed under the formation of product varieties, then $\mathcal
L_Z$ is definable by a Mal'tsev condition. (The latter condition
holds for $Z=[0,1]$, as proved in \S\ref{sec-filters}.) \vspace{0.1in}

\hspace{-\parindent}%
{\bf Problem 9.} If the previous problem has a positive solution
for $Z=[0,1]$, then $\mathcal L_Z$ is definable by a Mal'tsev
condition. The problem here is to give an explicit Mal'tsev
condition for $\mathcal L_{[0,1]}$. \vspace{0.1in}

\hspace{-\parindent}%
{\bf Problem 10.} Consider all finite equation-sets $\Sigma$ in a
fixed recursive similarity type that has infinitely many operation
symbols of each finite arity. Is the set $K_1$ of those $\Sigma$
that are analytic-compatible with $\Reals$ recursively inseparable
from the set $K_0$ of those $\Sigma$ that satisfy
$\lambda_{\Reals}(\Sigma)>0$? (In \cite{wtaylor-eri} we proved the
inseparability of $K_0$ from those $\Sigma$ that are not
compatible with $\Reals$.)\vspace{0.1in}

\hspace{-\parindent}%
{\bf Problem 11.} Is there any set $\Sigma$ of equations for which
the invariant $\mu_{\Sigma}$ (defined at the end of
\S\ref{sub-vision}) is an applicable topological invariant? Can
any connection be found between $\mu_{\Sigma}$ and classical
invariants such as cohomology and dimension?

\subsection{Acknowledgments} We should like to thank
Matt Insall, Benjamin Passer, Andrew Priest, Wlodek Charatonik
and Rob Roe of Missouri University of Science and Technology, who
supplied some useful comments on an early draft of this paper.

%
%
%
%
%
%
%
%
%

\section{The quantity $\lambda_A(\Sigma)$ for equations $\Sigma$ and a
 metric space $(A,d)$.}         \label{sec-def-lambda}

In \S\ref{sub-metric} --- indeed, in a large part of the paper
--- we shall begin with a topological space $A$ that has been
provided with a specific metric (or at least a pseudometric) $d$. At
the same time, we consider a set $\Sigma$ of equations (in an
arbitrary similarity type).

\subsection{A pseudometric $\lambda_{\mathbf A}$ for term operations
on a topological algebra $\mathbf A$.} \label{sub-lambda-alg}
\label{sub-metric}

For our given  $(A,d)$, we are generally interested in the
possible existence of topological algebras
 $\mathbf A = (A,\o F_t)_{t\in T}$ that model
$\Sigma$ approximately. In \S\ref{sub-lambda-alg} we
consider a single such $\mathbf A$, and define a
measure of how closely $\mathbf A$ approximates $\Sigma$. The
only constraint on
$\mathbf A$ is that its similarity type should include each
operation symbol occurring in the equations $\Sigma$ that we wish
to approximate.

For any term $\sigma=\sigma(x_0,x_1,\ldots)$ (of the appropriate
similarity type), let $\sigma^{\mathbf A}\FROM A^{\omega}\TO A$
(sometimes also denoted by $\o\sigma$) be the corresponding term
operation that is recursively defined in the usual way. For two such
terms $\sigma,\tau$, we define
\begin{align}                  \label{def-lambda-A}
         \lambda_{\mathbf A}(\sigma,\tau) \EQ
              \sup \bigl\{\,d(\sigma^{\mathbf A}({\mathbf a}),
              \tau^{\mathbf A}({\mathbf a}))
              \,:\, {\mathbf a}\in A^{\omega} \bigr\}
              \;\in\; \Reals^{\geq0}\cup\{\infty\}.
\end{align}
(If there is some need to specify $d$, we may write
$\lambda_{(\mathbf A,d)}(\sigma,\tau)$.) Clearly $\lambda_{\mathbf
A}$ is an $L_{\infty}$-type distance, and as such it is a
pseudometric on the set of all terms $\sigma$ (of the appropriate
type).

Clearly, if $d$ is a metric, then $\lambda_{\mathbf A}$ is a metric
(in the extended sense) as far as term operations are concerned:
$\lambda_{\mathbf A}(\sigma,\tau)=0$ iff $\sigma^{\mathbf A} =
\tau^{\mathbf A}$. From this, it is immediate that
\begin{align*}
  {\mathbf A}\;\models\;
\,\sigma\wavy\tau \quad\quad\text{if and only if}\quad\quad
   \lambda_{\mathbf A}(\sigma,\tau)\,=\,0.
\end{align*}
Thus $\lambda_{\mathbf A}(\sigma,\tau)$ may be considered a measure
of the failure of $\mathbf A$ to model $\sigma\wavy\tau$.

For a (finite or infinite) set $\Sigma$ of equations (of the
appropriate type), we may extend the previous definition to
\begin{align*}
    \lambda_{\mathbf A}(\Sigma) \EQ \sup\bigl\{\lambda_{\mathbf A}
      (\sigma,\tau)\,:\, \sigma\wavy\tau\,\in \Sigma \bigr\}
      \;\in\; \Reals^{\geq0}\cup\{\infty\}.
\end{align*}
As before,
\begin{align*}
  {\mathbf A}\;\models\;
\,\Sigma \quad\quad\text{if and only if}\quad\quad
   \lambda_{\mathbf A}(\Sigma)\,=\,0.
\end{align*}

\subsection{Definition of $\lambda_A$ for a space $A$.}
         \label{sub-def-lamb-space}

Considering now $(A,d)$ simply as a metric space---in other words
not supposing any operations given on $A$ in advance---we may
define the minimax quantity
\begin{align*}
     \lambda_A(\Sigma) \EQ \inf\bigl \{\lambda_{\mathbf
     A}(\Sigma) \,:\,{\mathbf A} =(A;\o F_t)_{t\in T},
       \;\text{$\o F_t$ any continuous operations} \bigr\}.
\end{align*}
(Again, we may write $ \lambda_{(A,d)}(\Sigma)$ if there is a need
to specify $d$.)

\def\upstrut{\rule{0mm}{7mm}}
\def\downstrut{\rule[-3.1mm]{0mm}{1mm}}
\begin{figure}       \label{fig-aa}
\begin{center}
 \begin{tabular}{l|m{1.5cm}|m{1.5cm}|m{1.5cm}|
                          m{1.8cm}|m{1.8cm}|m{1.5cm}|}
      \multicolumn{7}{l}{\hspace*{0.53in}$(0,1)$\hspace{0.49in}
                      $\Reals$\upstrut\downstrut
      \hspace{0.57in} $\Reals^2$ \hspace{0.54in}$Y$
      \hspace{0.64in}$Y_{\varepsilon}$\hspace{0.62in}
      $S^n$}\\ \cline{2-7}
 $\upstrut\Sigma\Lambda$ & $\models$ & $\models$ &
          $\models$ & $\models$ & $\models$ &
     $\lambda=1$\hspace{0.2cm}\hspace*{\fill}
          \S\ref{sub-sub-Sn}\\\cline{2-7}
 $\upstrut\Lambda$ & $\models$ & $\models$ & $\models$ &
               $\lambda \geq 0.144$\hspace{0.2cm}
               \hspace*{\fill}\S\ref{subsub-Y}
                        & $0<\lambda<\varepsilon$
               \hspace{0.2cm}\hspace*{\fill}\S\ref{subsub-Y-new-metric} &
               $\lambda=1$\hspace{0.2cm}
               \hspace*{\fill}\S\ref{sub-sub-Sn} \\\cline{2-7}
 $\upstrut\Sigma\Lambda_0$ & $\lambda=0$ \hspace{0.2cm}
 \hspace*{\fill} \S\ref{subsub-depend-metric}& $\lambda=\infty$\hspace{0.2cm}
          \hspace*{\fill}\S\ref{subsub-depend-metric}&  &   &  &
               $\lambda=1$\hspace{0.2cm}
               \hspace*{\fill}\S\ref{sub-sub-Sn} \\ \cline{2-7}
 $\upstrut\Gamma$ & $\models$ & $\models$ & $\models$ &
            &  &
               $\lambda=1$\hspace{0.2cm}
               \hspace*{\fill}\S\ref{sub-sub-Sn} \\\cline{2-7}
 $\upstrut\Gamma_2$ & &  $\lambda=\infty$\hspace{0.2cm}
 \hspace*{\fill}\S\ref{subsub-not-exp2}
          & \hspace*{0.1cm}{\bf ?}\hspace{0.4cm}
 \hspace*{\fill}\S\ref{subsub-square-not-exp2} &  &  &
               $\lambda=1$\hspace{0.2cm}
               \hspace*{\fill}\S\ref{sub-sub-Sn}  \\\cline{2-7}
 $\upstrut H$ & $\models$ & $\models$ & $\models$ & $\models$
                   & $\models$ &
               $\lambda=1$\hspace{0.2cm}\hspace*{\fill}
               \S\ref{sub-sub-Sn} \\
                \cline{2-7}
        \end{tabular}
 \end{center}

For the given spaces and theories, $A$ is compatible with $\Sigma$
if and only if $\models$ appears in the corresponding entry of the
chart. A blank spot means that we know that $A\models\Sigma$ fails
but we have not attempted to estimate the quantity
$\lambda_A(\Sigma)$. The single question-mark indicates a
$\lambda$-value that we have been unable to estimate.
\begin{tabbing}
Xxxx \= Xxxxxxxxxx  \=   \kill
  \> Theories $\Sigma$:  \> $\Sigma\Lambda$, semilattices;
      $\Sigma\Lambda_0$: with zero; $\Lambda$
        lattices; $\Gamma$ \\
        \>\> groups; $\Gamma_2$, groups of exponent 2; $H$,
           H-spaces. \\[0.2cm]
  \> Spaces A: \> $(0,1)$, $\Reals$, $\Reals^2$, with the
         usual metrics; $Y$ a 120-degree triode \\
         \> \>  with usual planar  metric; $Y_{\varepsilon}$,
         triode with a specific exotic  \\
         \>\> metric; $S^n$, the $n$-sphere (with
         $n\neq 0,1,3,7$). Except          for\\
         \>\>  $\Reals$ and $\Reals^2$, all metrics
                    scaled to a diameter of 1.
\end{tabbing}
   \caption{Some estimates of $\lambda_A(\Sigma)$.}
   \label{fig-chart-one}
 \end{figure}

This $\lambda_A$---defined for a metric space $A$ either here or
by properties (a) and (b) of \S\ref{sub-metric-approx-compat}---is
the central object of study in this paper. A significant part of
our work lies simply in estimating $\lambda$ for enough cases to
illustrate the sort of variation that $\lambda$ can exhibit, and to
illustrate the methods that can be used for such estimates.
Figures \ref{fig-chart-one} and \ref{fig-chart-two} summarize some
of our numerical findings. (For the exact statements of the
relevant results, see the sections referenced in the two Figures.)

\begin{figure}         \label{fig-bb}
 \begin{center}
 \begin{tabular}{l|m{1.8cm}|m{1.8cm}|m{1.8cm}|
                          m{1.8cm}|m{2.4cm}|}
      \multicolumn{6}{l}{\hspace*{1.05in}$\Reals$\hspace{0.59in}
                      $[0,1]$\upstrut\downstrut
      \hspace{0.44in} $[0,1]_{\mbox{?}}$ \hspace{0.42in}$[0,1]^k$
      \hspace{0.48in}$[0,1]^k_{\varepsilon}$}\\ \cline{2-6}
 $\upstrut I$ & $\lambda=0$ \hspace{0.2cm}
 \hspace*{\fill} \S\ref{sub-sub-no-wavy} & $\lambda=0$ \hspace{0.2cm}
 \hspace*{\fill} \S\ref{sub-sub-no-wavy}&
          $\lambda=0$ \hspace{0.2cm}
 \hspace*{\fill} \S\ref{sub-sub-no-wavy} & $\lambda=0$ \hspace{0.2cm}
 \hspace*{\fill} \S\ref{sub-sub-no-wavy} & $\lambda=0$ \hspace{0.2cm}
 \hspace*{\fill} \S\ref{sub-sub-no-wavy}\\ \cline{2-6}
 $\upstrut\Sigma_C$ & & $\lambda=0$ \hspace{0.2cm}
 \hspace*{\fill} \S\ref{subsub-lim-zero-consis} & &
               &
             \\\cline{2-6}
 $\upstrut\Gamma$ & $\models$ & $\lambda=0.5$\hspace{0.2cm}
          \hspace*{\fill}\S\ref{subsub=group-f.p.}& $\lambda\geq0.5$\hspace{0.2cm}
          \hspace*{\fill}\S\ref{subsub=group-f.p.} &$\lambda=0.5$\hspace{0.2cm}
          \hspace*{\fill}\S\ref{subsub=group-f.p.}   &$\lambda\geq0.5$\hspace{0.2cm}
          \hspace*{\fill}\S\ref{subsub=group-f.p.}
               \\ \cline{2-6}
 $\upstrut\Gamma_2$ & $\lambda=\infty$\hspace{0.2cm}
          \hspace*{\fill}\S\ref{subsub-not-exp2}
                        & $\lambda=0.5$\hspace{0.2cm}
          \hspace*{\fill}\S\ref{subsub=group-f.p.}
                   & $\lambda\geq0.5$\hspace{0.2cm}
          \hspace*{\fill}\S\ref{subsub=group-f.p.}
          &$\lambda=0.5$\hspace{0.2cm}
          \hspace*{\fill}\S\ref{subsub=group-f.p.}
           &$\lambda\geq0.5$\hspace{0.2cm}
          \hspace*{\fill}\S\ref{subsub=group-f.p.}
              \\\cline{2-6}
 $\upstrut INJ$ & & $\lambda=0.5$\hspace{0.2cm}
          \hspace*{\fill}\S\ref{subsub-inf-dim}&
                 $\lambda\geq0.5$\hspace{0.2cm}
          \hspace*{\fill}\S\ref{subsub-inf-dim}
           &$\lambda=0.354$\hspace{0.2cm}
          ($k=2$)
          \hspace*{\fill}\S\ref{subsub-square-no-inj}
             &$0<\lambda<\varepsilon
                       $\hspace{0.2cm}
          \hspace*{\fill}\S\ref{sub-delta-square-embed-square}
                \\\cline{2-6}
 $\upstrut INJ_{m,k}$ & & $\lambda=0.5$\hspace{0.2cm}
          \hspace*{\fill}\S\ref{sub-genzations-delta}
          & $\lambda\geq0.5$\hspace{0.2cm}
          \hspace*{\fill}\S\ref{sub-genzations-delta}
           &$\lambda=0.5$\hspace{0.2cm}
          \hspace*{\fill}\S\ref{sub-genzations-delta}
            &$0<\lambda<\varepsilon
          $\hspace{0.2cm}
          \hspace*{\fill}\S\ref{sub-genzations-delta}
               \\
                \cline{2-6}%
 $\upstrut Set^{[n]}$ & & &  \upstrut $\lambda=0.5$
               ($n=2$)\hspace{0.2cm} \hspace*{\fill}
               \S\ref{subsub-approx-cubes} & $\models$ if $n|k$
                   &  \begin{minipage}[t] {1.2in}$0<\lambda<\varepsilon$
                        \end{minipage}
               \hspace*{0.1cm} \begin{minipage}[t] {0.6in}
               \hspace*{\fill}($k\geq n$)\end{minipage}
                          \S\ref{subsub-approx-cubes}
               \\
                \cline{2-6}
        \end{tabular}
 \end{center}

For the given spaces and theories, $A$ is compatible with $\Sigma$
if and only if $\models$ appears in the corresponding entry of the
chart. A blank spot means that we know that $A\models\Sigma$
fails, but we have not attempted to estimate the quantity
$\lambda_A(\Sigma)$.

 \begin{tabbing}
Xxxx \= Xxxxxxxxxx  \=   \kill
  \> Theories $\Sigma$:  \> $I$ is an inconsistent theory from
  \S\ref{sub-sub-no-wavy}; $\Sigma_C$ is from
  \S\ref{subsub-lim-zero-consis};\\
  \>\> $\Gamma$ ($\Gamma_2$) is groups (of exponent 2); $INJ_{m,k}$
    comprises equations\\
    \>\> defining an embedding of $A^m$ into
     $A^k$ (with $m>k$); $INJ$ \\
     \>\> is the special case of
     $m=2$, $k=1$; $Set^{[n]}$ is the theory\\
     \>\> whose models are
     $n^{\mbox{th}}$ powers of sets.
   \\[0.2cm]
  \> Spaces A: \> $[0,1]_{\mbox{?}}$ is $[0,1]$ with an {\em
  arbitrary} diam-1  metric for   the usual\\
  \>\>  topology; $[0,1]^k_{\varepsilon}$ is $[0,1]^k$ with a
     specific exotic\\
     \>\> metric for the usual topology --- see
              \S\ref{sub-genzations-delta}.
\end{tabbing}
 \caption{Further estimates of $\lambda_A(\Sigma)$}\label{fig-chart-two}
\end{figure}

We clearly have: if $A\models\Sigma$, i.e. if $A$ is compatible
with $\Sigma$, then $\lambda_A(\Sigma)=0.$ As may be seen in
\S\ref{subsub-lim-zero}, \S\ref{subsub-depend-metric} and
\S\ref{subsub-lim-zero-consis}, the converse is false: {\em one
may find incompatible  $A$ and $\Sigma$ that nevertheless obey the
following property: for every real $\varepsilon>0$ there are
continuous operations $\o F^{\;\varepsilon}_t$ ($t\in T$) on the
space $A$ satisfying the equations $\Sigma$ within $\varepsilon$.}

\section{General and introductory remarks about $\lambda_A(\Sigma)$.}
                 \label{sec-general}
\subsection{$\lambda_A(\Sigma)$ and the radius and diameter of $A$.}
\label{subsub-lambda}

For $(A,d)$ any metric space, we let
\begin{align*}
           \text{diam}((A,d)) &\EQ \sup\,\{\,d(a_0,a_1)\,:\,
                          a_0,a_1\in A\,\}\\
           \text{radius}((A,d)) &\EQ \inf \,\{ r\,:\,
               A\subseteq\text{some ball of radius r}\};
\end{align*} \label{remarks-diam-radius}
When the context permits, we write $\text{diam}(A)$ and
$\text{radius}(A)$. Clearly $$\text{radius}(A)
\,\leq\,\text{diam}(A) \allowbreak \,\leq\,2\cdot
\text{radius}(A).$$ For some special spaces, such as an
$n$-dimensional cube, the inequality on the right is an equation;
for some other spaces, such as a sphere of dimension $n$, the
inequality on the left is an equation.

Suppose that $\Sigma$ is a set of equations (in variables $x_i$
($i\in\omega$)) that does not contain $x_i\wavy x_j$ for any $i\neq
j$, and that $(A,d)$ is any metric space of finite radius $R$. We
shall see that $\lambda_A(\Sigma)\leq R$.

For each $\varepsilon>0$, there exists $a\in A$ such that
$d(a,x)\leq R+\varepsilon$ for every $x\in A$.
Consider the topological algebra $\mathbf A \EQ (A,\o F_t)_{t\in T}$
of type appropriate to $\Sigma$ that is defined as follows: every
$\o F_t$ is a constant operation with value $a$. It is immediate
that every equation of $\Sigma$ holds within $R+\varepsilon$; hence
$\lambda_A(\Sigma)\,\leq\,R+\varepsilon$. Since this holds for every
positive $\varepsilon$, we have $\lambda_A(\Sigma)\,\leq\,R$; in
other words $\lambda_A(\Sigma)\,\leq\,\text{radius}(A)$.

If $\Sigma$ is any consistent\footnote{%
In the context of equational logic, we define a theory $\Sigma$ to
be {\em consistent} iff $\Sigma$ has a model of more than one
element. In
other words, $x_0\wavy x_1$ is not a consequence of $\Sigma$.} %
theory, then $\Sigma^{\star}$ will not contain any equation
$x_i\wavy x_j$ with $i\neq j$; hence if $\Sigma$ is consistent, then
$\lambda_A(\Sigma^{\star})\leq \text{radius}(A)$.

The reader may easily check that
\[
      \lambda_A(x_0\wavy x_1) \EQ \text{diam}(A).
\]

In particular, $\lambda_A$ cannot take any value strictly between
$\text{radius}(A)$ and $\text{diam}(A)$.

\subsection{Dependence upon deductions.} \label{sub-dep-deduct}

It is obvious from the definition that {\em if\/
$\Sigma\subseteq\Sigma'$, then
$\lambda_A(\Sigma)\,\leq\,\lambda_A(\Sigma')$.} In fact it is
possible to have
\begin{align}    \label{eq-ineq-deduce}
       \lambda_A(\Sigma)\,<\,\lambda_A(\Sigma'),
\end{align}
even though every equation of $\Sigma'$ is a logical consequence of
$\Sigma$. For a simple example, let $A$ be a space with
$\text{radius}(A)<\text{diam}(A)$, take $\Sigma$ to be an
inconsistent theory containing no equation of the form $x_i\wavy
x_j$, and let $\Sigma' = \Sigma\cup\{x_0\wavy x_1\}$. According to
\S\ref{subsub-lambda}, we have
\[
     \lambda_A(\Sigma)\;\leq\; \text{radius}(A)\;<\;\text{diam}(A)
                \EQ \lambda_A(\Sigma').
\]
Thus in fact $\lambda_A(\Sigma)$ is not a logical invariant of
$\Sigma$; in other words, not an invariant of the equational class
defined by $\Sigma$. (For another failure of invariance under
logical deduction, see \S\ref{subsub-dep-ded}.)

One may obtain by fiat a logical invariant of $\Sigma$, as
follows. Let $\Sigma^{\star}$ stand for the set of logical
consequences of $\Sigma$ (in a similarity type that is defined by
the operation symbols appearing in $\Sigma$). If we consider only
quantities of the form $\lambda_A(\Sigma^{\star})$, then we are
indeed considering a logical invariant. The results of
\S\ref{subsub-wavy}, most of \S\ref{sub-spheres-simple}, and
\S\ref{subsub-distort} turn out to be essentially of this form.

There are, however, good reasons to consider the logical
non-invariant $\lambda_A(\Sigma)$ in its own right. In some cases
we can make an upper estimate $\lambda_A(\Sigma)<K$ for a small
finite $\Sigma$, but are unable to extend that estimate to all the
logical consequences of $\Sigma$. (An example may be seen in the
proof of Part (i) of Theorem \ref{th-m=nk} in
\S\ref{sub-genzations-delta}.) This at least points to a further
problem.

As for lower estimates --- $\lambda_A(\Sigma)>K$ ---  for some
familiar finite axiom systems $\Sigma$ it happens that
\begin{itemize}
  \item[(i)] we do not know a positive lower bound on
           $\lambda_A(\Sigma)$;
  \item[(ii)] we can prove $\lambda_A(\Sigma\cup\Gamma)>K$ for a
           certain $K>0$ and a certain finite set $\Gamma$ of
           consequences of $\Sigma$;
  \item[(iii)]  therefore $\lambda_A(\Sigma^{\star}) >K$ by
           (\ref{eq-ineq-deduce}).
\end{itemize}
For instance, consider the lower bounds for $\lambda_Y(\Lambda)$
and $\lambda_{\Reals}(\Lambda_0)$ that appear in  \S\ref{subsub-Y}
and \S\ref{subsub-depend-metric}, respectively. (Here $\Lambda$
(resp. $\Lambda_0$) is lattice theory (resp. with $0$), and $Y$ is
a triode (one-dimensional Y-shaped compact subset of a plane).)
These results stem from Lemma \ref{lem-approx-lat} of
\S\ref{subsub-approx-lattice}; a close examination of that lemma reveals that
its assumptions involve equations that are redundant for lattice
theory. So, strictly speaking, \S\ref{subsub-Y} does not yield a
lower bound for the usual equational-axiomatic formulation of
lattice theory, although it does yield a lower bound for lattice
theory thought of as the equations true in all lattices. Similarly
the lower bound in \S\ref{subsub=group-f.p.}, requires certain
consequences of the axioms of group theory. (Meanwhile, the group
axioms appear to be neither necessary nor sufficient, in
themselves, to yield this estimate.) Similar remarks could be made
about the estimates appearing in \S\ref{subsub-not-exp2} and
\S\ref{subsub-inf-dim}.

Another example of a lower bound that apparently requires
a redundant extension of the original axiom system is Theorem
\ref{th-low-est-squares} of \S\ref{subsub-approx-cubes}; its proof
uses Equations
(\ref{eq-square-deriv-a}--\ref{eq-square-deriv-b}), which are
(logically) redundant in that context.

Among our non-trivial estimates on $\lambda_A(\Sigma^{\star})$ are
lower bounds that arise as in (iii) above, and the very crude
upper bound $\text{radius}(A)$ that occurs in
\S\ref{subsub-lambda} (for consistent $\Sigma$). Some more
interesting upper bounds on $\lambda_A(\Sigma^{\star})$ may be
found in \S\ref{subsub-distort} and \S\ref{subsub-Y-new-metric}.

For instance, in \S\ref{subsub-lambda-homotopy} and
\S\ref{subsub-depend-metric} below, we present examples where we
are able to compute that
$\lambda_A(\Sigma)\,\geq\,\text{diam}(A)$. In these cases, it
clearly follows that $\lambda_A(\Sigma)\EQ
\lambda_A(\Sigma^{\star})\allowbreak\EQ\allowbreak\text{diam}(A)$.

\subsection{Monotonicity of $\lambda_A$ under interpretability.}
               \label{sub-interpretable}
In 1974, W. D. Neumann \cite{neumann} introduced\footnote{%
A number of interpretability notions were already current,
especially due to A. Tarski and G. Birkhoff, including some that
were much like Neumann's. Nevertheless his emphasis on the
resulting ordering of theories (varieties) was new in the
equational context.} %
a quasi-ordering of equational theories, known as
interpretability. The interpretability of $\Gamma$ in $\Delta$,
denoted here $\Gamma\leq\Delta$, is defined as follows.

Let us suppose that the operation symbols of $\Gamma$ are $F_t$
($t\in T$). (These are not necessarily the operation symbols of
$\Sigma$.) Given terms\footnote{%
If $\Sigma$ and $\Gamma$ have disjoint sets of operation symbols,
we may express this relationship informally as ``$F_t\,=\,
\alpha_t\,.$'' For an example of this way of expressing an
interpretation, see Equations
(\ref{eq-interp-log-a}--\ref{eq-interp-log-c}) of
\S\ref{subsub-l-groups}.} %
 $\alpha_t$ ($t\in T$) in the language of $\Sigma$,
we define, for each term $\sigma$ in the language of $\Gamma$, a
term $\sigma^{\star}$ in the language of $\Sigma$. The definition
is by recursion in the length of $\sigma$:
\begin{align}                  \label{eq-star-variable}
                x^{\star} &\EQ x \quad\quad \text{($x$ any variable)},\\
                F_t(\tau_1,\ldots,\tau_{n(t)})^{\star}
                        &\EQ \alpha_t(\tau_1^{\star},\ldots,
                                 \tau_{n(t)}^{\star})
                     \label{eq-star-composite}
\end{align}
(where (\ref{eq-star-composite}) is formally defined by the
logical notion of {\em simultaneous substitution}). If necessary,
we may refer to $\alpha_t$ as the {\bf interpreting term for}
$F_t$ and to $\sigma^{\star}$ as the {\bf interpreting term for}
$\sigma$.

Following \cite{neumann} and \cite[page 1]{ogwt-mem} we now define
$\Gamma\,\leq\,\Sigma$ ($\Gamma$ is {\bf interpretable} in
$\Sigma$) to mean that there exist terms $\alpha_t$ ($t\in T$)
such that for all $\sigma\wavy\tau\in\Gamma$ we have
$\sigma^{\star}\wavy\tau^{\star}\in\Sigma$.

(It is important to point out that this definition differs subtly
from earlier versions. In most of the earlier contexts, such as
\cite{ogwt-mem}, one presented the definition of interpretability
in terms of the model classes (varieties) $\text{Mod }\Sigma$ and
$\text{Mod }\Gamma$, and thus one worked implicitly with the
deductively closed theories $\Sigma^{\star}$ and $\Gamma^{\star}$.
In such a context it is immaterial whether one says
``$\sigma^{\star}\wavy\tau^{\star}\in\Sigma$'' or
``$\sigma^{\star}\wavy\tau^{\star}$ is provable from $\Sigma$.''
(And thus, it may emphasized, all the older results may be
interpreted as valid under the definition proposed here.)  Here,
however, $\lambda_A(\Sigma)$ is not invariant under deductive
consequence (\S\ref{sub-dep-deduct}), and so we expressly require
$\sigma^{\star}\wavy\tau^{\star}$ to be a member of $\Sigma$. It
is only under this strict definition of interpretability that
Theorem \ref{th-interp-lesseq} holds.)

\begin{theorem} \label{th-interp-lesseq} Let $A$ be a metric space,
and let\/ $\Sigma$ and $\Gamma$ be sets of equations (finite
or infinite) in any similarity types. If\/ $\Gamma$ is
interpretable in $\Sigma$, then
$\lambda_A(\Gamma)\,\leq\,\lambda_A(\Sigma)$.
\end{theorem} \begin{Proof}
If $\lambda_A(\Sigma)=\infty$, there is nothing more to prove, so
we shall assume that $\lambda_A(\Sigma)<\infty$. We select a real
$\varepsilon>0$, and keep it fixed until the last two sentences of
the proof. By definition of $\lambda_A$, there exists a
topological algebra $\mathbf A$, of type appropriate to $\Sigma$
and based on $A$, such that
\begin{align}                       \label{eq-Sigma-within-eps}
            \lambda_{\mathbf A}(\gamma\wavy\delta) \;< \;\lambda_A(\Sigma)
                        \,+\,\varepsilon
\end{align}
for any equation $\gamma\wavy\delta$ in $\Sigma$.

Let the operation symbols of $\Gamma$ be $F_t$ ($t\in T$). Since
$\Gamma\,\leq\,\Sigma$, there exist $\Sigma$-terms $\alpha_t$
($t\in T$) such that $\sigma^{\star}\wavy\tau^{\star}$ is in
$\Sigma$ for each equation $\sigma\wavy\tau$ of $\Gamma$. We now
define a topological algebra $\mathbf A'\,=\,\la A,\o F_t\ra_{t\in
T}$, {\em via}
\begin{align}                    \label{eq-def-deriv-op}
            F_t^{\mathbf A'}\EQ \alpha_t^{\mathbf A}
\end{align}
for each $t\in T$. The rest of the proof is based on the claim
that for any term $\sigma$ in the language of $\Gamma$, we have
\begin{align}               \label{claim-interp-terms}
              \sigma^{\mathbf A'} \EQ ({\sigma^{\star}})^{\mathbf
                                A}
\end{align}
(i.e.\ the $A'$-interpretation of $\sigma$ is the same operation
as the $A$-interpretation of $\sigma^{\star}$). The proof is by
induction on the length of $\sigma$. If $\sigma$ is a variable,
this assertion is immediate from (\ref{eq-star-variable}). If
$\sigma$ is a composite term, say
$\sigma=F_t(\tau_1,\ldots,\tau_{n(t)})$, we first rewrite
(\ref{eq-star-composite}) as
\begin{align}                  \label{eq-repeat-star}
         \sigma^{\star} \EQ
           \alpha_t(\tau_1^{\star},\ldots,\tau_{n(t)}^{\star}),
\end{align}
and then calculate
\begin{align*}
     (\sigma^*)^{\mathbf A}
                       (a_1,a_2,\ldots) &\EQ
               \alpha_t^{\mathbf A}((\tau_1^{\star})^{\mathbf A}
     (a_1,a_2,\ldots),\ldots,(\tau_{n(t)}^{\star})^{\mathbf
A}(a_1,a_2,\ldots))\\
           &\EQ F_t^{\mathbf A'}(\tau_1^{\mathbf
           A'}(a_1,a_2,\ldots),\dots,\tau_{n(t)}^{\mathbf
                 A'}(a_1,a_2,\ldots))\\
             &\EQ \sigma^{\mathbf A'}(a_1,a_2,\ldots).
\end{align*}
(The first equation uses (\ref{eq-repeat-star}) together with the
usual recursive definition of the term-function
$(\sigma^{\star})^{\mathbf A}$.  The second equation uses
(\ref{eq-def-deriv-op}); it also uses (\ref{claim-interp-terms})
inductively for $\tau_j^{\mathbf A'}=(\tau^{\star})^{\mathbf A}$
($1\leq j\leq n(t)$). The third line is based on the recursive
definition of $\sigma^{\mathbf A'}$.) This completes our proof of
(\ref{claim-interp-terms}).

From (\ref{claim-interp-terms}) it is immediate that if
$\sigma\wavy\tau$ is any equation in the operation symbols $F_t$,
then
\begin{align}                     \label{eq-same-lambdas}
         \lambda_{\mathbf A'}(\sigma\wavy\tau) \EQ
          \lambda_{\mathbf A}(\sigma^{\star}\wavy\tau^{\star}).
\end{align}
If, moreover, $\sigma\wavy\tau\in\Gamma$, then (by the first part
of this proof) $\sigma^{\star}\wavy\tau^{\star}\in\Sigma$. Thus
(\ref{eq-Sigma-within-eps}) and (\ref{eq-same-lambdas})
immediately yield
\begin{align}                          \label{eq-estimate-A-prime}
            \lambda_{\mathbf A'}(\sigma\wavy\tau) \;< \;\lambda_A(\Sigma)
                        \,+\,\varepsilon.
\end{align}
In other words, we have now shown that, for every positive real
$\varepsilon$, there is a topological algebra $\mathbf A'$ based
on $A$ such that (\ref{eq-estimate-A-prime}) holds for every
$\sigma\wavy\tau\in\Gamma$. Hence by \S\ref{sub-lambda-alg}, for
every positive real $\varepsilon$ there is a topological algebra
$\mathbf A'$ based on $A$ such that
\begin{align*}
           \lambda_{\mathbf A'}(\Gamma) \;\leq \;\lambda_A(\Sigma)
                        \,+\,\varepsilon.
\end{align*}
Since $\lambda_A$ (for a space $A$) is defined as the $\inf$ of
values of $\lambda_{\mathbf A}$ (for topological algebras based on
$A$---see \S\ref{sub-def-lamb-space}), we finally have
$\lambda_A(\Gamma)\leq\lambda_A(\Sigma)$.
\end{Proof}

Many applications of Theorem \ref{th-interp-lesseq} are fairly
obvious, and will not require any emphasis. For example, in
\S\ref{subsub-not-exp2} we prove that
$\lambda_{\Reals}(\Gamma_2)\geq \text{radius}(\Reals)/2$, for
$\Gamma_2$ a version of the theory of groups of exponent 2. It is
not hard to see (either directly or through Theorem
\ref{th-interp-lesseq}) that the same estimate holds for a (strong
enough; perhaps redundant) version of Boolean algebra. Many such
observations are possible as we go along; we generally will not
mention them. Nevertheless Theorem \ref{th-interp-lesseq} has
provided some guidance for our exposition: generally speaking we
have striven to attach estimates from below to theories that are
low in our quasi-ordering, and to attach estimates from above to
high theories.

We shall apply Theorem \ref{th-interp-lesseq} in
\S\ref{subsub-l-groups}, for (a certain set $\Lambda\Gamma$ of
equations of) the theory of lattice-ordered groups and $A$ an
arbitrary compact metric space. We could of course estimate
$\lambda_A(\Lambda\Gamma)$ directly, but it seems more informative
to take note of a certain theory $\Sigma_2\,\leq\,\Lambda\Gamma$,
for which an estimate of $\lambda_A$ seems more natural.
$\Sigma_2$ in effect represents an interesting Mal'tsev condition
satisfied by LO-groups, and the estimate seems to relate naturally
to this condition.

\subsection{Dependence on the choice of metric.}
                 \label{sub-dep-metric}
It is apparent from the definitions in \S\ref{sub-metric} and
\S\ref{sub-def-lamb-space} that $\lambda_A$ apparently depends on
our choice of metric to represent the topology of $A$, and hence
is not an invariant of the topological space $A$. Indeed numerous
examples will confirm that that this dependence is very real, even
among metrics that are normalized to be of diameter $1$. (See
e.g.\ \S\ref{subsub-depend-metric}, \S\ref{subsub-Y} and
\S\S\ref{subsub-square-no-inj}--\ref{sub-delta-square-embed-square}
below. In fact \S\ref{subsub-distort} contains an example where
the choice of (diameter-1) metric can yield any value strictly
between 0 and 1.)

Nevertheless, the metric-based quantity $\lambda_{(A,d)}$ can
convey interesting and useful information. We will return in
\S\ref{sec-top-inv} to some definitions (using $\inf$ and $\sup$)
that extract a topological invariant from the function
$d\GOESTO\lambda_{(A,d)}$. Here we shall only prove that, for
compact $A$, the condition $\lambda_A(\Sigma)>0$ is a topological
invariant. (Compactness is essential, as one may see from the two
$\lambda$-values that are calculated in
\S\ref{subsub-depend-metric} below.)

\begin{lemma}  \label{lemma-comp-two-metrics}
If\/ $\rho$ and $d$ are metrics defining one and the
same topology on $A$, if that topology is compact, and if\/
$\lambda_{(A,\rho)}(\Sigma)>0$, then $\lambda_{(A,d)}(\Sigma)>0$.

\end{lemma} \begin{Proof}
The compact space $A$ has finite diameter, and clearly both
$\lambda$-values are limited by the diameter
(\S\ref{subsub-lambda}), and hence both are finite. Suppose that
$\lambda_{(A,\rho)}(\Sigma)=\varepsilon_2>0$, and let $d$ be a
metric for the same topology on $A$. Take $\varepsilon_1>0$ to be
a Lebesgue number of the metric space $(A,d)$ for its covering by
all $\frac{\varepsilon_2}{2}$-balls of the metric space
$(A,\rho)$. (In other words, for each $a\in A$ there exists $a'\in
A$ such that $B_d(a,\varepsilon_1)\subseteq
B_{\rho}(a',\varepsilon_2/2)$.) We immediately have\footnote%
{It may already be in the literature that for all
$\varepsilon_2>0$ there exists $\varepsilon_1>0$ satisfying
(\ref{eq-comp-metrics}); we have not seen it.}
\begin{align}                       \label{eq-comp-metrics}
   \mbox{if}\quad d(a,b)\leq\varepsilon_1 \quad\mbox{then}\quad
   \rho(a,b)\leq\varepsilon_2.
\end{align}
To complete the proof, it will be enough to establish that
\begin{equation*}
\lambda_{(A,d)}(\Sigma)<\varepsilon_1
\end{equation*}
is impossible. If it held, then we would have $((A,d),\o
F_t)_{t\in T}\models_{\varepsilon_1}\Sigma$ for some continuous
operations $\o F_t$. By (\ref{eq-comp-metrics}), these same
operations $\o F_t$ would establish that
$\lambda_{(A,\rho)}(\Sigma)<\varepsilon_2$, which contradicts our
original assumption about $\varepsilon_2$. This contradiction
completes the proof of the Lemma.
\end{Proof}

\subsection{The usual $\Reals$-metric yields only extreme
values for $\lambda_{\Reals}$.}
\begin{theorem}\label{th-shrink-metric} Consider a metric
space $(A,d)$ for which there exist functions $\phi$ and $\gamma$
with the following properties: $\phi\FROM[0,\infty)\TO[0,\infty)$
is a monotone increasing function such that, for all $x$, the
sequence $\phi^n(x)$ approaches $0$ in $\Reals$. And suppose that
$\gamma\FROM A\TO A$ is a homeomorphism satisfying
\begin{align}             \label{eq-shrink-rate}
             d(\gamma(a),\gamma(b)) \;\leq\; \phi(d(a,b))
\end{align}
for all $a,b\in A$. Then for any $\Sigma$, $\lambda_{(A,d)}(\Sigma)
\,=\, 0$ or $\infty$.
\end{theorem} \begin{Proof}
Suppose $\lambda_{(A,d)}(\Sigma) \,<\,\infty$. Then there exist
$K<\infty$ and continuous operations $\o F_t$ on $A$ such that
each equation of $\Sigma$ holds for these operations within $K$.
If $K=0$, then clearly we are done. Otherwise, consider the
operations ${\o F_t}'$, where
\[
      {\o F_t}'(a_1,a_2, \cdots) \EQ
                   \gamma\,(\o F_t(\INV{\gamma}(a_1),
                   \INV{\gamma}(a_2), \cdots)).
\]
We claim that these operations satisfy $\Sigma$ within $\phi(K)$.

For any term $\sigma$, let $\o\sigma$ denote the usual term function
associated to $\sigma$ in the algebra $(A,\o F_t)_{t\in T}$, and let
${\o\sigma\,}'$ denote the corresponding term function in $(A,{\o
F_t}')_{t\in T}$. A straightforward inductive argument establishes
that
\[
             {\o \sigma\,}'(a_1,a_2, \cdots) \EQ
                   \gamma\,(\o \sigma(\INV{\gamma}(a_1),
                   \INV{\gamma}(a_2), \cdots)).
\]
Now, for any equation $\sigma\wavy\tau$ of $\Sigma$, we may
calculate
\begin{align*}
      d({\o \sigma\,}'(a_1, \cdots),{\o \tau\,}'(a_1,
      \cdots)) &\EQ d(\gamma\,(\o \sigma(\INV{\gamma}(a_1),
                   \cdots)),\gamma\,(\o \tau(\INV{\gamma}(a_1),
                   \cdots)))\\
        &\;\leq\; \phi\,(d(\o \sigma(\INV{\gamma}(a_1),
                   \cdots),\o \tau(\INV{\gamma}(a_1),
                   \cdots)))\\
        &\;\leq\; \phi(K),
\end{align*}
where the second line follows from (\ref{eq-shrink-rate}). Thus our
claim is verified.

Now this construction may clearly be iterated: there are continuous
operations satisfying $\Sigma$ within $\phi^2(K)$, within
$\phi^3(K)$, and so on. Since $\phi^n(K)\rightarrow0$, it is clear
that $\lambda_{(A,d)}(\Sigma) \,=\,0.$
\end{Proof}

Our main example of a metric space satisfying the hypotheses of
Theorem \ref{th-shrink-metric} is $\Reals^n$ with the usual
metric. For this space, $\gamma$ and $\phi$ can be taken both to
be scalar multiplication by $1/2$. Thus
\begin{corollary}           \label{cor-zero-infty}
If $d$ is the usual metric on $\Reals^n$, and $\Sigma$ is any set
of equations, then $\lambda_{(\Reals^n,d)}(\Sigma) \,=\, 0$ or
$\infty$.\ENDPROOF
\end{corollary}

On the other hand, if $d$ is the usual metric on $(0,1)$ (an
obvious homeomorph of $\Reals$), then such $\phi$ and $\gamma$ do
not exist. (As may be seen from a simple measure-theoretic
argument.)\vspace{0.1in}

\hspace*{-\parindent}%
{\bf Problem} --- is Corollary \ref{cor-zero-infty} true on
$(0,1)$?

\subsection{An inequality for retractions.}  \label{sub-inequal-retract}
Theorem \ref{thm-inequal-retract}, which follows, is typically
applied in a context where $\lambda_{(B,d)}$ is known (either
exactly or approximately), and an appropriate $K$-value is known
or can be found. Then the inequality (\ref{eq-inequal-retract}) is
used to supply an upper estimate on $\lambda_{(A,d)}$. Such use of
the theorem may be found in \S\ref{sub-sub-S2-epsilon},
\S\ref{subsub-distort}, \S\ref{subsub-S1-different},
\S\ref{subsub-Y-new-metric}, \S\ref{subsub-approx-cubes} and
\S\ref{sub-circle-revisit}.

\begin{theorem}               \label{thm-inequal-retract}
Let\/ $(A,d)$ be a metric space, and $B$ a subset
of\/ $A$. Suppose that there is a continuous map $\psi\FROM A\TO
B$ such that $\psi\upharpoonright B$ is the identity on $B$ (i.e.,
$\psi$ retracts $A$ onto $B$), and such that\/ $d(a,\psi(a))\leq
K$ for all\/ $a\in A$. Then
\begin{align}                 \label{eq-inequal-retract}
             \lambda_{(A,d)}(\Sigma)\;\leq\;\lambda_{(B,d)}(\Sigma)
                     \,+\,K
\end{align}
for every consistent $\Sigma$.
\end{theorem} \begin{Proof}
Let us assume that the operations appearing in $\Sigma$ are $F_t$
($t\in T$). For each $\varepsilon>0$, there are continuous
operations $F_t^{\mathbf B}$ on $B$, forming an algebra $\mathbf
B$ such that $\lambda_{\mathbf B}(\Sigma)<\lambda_{(B,d)}(\Sigma)
+ \varepsilon$. We define operations $F_t^{\mathbf A}$ on $A$ as
follows:
\begin{align}             \label{eq-define-approx-A}
            F_t^{\mathbf A}(a_1,\cdots,a_{n(t)}) \EQ
                     F_t^{\mathbf B}(\psi(a_1),\cdots,\psi(a_{n(t)})).
\end{align}
Since $B$ is closed under the operations $F_t^{\mathbf B}$, and
since $\psi$ retracts $A$ onto $B$, we readily see that
\begin{align*}
            \sigma^{\mathbf A}(a_1,\cdots) \EQ
                     \sigma^{\mathbf B}(\psi(a_1),\cdots))
\end{align*}
for any term $\sigma$ in our language (other than $\sigma$ a
single variable).

To estimate $\lambda_{\mathbf A}(\Sigma)$, we consider a single
equation $\sigma\wavy\tau$ of $\Sigma$. If $\sigma$ and $\tau$ are
both variables, they must be the same variable, since $\Sigma$ is
consistent. In this case $\lambda_A$ and  $\lambda_B$ both take
the value $0$, and so (\ref{eq-inequal-retract}) holds, for this
one equation.

If neither $\sigma$ nor $\tau$ is a variable, then for any
$a_1,a_2,\cdots\in A$, we have
\begin{align*}
           d(\sigma^{\mathbf A}(a_1,\cdots),\tau^{\mathbf
           A}(a_1,\cdots)) &\EQ d(\sigma^{\mathbf
           B}(\psi(a_1),\cdots)),\tau^{\mathbf
           B}(\psi(a_1),\cdots))\\
           \;&<\;\lambda_{(B,d)}(\Sigma)\,+\,\varepsilon.
\end{align*}
If, say, $\tau$ is the variable $x_1$, then we have
\begin{align*}
           d(\sigma^{\mathbf A}(a_1,\cdots),\tau^{\mathbf
           A}(a_1,\cdots)) &\EQ d(\sigma^{\mathbf
           B}(\psi(a_1),\cdots)),a_1)\\
              &\leq \;d(\sigma^{\mathbf
              B}(\psi(a_1),\cdots)),\psi(a_1)) \,+\,
                       d(\psi(a_1),a_1)\\
           \;&<\;\lambda_{(B,d)}(\Sigma)\,+\,\varepsilon \;+\; K.
\end{align*}
It is now clear that for all $\sigma\wavy\tau\in\Sigma$ we have
\begin{align}        \label{eq-estimate-A}
            \lambda_{\mathbf A}(\sigma,\tau) \leq
                  \lambda_{(B,d)}(\Sigma)\,+\,\varepsilon
                  \,+\, K,
\end{align}
which is to say
\begin{align*}
            \lambda_{\mathbf A}(\Sigma) \leq
                  \lambda_{(B,d)}(\Sigma)\,+\,\varepsilon
                  \,+\, K.
\end{align*}
Since this estimate can be made true for every $\varepsilon>0$,
and since $\lambda_{(A,d)}(\Sigma)$ is the $\inf$ of all possible
$ \lambda_{\mathbf A}(\Sigma)$ values, we clearly have established
Equation (\ref{eq-inequal-retract}).
\end{Proof}

Let $A$ be an interval $[a,b]$, and let $B$ be the subinterval
$[a+K, b-K]$. Map $\psi\FROM A \TO B$ via
\[
 \psi(x) \EQ (a+K)\vee[x\wedge(b-K)].
\]
Then $d(a,\psi(a))\leq K$ for all $a\in A$, as was assumed for
Theorem \ref{thm-inequal-retract}. Equation
(\ref{eq-inequal-retract}) fails, however, if we take an
inconsistent $\Sigma$ consisting of the single equation $x_0\wavy
x_1$. In this case $\lambda_A$ and $\lambda_B$ evaluate to the
diameters of the respective intervals.

\section{Some estimates of $\lambda_A(\Sigma)$.}\label{sub-num-ex}
As we remarked at the end of \S\ref{sub-def-lamb-space}, it is
easy to compute $\lambda_A(\Sigma)$ when $A\models\Sigma$. In most
other cases it seems to be difficult, but we are able to get some
approximate results. We are able to compute exact values of
$\lambda_A(\Sigma)$  in \S\ref{subsub-lim-zero},
\S\ref{subsub-dep-ded},
\S\ref{sub-sub-Sn}, \S\ref{sub-sub-S1}, \S\ref{sub-sub-S3},
\S\ref{sub-sub-S7}, 
\S\ref{subsub-distort} and \S\ref{subsub-lambda-homotopy}.

The topics of \S\ref{sub-num-ex} are gathered mostly according to
a rough classification of the $\Sigma$ appearing, which is to say,
according to the rows of the tables in Figures \ref{fig-chart-one}
and \ref{fig-chart-two} of \S\ref{sub-def-lamb-space}.
(\S\ref{sub-spheres-simple} refers to spheres, but is then itself
subdivided according to $\Sigma$.) It would of course be possible
to arrange the material according to the space $A$ involved
(columns of the tables), or according to method of proof. Perhaps
the following remarks will be of some use to the reader who wishes
to see some organization according to method.

For proving lower estimates, i.e.\ estimates of the form
$\lambda_A(\Sigma)\geq K$, we have essentially two methods. These
are, roughly
\begin{itemize}
\item[(i)] Homotopy. In some contexts, such as $S^n$, two maps near
           to each other must be homotopic. If we already know that
         $\Sigma$ cannot be satisfied on $A$ within homotopy, this
          knowledge may be applied to yield an estimate on $\lambda_A
          (\Sigma)$. This
        method may be found throughout \S\ref{sub-spheres-simple},
         and appears again in \S\ref{subsub-group-cyl}.
\item[(ii)] We start with some $A$-points $P$, $Q$, etc., that are a known
        distance apart, $d(P,Q)\geq d\ldots$ etc.
        Some auxiliary points $E,F,\ldots$ are then located, using
        both the operations and various topological properties, such as
                \begin{itemize}
         \item[(a)] The fixed-point theorem of Brouwer et al. (For an
                   example, see \S\ref{subsub=group-f.p.}.)
         \item[(b)] Certain homotopy classes of maps must be onto. (For an
                   example, see \S\ref{subsub-distort}.)
         \item[(c)] Compactness; in particular, the theorem (of Bolzano
            and Weierstrass) that
            every sequence in a compact metric space has a convergent
              susequence. (This is applied in
               \S\ref{subsub-aux-to-log}, which has a corollary in
                   \S\ref{subsub-l-groups}---an estimate
                        for lattice-ordered groups.)
      \item[(d)] The Intermediate Value Theorem. (For examples,
                   see \S\ref{subsub-not-exp2},
           \S\ref{subsub-depend-metric},
                \S\ref{subsub-inf-dim}. (Lemma \ref{lem-special-acyclic}
                 of \S\ref{subsub-approx-lattice} below supplies a
                  slight generalization of the IVT, which is
                   applied to lattice theory in Lemma \ref{lem-approx-lat}
                 of \S\ref{subsub-approx-lattice}, and then to an
                 estimate of $\lambda$ in \S\ref{subsub-Y}.)
        \item[(e)] The Theorem of Borsuk and Ulam. (For examples,
                   see \S\ref{subsub-square-no-inj},
                 \S\ref{sub-genzations-delta},
                 \S\ref{subsub-approx-cubes}.)
        \end{itemize}
         Then the approximate satisfaction of $\Sigma$ within $K$
         would force some of the
         points $P$, $Q$, $E$, $F$, \ldots, to be close to one another,
         and the triangle inequality would yield e.g. $d(P,Q)<d$. This
         contradiction
        establishes that $\lambda_A(\Sigma)>K$.
\end{itemize}

For proving an upper estimate on $\lambda_A(\Sigma)$, we generally
employ a direct, constructive method. We propose a topological
algebra $\mathbf A$ based on $A$, and then directly estimate the
supremum appearing in (\ref{def-lambda-A}) from some knowledge of
the geometry of $A$ and some knowledge of the term-operations
$\sigma^{\mathbf A}$ and $\tau^{\mathbf A}$. For some direct
applications of this method, see \S\ref{subsub-depend-metric},
\S\ref{subsub-lim-zero-consis}, \S\ref{subsub-inf-dim},
\S\ref{subsub-square-no-inj} and
\S\ref{sub-delta-square-embed-square}. We sometimes also apply
this method indirectly through the use of Theorem
\ref{thm-inequal-retract}: in our definition of the topological
algebra $\mathbf A$ is accomplished by Equations
(\ref{eq-define-approx-A}), and the required estimates are
accomplished by Equations (\ref{eq-estimate-A}). As mentioned
earlier, such applications of Theorem \ref{thm-inequal-retract}
can be seen in \S\ref{sub-sub-S2-epsilon}, \S\ref{subsub-distort},
\S\ref{subsub-S1-different}, \S\ref{subsub-Y-new-metric},
\S\ref{subsub-approx-cubes} and \S\ref{sub-circle-revisit}.

In anticipation of \S\ref{sec-top-inv}, at some places in
\S\ref{sub-num-ex}---see e.g. \S\ref{sub-sub-S2-epsilon},
\S\ref{sub-sub-infty}, \S\ref{subsub-S1-different},
\S\ref{subsub-Y-new-metric}, \S\ref{sub-delta-square-embed-square}
and \S\ref{subsub-approx-cubes} (Theorem
\ref{th-nthpow-on-inttok})---we will introduce an unfamiliar
metric $d$ for a topological space $A$ and then estimate
$\lambda_{(A,d)}(\Sigma)$ from above, by one of the methods
mentioned here.

Throughout \S\ref{sub-num-ex}, we allow the use of $a\wavy b$, for
$a,b$ in a metric space, to mean that $d(a,b)<\varepsilon$ (where
$\varepsilon$ will be understood in context). We have thus
overloaded the symbol ``$\wavy,$'' using it to denote both formal
equations and approximate equality. (The former usage applies only
to terms; the latter applies only to elements of metric spaces.)

Although a lot of detailed information arises in
\S\ref{sub-num-ex}, we remind the reader to consult the charts in
Figures \ref{fig-chart-one} and \ref{fig-chart-two} of
\S\ref{sub-def-lamb-space} for an overview.

 \subsection{Inconsistent $\Sigma$} \label{subsub-lim-zero}
 \subsubsection{$\Sigma$ contains $x_i\wavy x_j$ for $i\neq j$.}
               \label{subsub-wavy}
 As we saw in \S\ref{sub-dep-deduct}, in this case $\lambda_A(\Sigma) =
 \text{diam}(A)$.
 \subsubsection{$\Sigma$ contains no $x_i\wavy x_j$ for $i\neq j$.}
                     \label{sub-sub-no-wavy} \label{subsub-dep-ded}
As we saw in \S\ref{sub-dep-deduct}, in this case
$\lambda_A(\Sigma) \leq \text{radius}(A)$. Here in
\S\ref{sub-sub-no-wavy} we give an example of $\Sigma$ and
$\Sigma'$ of this type, with $\Sigma'\subseteq\Sigma^{\star}$, and
where
\[
\lambda_{A}(\Sigma)\EQ0,\;\;\;\text{and}\;\;\;
\lambda_{A}(\Sigma')\EQ\text{diam}(A)/2\EQ\text{radius}(A).
\]

Take $(A,d)$ to be a simplex of finite dimension, where, for now,
$d$ is any metric that defines the usual topology. We will think
of $A$ as equal to $[0,1]^k$ for some $k$. If $K$ is a closed
subset of $A^n$ for some $n$, and $\o F\FROM K\TO A$ is a
continuous map, then $\o F$ may be extended to a continuous map
$\o F\FROM A^n\TO A$. (To see this, apply Tietze's Theorem to each
coordinate map $\pi_j\Compos\o F\FROM K\TO[0,1]$, where $\pi_j$ is
the $j$-th co-ordinate projection $[0,1]^k\TO[0,1]$.)

For a single ternary operation symbol $F$ and two constants $a$
and $b$, we take $\Sigma$ to comprise the equations
\begin{align*}
        F(a,x_0,x_1)  &\WAVY  x_0  \\
        F(b,x_0,x_1)  &\WAVY  x_1  \\
                a &\WAVY b.
\end{align*}
This $\Sigma$ is obviously inconsistent, hence not compatible with
$A$. On the other hand, we shall easily see that
$\lambda_A(\Sigma)=0$. Moreover, for the sake of future reference
(\S\ref{subsub-delta+0}), we point out that our calculation for
$\Sigma$ applies to any metric that defines the given topology on
$A$.

It will suffice, given $\varepsilon>0$, to define $\o F$, $\o a$
and $\o b$ on $A$, so that $\lambda_{\mathbf
A}(\Sigma)<\varepsilon$ for $\mathbf A = (A;\o F,\o a,\o b)$.
Since $A$ is not discrete, there exist $\o a, \o b\in A$ with
$0<d(\o a,\o b)<\varepsilon$. The third equation on $\Sigma$ now
holds within $\varepsilon$; to complete the job we shall define a
continuous ternary operation $\o F$ for which the first two
equations hold {\em exactly}. In this way we shall have all three
equations holding within $\varepsilon$, as desired.

The first two equations amount to a definition of $\o F$ on a
closed subset of $A^3$. By the corollary of Tietze's Theorem that
is described above, $\o F$ may be extended to a continuous map
with domain $A^3$, and the construction is complete. This
completes our derivation of $\lambda_A(\Sigma)=0.$

Let us now enlarge $\Sigma$ as follows:
\begin{align*}
      \Sigma' \EQ \Sigma\,\cup\,\{F(a,x_0,x_1)  &\WAVY  x_1\}.
\end{align*}
Notice that $\Sigma'$ differs from $\Sigma$ only in including one
logical consequence of $\Sigma$ (derived from $\Sigma$ by
substitution). Nevertheless,  as was mentioned in
\S\ref{sub-dep-deduct}, the value of $\lambda_A$ may increase
under such an enlargement of $\Sigma$. In this case, as will be
apparent from Lemma \ref{lemma-comp-two-metrics},
$\lambda_A(\Sigma')>0$ for any appropriate choice of the metric
$d$.

For a precise calculation of $\lambda_A(\Sigma')$, we will now
take $d$ to be the usual Euclidean metric on $A^k$. Under this
assumption, we shall now prove that
$\lambda_A(\Sigma')=\text{diam}(A)/2$.

By definition, there exist a constant $\o a$ and a ternary
operation $\o F$ on $A$ so that
\begin{align*}
                    d(\o F(\o a, a_1, a_2), a_1) \,&\leq\,
                    \lambda_A(\Sigma') \\
                    d(\o F(\o a,a_1,a_2 ),a_2) \,&\leq\,
                     \lambda_A(\Sigma'),
\end{align*}
for all $a_1,a_2\in A$, which implies $d(a_1,a_2) \,\leq\,
2\lambda_A(\Sigma')$ or $\lambda_A(\Sigma') \,\geq\,
\text{diam}(A)/2$. On the other hand if we define $\o
F(a_1,a_2,a_3)$ to be the midpoint of the segment joining $a_2$
and $a_3$, then the resulting $\o F$ is clearly continuous, and
the equations of $\Sigma'$ are clearly satisfied within
$\text{diam}(A)/2$; hence $\lambda_A(\Sigma') \,\leq\,
\text{diam}(A)/2$.
%


\subsection{Some exact formulas for spheres.}
              \label{sub-spheres-simple}

An algebra {$\mathbf A$} is {\em trite} iff every operation of
$\mathbf A$ is either a projection map or constant. An equational
theory $\Sigma$ is {\em easily satisfied\/} or {\em undemanding\/}
iff it has a trite model of more than one element---and hence has
trite models of every cardinality. All other theories are {\em
demanding.} (These definitions come from \cite{wtaylor-sae}.) It
is easily seen that if $\Sigma$ is undemanding and $A$ is any
space, then $A\models\Sigma$, and so $\lambda_A(\Sigma)=0$. Hence
every non-trivial calculation or estimation of $\lambda_A(\Sigma)$
in this paper will concern a demanding theory $\Sigma$.

Obviously $\Sigma$ and $\Sigma^{\star}$ have the same models, and
hence $\Sigma$ is demanding if and only if $\Sigma^{\star}$ is
demanding. Thus any evaluation of $\lambda_A(\Sigma)$ that is
based on the demanding property, such as (\ref{eq-simple-lambda})
below, holds equally well for $\lambda_A(\Sigma^{\star})$.

We showed in \cite{wtaylor-sae} that {\em certain spaces $A$ are
compatible only with the undemanding theories.}\footnote{%
In  \cite{wtaylor-sae} we trace one case of this back to work of
J. F. Adams \cite{adams}. In \cite{wallace} A. D. Wallace
attributes it to E. Cartan in the
special case of $\Sigma=$ groups.} %
(These include the spaces $A$ mentioned in \S\ref{sub-sub-Sn} and
\S\ref{sub-sub-infty} below.) In other words, if $\Sigma$ is
demanding then $A\not\models\Sigma$. In \S\ref{sub-sub-Sn}
we will have the stronger conclusion that $\lambda_A(\Sigma)>0$
for certain such $A$ and $\Sigma$. Indeed, in certain cases of
interest for \S\ref{sub-spheres-simple} the function $\lambda_{A}$
takes a
particularly simple form\footnote%
{By the estimates found in \S\ref{subsub-lambda}, this $\lambda_A$
is as large as possible among all metric spaces $A$---such as
spheres---that have radius and diameter both equal to $1$.}, %
 namely
\begin{align}              \label{eq-simple-lambda}
       \lambda_{A}(\Sigma) \EQ
            \begin{cases}
                  0 & \text{if $\Sigma$ is undemanding}\\
                  1 & \text{if $\Sigma$ is demanding}.
            \end{cases}
\end{align}

As we shall see in \S\ref{sub-sub-Sn}, Equation
(\ref{eq-simple-lambda}) holds for $A=S^n$, a diameter-$1$ sphere
of dimension $n$, where $n\neq 1,3, 7$. Moreover the proof of
(\ref{eq-simple-lambda}) will be almost immediate from the main
theorem of \cite{wtaylor-sae} (which is based on serious results
of algebraic topology). Notice also that Equation
(\ref{eq-simple-lambda}) implies invariance under deduction (since
``demanding'' depends only on the class of models of $\Sigma$) --
and hence the difficulties of \S\ref{sub-dep-deduct} do not arise.
The same method will also yield some values of $\lambda_A(\Sigma)$
for $A=S^1$ (see \S\ref{sub-sub-S1}), for $A=S^3$ (see
\S\ref{sub-sub-S3}) for $A=S^7$ (see \S\ref{sub-sub-S7}), and for
various $\Sigma$.

\subsubsection{$\Sigma$ demanding,  $A=S^n$
($n\neq1,3,7$).}
                           \label{sub-sub-Sn}
$S^n$ is the usual $n$-dimensional sphere, which may be explicitly
realized as the set of points in $\Reals^{n+1}$ having Euclidean
distance $1/2$ from the origin. Equation (\ref{eq-simple-lambda})
actually holds for any diameter-1 metric $d$ that has the
following property: the diameter of $(A,d)$ is realized by each
pair of antipodal points, and by no other pairs of points. Scaled
geodesic distance has this property, as does, more simply, scaled
Euclidean distance as measured in the ambient space
$\Reals^{n+1}$. (Given such a metric, we re-scale it so that the
diameter is $1$.)

The first alternative of (\ref{eq-simple-lambda}) is immediate: if
$\Sigma$ is undemanding, then $\lambda_{S^n}(\Sigma)=0$.

For the second alternative, we prove the contrapositive: given
$\Sigma$ with $\lambda_{S^n}(\Sigma)\neq 1$, we shall prove that
$\Sigma$ is undemanding. It is immediate from
\S\ref{subsub-lambda} that $\lambda_{S^n}(\Sigma)<1$.  By the
definitions in \S\ref{sub-lambda-alg}, there exist a constant
$K<1$ and a topological algebra $\mathbf A$ based on $A=S^n$, such
that $\lambda_{\mathbf A}(\sigma,\tau)\leq K$ for each equation
$\sigma\wavy\tau$ of $\Sigma$.
According to Equation (\ref{def-lambda-A}), for each
$\sigma\wavy\tau$ in $\Sigma$, and for each $\mathbf
a\in(S^n)^{\omega}$, we have
\begin{align}                                \label{eq-sn-close}
          d(\sigma^{\mathbf A}(\mathbf a),\tau^{\mathbf A}(\mathbf
          a))      <  K.
\end{align}

We now define a homotopy
\begin{align*}
          \Gamma\FROM(S^n)^{\omega}\times[0,1]\TO S^n
\end{align*}
between $\sigma^{\mathbf A}$ and  $\tau^{\mathbf A}$, as follows.
Let us be given $\mathbf a\in (S^n)^{\omega}$ and $t_0\in[0,1]$.
By Equation (\ref{eq-sn-close}), $\sigma^{\mathbf A}(\mathbf a)$
and $\tau^{\mathbf A}(\mathbf a)$ are not antipodes; hence there
is a unique shortest geodesic $\gamma(t)$, whose parameter $t$ is
proportional to arc length along $\gamma$, and with
$\gamma(0)=\sigma^{\mathbf A}(\mathbf a)$ and
$\gamma(1)=\tau^{\mathbf A}(\mathbf a)$. We then define
\[
              \Gamma(\mathbf a,t_0)\EQ \gamma(t_0).
\]
Clearly $\Gamma$ is continuous and is a homotopy between
$\sigma^{\mathbf A}$ and $\tau^{\mathbf A}$.

It is now clear that $\mathbf A$ satisfies $\Sigma$ up to
homotopy; in other words, $S^n$ is compatible with $\Sigma$ up to
homotopy. By Theorem 1 of \cite{wtaylor-sae}, $\Sigma$ is
undemanding. This concludes our proof of the contrapositive of
$\lambda_{S^n}(\text{demanding})=1$.

\subsubsection{$\lambda_{S^2}(\Sigma)$ in a different metric, for $\Sigma\,=$
                    H-spaces.}
                 \label{sub-sub-S2-epsilon}
By Lemma \ref{lemma-comp-two-metrics} of \S\ref{sub-dep-metric},
and by \S\ref{sub-sub-Sn} just above, we clearly have
\begin{align}              \label{eq-simple-lambda-b}
       \lambda_{(S^n,d)}(\Sigma)\;
            \begin{cases}
                \;=\;  0 & \text{if $\Sigma$ is undemanding}\\
                \;>\;  0 & \text{if $\Sigma$ is demanding},
            \end{cases}
\end{align}
for {\em any metric} $d$ on the sphere $S^n$ ($n\neq 1,3,7$). Here
we will see that the non-zero value associated to a demanding
theory can, in some cases, be arbitrarily small, even when the
diameter is constrained to be $1$. (Similar remarks are made for
$S^1$ in \S\ref{subsub-S1-different}.)

 In \S\ref{sub-sub-S2-epsilon} we take
$\Sigma$ to be the theory of H-spaces, otherwise known as the
theory of a two-sided unit element. It consists of the two
equations
\begin{align}            \label{eq-H-space}
            F(e,x) \Wavy x,\quad\quad F(x,e)\Wavy x,
\end{align}
for a binary operation $F$ and a nullary operation $e$. One easily
checks that this $\Sigma$ is demanding.

Given real $\varepsilon>0$, we shall place a non-standard metric
$d$ (of diameter 1) on the ordinary sphere $S^2$ and calculate
that $\lambda_{(S^2,d)}(\Sigma^{\star})\leq\varepsilon$. To define
the metric $d$, we consider the realization of $S^2$ as the
prolate ellipsoid that is the locus of
\begin{align*}
      4\varepsilon^2 x^2 \,+\, 4y^2 \,+\, 4 z^2 \EQ
                        \varepsilon^2
\end{align*}
in ordinary 3-space. Then $d$ is defined to be the ordinary
Euclidean distance of $\Reals^3$, restricted to the ellipsoid. One
easily checks that $(S^2,d)$ has diameter $1$.

We define a continuous function $\psi \FROM S^2\TO S^2$ via
\begin{align*}
         \psi(x,y,z) \EQ (x,\sqrt{y^2+z^2},0).
\end{align*}
One easily checks that, for all $\mathbf x= (x,y,z)\in S^2$, the points
$\mathbf x$ and $\psi(\mathbf x)$ lie on a circle of radius $\leq
\varepsilon/2$ (in a plane perpendicular to the $x$-axis); hence
\begin{align}                \label{eq-fors2-a}
           d(\mathbf x,\psi(\mathbf x))\;\leq\;\varepsilon.
\end{align}
The image of $\psi$ is the semi-ellipse $E$ in the $x,y$-plane
that is bijectively parametrized by
\begin{align*}
         e(t) \EQ(\frac{1}{2}\,\sin \frac{\pi t}{2},\,
         \frac{\varepsilon}{2}\,\cos \frac{\pi t}{2},\,0)
\end{align*}
for $-1\leq t\leq 1$. One easily sees that
\begin{align}                     \label{eq-fors2-b}
                    \psi\upharpoonright E \EQ\text{identity}.
\end{align}
The reader may check that we now have the hypotheses of Theorem
\ref{thm-inequal-retract} of \S\ref{sub-inequal-retract} (with
$S^2$ for $A$, $E$ for $B$, and $\varepsilon$ for $K$). In this
context, the conclusion of the theorem is
\begin{align*}
        \lambda_{(S^2,d)}(\Sigma^{\star})
        \;\leq\;\lambda_{(E,d)}(\Sigma^{\star})
                     \,+\,\varepsilon.
\end{align*}
Moreover $E$ is topologically a closed segment, which can be made
into an H-space (in several interesting ways); hence
$\lambda_{(E,d)}(\Sigma^{\star})=0$. We now have the desired
conclusion that $\lambda_{(S^2,d)}(\Sigma^{\star})<\varepsilon$.

\subsubsection{$\Sigma$ not Abelian,   $A=S^1$.}
                           \label{sub-sub-S1}
We call a set $\Sigma$ of equations {\em Abelian} iff it is
interpretable (in the sense of \cite{ogwt-mem}) in the equational
theory of Abelian groups (or $\Integers$-modules). Equivalently,
$\Sigma$ is Abelian if and only if it has a model based on
$\Integers$ with operations of the form
\begin{align}          \label{eq-def-op-abelian}
   \o F(x_1,\cdots,x_n)\EQ m_1x_1+\cdots+m_nx_n,
\end{align}
where each $m_i\in\Integers$.

It was proved in Theorem 41 on page 234 of \cite{wtaylor-sae} that
if $\Sigma$ is compatible with
$S^1$, even if compatible only within homotopy, then
$\Sigma$ is Abelian, and conversely if $\Sigma$ is Abelian, then
$\Sigma$ is compatible with $S^1$. Therefore,
$\lambda_{S^1}(\Sigma)=0$ for every Abelian
$\Sigma$. In fact, we have the stronger result that
\begin{align}              \label{eq-simple-lambda-1}
       \lambda_{S^1}(\Sigma) \EQ
            \begin{cases}
                  0 & \text{if $\Sigma$ is Abelian}\\
                  1 & \text{otherwise}.
            \end{cases}
\end{align}
We sketch the proof. We have already established the first
alternative of (\ref{eq-simple-lambda-1}). To prove the second
alternative through its contrapositive, we begin with
$\lambda_{S^1}(\Sigma)\neq1$. Proceeding as in \S\ref{sub-sub-Sn},
we get $\sigma^{\mathbf A}$ homotopic to $\tau^{\mathbf A}$ for
each equation $\sigma\wavy\tau$ of $\Sigma$. It is immediate from
Theorem 41 on page 234 of \cite{wtaylor-sae} that $\Sigma$ is
Abelian. This completes the contrapositive proof of the second
alternative.

If $\Sigma$ is Abelian, then $\Sigma$ is modeled by an algebra
$\mathbf Z\ISO (\Integers, \cdots)$ of the aforementioned type;
clearly $\mathbf Z$ also models $\Sigma^{\star}$; hence
$\Sigma^{\star}$ is also Abelian. Thus (\ref{eq-simple-lambda-1})
yields the same values for $\lambda_{S^2}(\Sigma^{\star})$ and for
$\lambda_{S^2}(\Sigma)$.

For some specific instances of non-Abelian equational theories,
the reader is referred to \S\ref{subsub-lambda-homotopy}.

\subsubsection{$\Sigma=\text{Abelian groups}$, $A=S^3$.}
                \label{sub-sub-S3}
$S^3$ differs from the general $S^n$ of \S\ref{sub-sub-Sn}, and
from $S^7$, in that $S^3$ is compatible with $\Gamma$, the
equational theory of groups, and hence $\lambda_{S^3}(\Gamma)=0$.
In \S\ref{sub-sub-S3} we shall also prove that
$\lambda_{S^3}(\Gamma_0)=1$, where $\Gamma_0$ stands for any set
of equations axiomatizing Abelian groups.

The proof follows that of \S\ref{sub-sub-Sn} and
\S\ref{sub-sub-S1}, and proceeds by contradiction. If
$\lambda_{S^3}(\Gamma_0)<1$, then there is a topological algebra
$\mathbf A$ on $S^3$ that substantiates $\lambda<1$. As in the
earlier proofs, we obtain $\sigma^{\mathbf A}$ homotopic to
$\tau^{\mathbf A}$ for each equation $\sigma\wavy\tau$ of
$\Gamma_0$. In other words $\mathbf A$ is an Abelian group on
$S^3$, up to homotopy. A 1953 theorem of R. Bott \cite{bott} says
that no such $\mathbf A$ exists.

\subsubsection{$\Sigma=\text{Monoids}$, $A=S^7$.}
                \label{sub-sub-S7}

$S^7$ differs from the general $S^n$ of \S\ref{sub-sub-Sn} in that
$S^7$ is compatible with $K$, the equational theory of H-spaces,
as defined in Equation (\ref{eq-H-space}). (Multiplication of unit
Cayley numbers provides the requisite topological algebra.) Hence
$\lambda_{S^7}(K)=0$. In \S\ref{sub-sub-S7} we shall also see that
$\lambda_{S^7}(M)=1$, where $M$ stands for any set of equations
axiomatizing {\em monoids} (otherwise known as {\em associative H-spaces\/}).
One possible form of $M$ comprises the two equations of
(\ref{eq-H-space}) together with the associativity law for $F$.
The proof (details omitted) is like the proof in
\S\ref{sub-sub-S3}, except that it uses I. M. James' 1957 result
\cite{james} that there do not exist $\o e\in S^7$ and a
continuous binary operation $\o F$ on $S^7$ that satisfy $M$
within homotopy.

\subsubsection{{$\lambda_{\infty}(K)$  for $K\,=$
                    H-spaces.}}   \label{sub-sub-infty}
We briefly consider values of $\lambda_{\infty}(\Sigma)$, where
$\infty$ stands for a union of two topological
circles $S_1$, $S_2$ joined at
a single point $P$. (For example a lemniscate is such a space.) In
particular we shall examine values of  $\lambda_{\infty}(H)$,
where $H$ denotes the theory of H-spaces, defined by Equations
(\ref{eq-H-space}) of \S\ref{sub-sub-S2-epsilon}. We will examine
the variation in the value of $\lambda_{(\infty,d)}(H)$ that
occurs among diameter-1 metrics $d$ that define the usual topology
on $\infty$.

For the moment we work only on making $\lambda_{\infty}(H)$ small
by a suitable choice of the metric $d$. Consider the ellipses
$E_1(\varepsilon)$ and $E_2(\varepsilon)$ in 3-space, parametrized
respectively by
\begin{align*}
  \gamma_1(t) &\EQ \frac{1}{2}\, (\cos t,\,\sin
                                    t,\,\varepsilon(1+\cos t))
                    \;\;\;\;\;\text{and} \\
 \gamma_2(t) &\EQ \frac{1}{2}\, (\cos t,\,\sin
                                    t,\,0).
\end{align*}
The union $E\EQ E_1(\varepsilon)\cup E_2(\varepsilon)$ is
homeomorphic to $\infty$, with $(-1/2,0,0)$ corresponding to the
point $P$ common to $S_1$ and $S_2$, and so we give $E\EQ\infty$
the Euclidean metric of $\Reals^3$ in which it is embedded. In
this $d$, the diameter of $\infty$ is $\sqrt{1+\varepsilon^2}$,
approximately
$1 + \frac{1}{2}\varepsilon^2$, with a limit of diameter
$1$ as $\varepsilon\rightarrow0$. (The reader may re-scale to
diameter 1 if desired.)

To estimate $\lambda_{(\infty,d)}(H)$, we define a continuous
function $\psi \FROM E \TO E$ via
\begin{align*}
         \psi(x,y,z) \EQ (x,y,0).
\end{align*}
It is immediate from the geometry that, for all $\mathbf x\in E$,
\begin{align}
           d(\mathbf x,\psi(\mathbf x))\;\leq\;\varepsilon.
\end{align}
The image of $\psi$ is the ellipse $E_2$, and one easily sees that
\begin{align}
                    \psi\upharpoonright E_2 \EQ\text{identity}.
\end{align}

The reader may check that we now have the hypotheses of Theorem
\ref{thm-inequal-retract} of \S\ref{sub-inequal-retract} (with
$E\EQ\infty$ for $A$, $E_2$ for $B$, and $\varepsilon$ for $K$).
In this context, the conclusion of the theorem is
\begin{align*}
        \lambda_{(\infty,d)}(H^{\star})
        \;\leq\;\lambda_{(E_2,d)}(H^{\star})
                     \,+\,\varepsilon.
\end{align*}
Moreover $E_2$ is topologically a circle, which can be made into
an H-space; hence $\lambda_{(E_2,d)}(H^{\star})=0$. We now have
the desired conclusion that
$\lambda_{(\infty,d)}(H^{\star})<\varepsilon$ (which obviously
remains true after rescaling to diameter $1$).

\subsection{$\Sigma$ related to group theory}   \label{sub-groups}

\subsubsection{Group theory on spaces with the
fixed-point property.} \label{subsub=group-f.p.}

Let $\Sigma$ be a consistent set of equations, in unary $-$ and
binary $+$, containing
\begin{align}         \label{eq-gp-fp-a}
           x_1 &\WAVY \{x_0 + [((-x_0) + x_1) + x_2]\} + (- x_2)\\
           x_0 &\Wavy \{x_0 + x_1\} + (-x_1).\label{eq-gp-fp-b}
\end{align}
For instance, $\Sigma$ could be a conservative extension of group
theory. Notice that whatever we can prove about such a $\Sigma$
will also be true for $\Sigma^{\star}$; hence
\S\ref{subsub=group-f.p.} applies also to $\Sigma^{\star}$.

We now take $A$ to be any space that has the {\em fixed-point
property\/}; i.e., we shall assume that for any continuous
function $F\FROM A\TO A$, there exists $e\in A$ with $F(e)=e$. (By
the Theorem of Brouwer, every finite-dimensional closed simplex
has this property. More generally, the class of such spaces is
closed under retraction, and every compact convex subset of a
locally convex topological vector space has this property
\cite[Theorem 15.4]{hu}.) We will
show that\footnote{%
In this context, we define $\infty/2$ to be $\infty$.} %
$\lambda_A(\Sigma)\geq \text{diam}(A)/2$.

We will describe what happens for $\text{diam}(A)<\infty$;
the other case can be left to the reader.
Let us fix an $\varepsilon>0$. By definition of $\lambda$,
there exist continuous operations $+$ and $-$ that satisfy
our equations within $\lambda_A(\Sigma) + \varepsilon$. By definition
of the diameter, there exist $a,b\in A$ with
$d(a,b)>\text{diam}(A)-\varepsilon$. Let $e$ be a fixed point
to the continuous function
\begin{align*}
            x \;\GOESTO ((-a)+b)+x.
\end{align*}
Then we have
\begin{align*}
    b \WAVY \{a + [((-a) + b) + e]\}+ (- e) \EQ \{ a + e\}+ (-e)
                \WAVY a,
\end{align*}
from which it follows that
\begin{align*}
      d(a,b) \;\leq \;2(\lambda_A(\Sigma)+\varepsilon).
\end{align*}
Thus we have
\begin{align*}
     \text{diam}(A) -\varepsilon
                  \;< \;d(a,b) \;&\leq\; 2\lambda_A(\Sigma)\,+\,2\varepsilon\\
     2\lambda_A(\Sigma) \;&\geq\; \text{diam}(A) \,-\,3\varepsilon.
\end{align*}
Letting $\varepsilon$ approach zero, we obtain $\lambda_A(\Sigma)
\geq \text{diam}(A)/2$.

If $A$ is a space whose radius is half its diameter (e.g. $[0,1]^k$
for finite $k$), then by \S\ref{subsub-lambda} we have
$\lambda_A(\Sigma)=\text{radius}(A)$.

\subsubsection{Group theory on spheres $S^n$ for various $n$.}
 Here we recapitulate
the results of \S\ref{sub-spheres-simple} that have to do with
group theory. Recall that our spheres all have diameter $1$. If
$n\neq1,3,7$, then $\lambda_{S^n}(\Sigma) = 1$, for every
demanding theory $\Sigma$. Thus $\lambda_{S^n}$ takes the value
$1$, for even the simplest interesting generalizations of group
theory, such as $H$-spaces.

As for the exceptional values $1,3,7$, $\lambda_{S^1}$ takes the
value $0$ on Abelian groups. $\lambda_{S^3}$ takes the value $1$
on Abelian groups, but $0$ on groups. $\lambda_{S^7}$ takes the
value $1$ on groups, even on monoids, but takes the value $0$ on
H-spaces.

\subsubsection{Group theory on $S^1\!\times Y$.}
                    \label{subsub-group-cyl}
For \S\ref{subsub-group-cyl} we take $\Gamma$ to be any set of
equations in a binary operation $+$ and a unary $-$ that contains
the three equations
\begin{align}          \label{eq-cyl-groups-a}
       x + ((-y) + y) \Wavy x   \\
        x + ((-x)+y) \Wavy y \label{eq-cyl-groups-b}\\
            (x+y)+(-y) \Wavy x. \label{eq-cyl-groups-c}
\end{align}
(These three are not the same as Equations
(\ref{eq-ext-group-thy}--\ref{eq-ext-group-thy-b}) below, and neither
are they the same as (\ref{eq-gp-fp-a}--\ref{eq-gp-fp-b}) above,
but all
three sets are consequences of any version of group theory based
on $+$ and $-$.)

Let $S^1$ be the unit circle with distance defined by arc length,
scaled (as in \S\ref{sub-spheres-simple}) to give $S^1$ a diameter
of $1$. Let $Y$ be the triode that is also discussed in
\S\ref{subsub-Y} below. More precisely, $Y$ is the union of three
unit-length segments in the plane, meeting at one point and at
120-degree angles. $Y$ is given the metric inherited from the
plane. We equip their product with the $L^1$-metric
\[
             d(\,(a,b),\,(c,d)\,) \EQ d_{S^1}(a,c) \,+\, d_Y(b,d).
\]

Let $e$ denote the center point of $Y$---the meeting point where
three angles occur. The set $S^1\times\{e\}$ will be denoted $C$,
and called the {\em central ring} of $S^1\times Y$. If $S$ is one
of the three segments forming $Y$, then $S^1\times S$ is called a
{\em flange} of $S^1\times Y$. Each two of the three flanges have
$C$ as their intersection. The reader who desires to do so may
visualize $S^1\times Y$ as a cylinder with an added flange (since
two of the flanges make a cylinder). Locally
---that is, if one looks at $U\times Y$ for $U$ a segment in
$S^1$---it may be seen as three rectangles in space (each of the
form $(U\times Y)\cap S$), concurrent along the single segment
$(U\times Y)\cap C$ (thereby forming three dihedral angles of 120
degrees\footnote%
{Of course the dihedral angle measure is unimportant topologically;
it is, however, important for the precise metric, and perhaps
it helps for making a mental picture.}%
). Clearly $S^1\times Y$ is
non-homogeneous as a topological space, and hence incompatible
with group theory. Here we will prove the sharper result that
$\lambda_{S^1\times Y}(\Gamma)\,\geq\,0.1$.

For a contradiction, we assume that $\lambda_{S^1\times
Y}(\Gamma)\,<\,0.1$. Thus there are continuous operations,
$\boxplus$ binary and $\boxminus$ unary, satisfying (\ref
{eq-cyl-groups-a}--\ref{eq-cyl-groups-b}) within 0.1 on $S^1\times
Y$.

For $a\in A$, we consider the continuous translation $\tau_a\FROM
A\TO A$ given by $\tau_a(x)=x\boxplus a$. We first prove that
$\tau_a$ {\em is one-one up to} $0.2$, by which we mean that
\begin{align}                 \label{eq-almost-one-one}
      \text{if}\quad \tau_a(x)=\tau_a(y),
                     \quad \text{then} \quad d(x,y)<0.2.
\end{align}
 Well, given
$\tau_a(x)=\tau_a(y)$, we may calculate
\begin{align*}
           x \Wavy (x\boxplus a)\boxplus (\boxminus a) \EQ
                (y\boxplus a)\boxplus (\boxminus a) \Wavy y,
\end{align*}
by two applications of (\ref{eq-cyl-groups-c}). Now
(\ref{eq-almost-one-one}) follows by the triangle inequality.

Let $S_i$ ($1\leq i\leq3$) denote the three flanges of $S^1\times
Y$; for $i=1,2$ choose $a_i\in S_i$ with $d(a_i,C) > 0.8$ for each
$i$. Let us take $a$ to be $(\boxminus a_1)\boxplus a_2$, and
consider $\tau_a$ as defined above. By (\ref{eq-cyl-groups-b}) we
have
\begin{align*}
      \tau_a(a_1) \EQ a_1\boxplus ((\boxminus a_1)\boxplus a_2)
           \Wavy a_2;
\end{align*}
in other words, $d(a_2,\tau_a(a_1))<0.1$, and so by the triangle
inequality we have $d(\tau_a(a_1),C)>0.7$.

Let $R_{12}$ denote the cylinder $S_1\cup S_2$. (It is isometric
to $S^1\times [-1,1]$.) Let $\gamma(t)$ be the parametrized
straight path in $R_{12}$ that goes in constant speed from
$\gamma(0)=a_2$ to $\gamma(1)=a_1$. Consider now the homotopy
\[
       \theta_t=\tau_{(\boxminus\gamma(t))\boxplus a_2}
          \FROM S^1\times Y\TO S^1\times Y.
\]
Clearly $\theta_1=\tau_a$, and by (\ref{eq-cyl-groups-a}) we have
$\theta_0$ within $0.1$ of the identity function.

Let $c$ be the midpoint of segment $[a_1,a_2]$ in the cylinder
$R_{12}$. Clearly $c$ lies on the central ring $C$. Moreover
$\gamma(0.5)=c$, so that $\theta_{0.5}$ moves $c$ to within 0.1 of
$a_2$ (reasoning as above). Following the trajectory of the ring
$C$ under $\theta_t$ as $t$ goes from $0$ to $1$, we see one point
$c$ moves into $S_2$ at $0.5$. Then $\theta_t(C)$ is constrained
to the far side of $\theta_t(a_1)$ in $S_2$, by
(\ref{eq-almost-one-one}). Thus {\em part of the curve
$\theta_1(C)=\tau_a(C)$ is farther than $0.7$ from $C$}.

By continuity, some of $\tau_a(S_3)$ must lie in $S_2$ farther
than $0.7$ from $C$. Let $Q$ be the outer edge of $S_3$. In
considering a line from $\tau_a(c)$ to $\tau_a(a_1)$, let us
consider whether we meet a point $\tau_a(q)$ for $q\in Q$. There
are actually three cases here.\vspace{0.1cm}

\hspace{-\parindent}%
{\bf Case 1.} $\tau_a(b_1) =\tau_a(q)$ for some $b_1\in S_1$ and
some $q\in Q$. This clearly contradicts (\ref{eq-almost-one-one}),
and so the proof is complete.  \vspace{0.1cm}

\hspace{-\parindent}%
{\bf Case 2.} $\tau_a[a_1,c]\subseteq\tau_a[S_3]$. In this case
$\tau_a(a_1)=\tau_a(b_3)$ for some $b_3\in S_3$. Again this
contradicts (\ref{eq-almost-one-one}). \vspace{0.1cm}

\hspace{-\parindent}%
{\bf Case 3.} $\tau_a[a_1,c]\cap\tau_a[S_3]\,=\,\emptyset$. In
considering a line from $\tau_a(c)$ to $\tau_a(a_2)$, let us
consider whether we meet a point $\tau_a(q)$ for $q\in Q$. Here
there are two cases analogous to Cases 1 and 2 above; details
omitted. This completes the proof that $\lambda_{S^1\times
Y}(\Gamma)\geq 0.1$.

These $\lambda$-values---more precisely, the fact that they are
non-zero---will be needed for Theorem
\ref{th-prod-smaller-factors} of \S\ref{sub-prod-on-prod}.

\subsubsection{Group theory on a thickening of $S^1$.}
           \label{subsub-distort}
It is occasionally possible to identify $\lambda_A(\Sigma)$ as the
size of some geometric feature of the metric space $A$. In the
example that follows, $A_{\alpha}$ is a certain metric subspace of the
right circular cylinder of \S\ref{subsub-group-cyl}, with $2\alpha$
its height in the axial direction. We shall describe this space more
thoroughly, and then prove that
$\lambda_{A_{\alpha}}(\Gamma)=\alpha$, where  $\Gamma$ is given by the
group-theoretic equations appearing in
(\ref{eq-ext-group-thy}--\ref{eq-ext-group-thy-b}) below.

For $\alpha$ any real number with $0<\alpha<1$, we define
\begin{align*}
    A_{\alpha}\EQ \{\,(x,y,z)\,:\; &x^2+y^2=1 \;\&\\
        &(-\alpha x\leq z\leq\alpha x\;\;\text{or}\;\;
           z=0)\, \}\;\subseteq\;\Reals^3.
\end{align*}
This space may easily be sketched as a subset of a cylinder in
$\Reals^3$---such as was defined in \S\ref{subsub-group-cyl}. We
give it the rectangular or taxicab metric in that space:
$d(\mathbf x,\mathbf y)=\sum|x_i-y_i|$. (Notice that the spaces
$A_{\alpha}$ are all homeomorphic one to another, but the
homeomorphisms are not isometries.)

Notice that for $(x,y,z)\in A_{\alpha}$ with $x<0$, the definition
yields $z=0$ as the only possible value for $z$. Thus $A_{\alpha}$
contains the circle
\begin{align*}
           C\EQ\{\,(x,y, 0)\,:\, x^2 + y^2 = 1\,\},
\end{align*}
and for negative $x$, these are the only points in $A_{\alpha}$.
For positive $x$, there are other points $(x,y,z)$. The farthest
of these from the circle $C$ are $(1,0,\alpha)$ and
$(1,0,-\alpha)$. Thus $\alpha$ is a measure of how far
$A_{\alpha}$ extends away from the circle $C$.

The map $(x,y,z)\GOESTO(x,y,0)$ is a retraction that moves no
point farther than $\alpha$. Its image is $C$, which is
homeomorphic to $S^1$, and thus compatible with group theory, and
so by Theorem \ref{thm-inequal-retract} of
\S\ref{sub-inequal-retract}, if $\Gamma$ is any consistent set of
equations holding in group theory (regardless of the precise
similarity type used to express group theory, and clearly
regardless of any distinction between $\Gamma$ and
$\Gamma^{\star}$), then
\begin{align}                  \label{eq-A-alpha-lower}
        \lambda_{A_{\alpha}}(\Gamma)\;\leq\;\alpha.
\end{align}

For the reverse inequality, we restrict $\Gamma$ to be a set of
equations in the operations $+$ and $-$ that contains the three
equations
\begin{gather}   \label{eq-ext-group-thy}
          (x + y) + (-y) \WAVY x \\
              x + 0 \WAVY x; \quad \quad 0+x\WAVY x.
                  \label{eq-ext-group-thy-b}
\end{gather}
We shall prove that for such a $\Gamma$,
\begin{align}             \label{eq-A-alpha-upper}
        \lambda_{A_{\alpha}}(\Gamma)\;\geq\;\alpha.
\end{align}
Now if $\lambda_{A_{\alpha}}(\Gamma)\geq 1$, then
(\ref{eq-A-alpha-upper}) holds by our assumption on $\alpha$.
Hence we will assume from now on that
$\lambda_{A_{\alpha}}(\Gamma)<1$. We will consider an arbitrary
$\varepsilon$ with $0<\varepsilon<1$, and will prove that
$\lambda_{A_{\alpha}}(\Gamma)+\varepsilon
>\alpha$.

We define a closed curve $f$  in $A_{\alpha}$ (for $0\leq t\leq
2\pi$), as follows:
\begin{align*}
    f(t) &\EQ \begin{cases}
              \;\; (\cos t,\,\sin t,\,-\alpha\cos t)  &\text{if
                           $\cos t\geq0$} \\
               \;\;(\cos t,\,\sin t,\,0)
                &\text{if $\cos t\leq0$}.
          \end{cases}
\end{align*}
($f$ maps, so to speak, to the lower periphery of $A_{\alpha}$.)
Concerning the point
$
        A \EQ (1,0,\alpha)\in A_{\alpha},
$ let us check that $A$ has distance at least $2\alpha$ from every
point $B=(\cos t,\cdots)$ in the image of $f$. If $\cos t<0$, then
\begin{align*}
     d(A,B) \EQ (1-\cos t) \,+ \,\cdots  \,+\, (\alpha-0)\;>\;
                  1\,+\,\alpha\;>\;2\alpha.
\end{align*}
 On the other hand, if $\cos t\geq
0$, then
\begin{align*}
          d(A,B) &\EQ |1-\cos t| \,+\, |\sin t| \,+\,
                        |\alpha+\alpha\cos t|\\
              &\;\geq\; |1-\cos t| \,+\, |\alpha+\alpha\cos t|\\
              &\EQ 1-\cos t \,+\, \alpha+\alpha\cos t\\
              &\EQ 2\alpha \,+\, (1-\cos t)(1-\alpha)\;\geq\;
              2\alpha,
\end{align*}
where the final inequality comes from our assumption that
$\alpha<1$.

For this proof, when we refer to a closed curve, we mean a
continuous map with domain $S^1$. We view our $f(t)$ as such a
closed curve, by representing $S^1$ as $\Reals/2\pi$, and relying
on the periodicity of the trigonometric functions. The same
applies to curves constructed from $f(t)$, such as in
(\ref{eq-a+0}) and (\ref{eq-a+B}) below. Finally when we say that
closed curves $g_0(t)$ and $g_1(t)$ in $A$ are homotopic, we mean
that there exists a map $G\FROM S^1\times [0,1]\TO A$ such that
$G(t,i) = g_i(t)$ for $t\in S^1$ and $i\in\{0,1\}$.

By definition of $\lambda$, there exist continuous operations
$\boxplus$ and $\boxminus$ satisfying $\Gamma$ within
$\lambda_{A_{\alpha}}(\Gamma)+\varepsilon$ on $A_{\alpha}$. As a
consequence of (\ref{eq-A-alpha-lower}), these operations satisfy
$\Gamma$ within $1+\varepsilon<2$. We consider the closed curve
\begin{align}           \label{eq-a+0}
             t \GOESTO 0 \boxplus f(t)
\end{align}
By the second equation of (\ref{eq-ext-group-thy-b}) it must stay
within $\lambda_{A_{\alpha}}(\Gamma)+\varepsilon$ of $f(t)$, and
hence within $<2$ of $f(t)$.
Therefore the curve (\ref{eq-a+0}) is homotopic to $f(t)$, by the
homotopy that was introduced in \S\ref{sub-sub-Sn}. Since there is
a path connecting $0$ to $A$, the path
\begin{align}           \label{eq-a+B}
             t \GOESTO A \boxplus f(t)
\end{align}
is also homotopic to $f(t)$. Thus (\ref{eq-a+B}) maps onto
$\{f(t):\pi/2\leq t\leq3\pi/2\}$, and in particular, there exists
$t_0$ such that
\begin{align*}
           A \boxplus f(t_0) \EQ (-1,0,0).
\end{align*}
By reasoning similar to that for  (\ref{eq-a+B}), except using the
first equation of (\ref{eq-ext-group-thy-b}), the map
\begin{align*}
             s \GOESTO f(s) \boxplus f(t_0)
\end{align*}
is homotopic to $f(s)$, and so there exists $s_0$ with
\begin{align*}
           f(s_0) \boxplus f(t_0) \EQ (-1,0,0).
\end{align*}
Thus in particular we have
\begin{align*}
          f(s_0) \boxplus f(t_0) \EQ A \boxplus f(t_0) .
\end{align*}
We now calculate, using (\ref{eq-ext-group-thy}):
\begin{align*}
    A \WAVY (A \boxplus f(t_0)) \boxplus (\boxminus f(t_0)) \EQ
    (f(s_0) \boxplus f(t_0)) \boxplus (\boxminus f(t_0)) \WAVY
    f(s_0),
\end{align*}
where $\wavy$ refers to approximate equality within
$\lambda_{A_{\alpha}}(\Gamma)+\varepsilon$.
We proved above that $d(A,f(s))\geq2\alpha$. Thus the triangle
inequality now yields
\begin{align*}
     2\alpha&\;\leq\; d(A,f(s))\;\leq\; 2
     (\lambda_{A_{\alpha}}(\Gamma)+\varepsilon)\\
     &\EQ 2\lambda_{A_{\alpha}}(\Gamma)\,+\,2\varepsilon.
\end{align*}
Since $\varepsilon$ can be taken arbitrarily small and positive,
we now have the estimate (\ref{eq-A-alpha-upper}).

Combining (\ref{eq-A-alpha-lower}) with (\ref{eq-A-alpha-upper}),
we see that if $\Gamma$ is a set of equations in $+$ and $-$ that
contains Equations
(\ref{eq-ext-group-thy}--\ref{eq-ext-group-thy-b}), and if each
equation of $\Gamma$ holds in group theory, then
$\lambda_{A_{\alpha}}(\Gamma)=\alpha$.

Also note that $A_{\alpha}$ is compatible with H-space theory.
(The proof is left to the reader.)




\subsubsection{Boolean algebra}    \label{subsub-boolean}
In 1947, I. Kaplansky proved \cite{kapl-amerj}, inter alia,
that if $\mathbf A$
is a compact topological Boolean ring (equivalently, a compact
topological Boolean algebra), then $\mathbf A$ is isomorphic (as a
topological algebra) to a power of the two-element discrete
Boolean ring. In particular, this result says that if $A$ is a
compact metrizable space that is not homeomorphic to a power of
the two-element discrete space, then $A$ is not compatible with
Boolean rings.

Kaplansky's 1947 proof\footnote{%
The result was independently discovered, with a similar proof,
by P.S. Rema in 1964 \cite{rema}.} %
 is probably the most sophisticated proof of
incompatibility on record: it relies on the result that if $G$ is
a locally compact Abelian group and $g\in G$ with $g$ not the
identity element, then there exists a group character $f\FROM
G\TO\Complex$ with $f(g)\not=1$. The present author has no clue
how such a method could be extended to yield a lower estimation
for $\lambda$ in this context. The trouble is that when the group
axioms are relaxed so as to hold only approximately, the entire
apparatus of group duality becomes unavailable. (Some approximate
version of it might be possible, but that would appear to be a
difficult prospect indeed, requiring serious theory-building.)

Therefore, it seems unlikely that we shall soon see a theorem of
the form
\begin{align}           \label{eq-est-ba2}
          \lambda_A(BA)>K
\end{align} for $BA$ taken
to be Boolean ring theory. There is indeed another reason
(\ref{eq-est-ba2}) is unlikely: for most spaces $A$ it would
probably be more practical and more meaningful to obtain an
estimate of $\lambda_A(\Gamma)$ or of $\lambda_A(\Lambda)$ (group
theory or lattice theory), since the exclusion of either of these
will exclude a Boolean ring structure. Therefore, an inequality
like (\ref{eq-est-ba2}) would not have a significant role unless
$A$ is a compact metric space that is compatible with group theory
and with lattice theory, but is not a Cantor space. Such spaces
seem rare.

In 1969, T. H. Choe gave a simpler proof \cite{choe} of
Kaplansky's result, but only for spaces of finite dimension. It is
not clear how Choe's proof could be turned into an estimate of
$\lambda_A(BA)$ for some space $A$.
\subsubsection{$A=\Reals$ ; $\Sigma=$ groups of exponent 2.}
   \label{subsub-not-exp2}
Here we suppose that $\Sigma$ is a finite set of equations in $+$
and $0$ that includes
\begin{align*}
            x_0 + x_0 &\WAVY 0 \\
            0 + x_0   &\WAVY x_0\\
            x_0 + (x_0 + x_1) &\WAVY x_1.
\end{align*}
For instance, $\Sigma$ could be an axiomatization of the theory
of groups of exponent $2$. Here we shall prove that, in any metric
for $\Reals$,
$\lambda_{\Reals}(\Sigma)\geq\text{radius}(\Reals)/2$. (We mean
this to include the assertion that
$\lambda_{\Reals}(\Sigma)=\infty$ whenever $\Reals$ is metrized
with infinite radius.)

For a contradiction, let us suppose that $\lambda_{\Reals}(\Sigma)
=K <\text{radius}(\Reals)/2$. Therefore, there exist a constant
$\o0\in\Reals$ and a continuous binary operation $\boxplus$ on
$\Reals$ that satisfy $\Sigma$ within $K$. Since
$2K<\text{radius}(\Reals)$, the $2K$-ball centered at $\o0$ is not
all of $\Reals$. In other words, there exists $a\in\Reals$  with
$d(a,\o0)\,>\,2K$. We may assume, without loss of generality, that
$\o0<a$; in this case, we of course have
\begin{align}               \label{eq-02K}
      \o0+K<a-K.
\end{align}

We consider the continuous real-valued function
\begin{align*}
              x \GOESTO \phi(x) \EQ x\boxplus a.
\end{align*}
It follows readily from $\Sigma$ that
\begin{align*}
            \phi(\o0) &\EQ \o0\boxplus a \WAVY
            a,\quad\quad\quad\text{and}\\
            \phi(a)   &\EQ a\boxplus a \WAVY\o0.
\end{align*}
From these estimates, and (\ref{eq-02K}), we deduce
\begin{align*}
             \phi(\o0) \;&>\; a -K \;>\; \o0 +K \,>\,\o0\\
             \phi(a) \;&<\; \o0 +K \;<\; a -K \,<\,a
\end{align*}

These last inequalities display a sign change for the function
$\phi(x)-x$; hence, by the IVT, $\phi$ has a fixed point, i.e.
$c\,=\,\phi(c)\,=\,c\boxplus a\,$ for some $c\in\Reals$. From
$\Sigma$, we have
\begin{align*}
          a \WAVY c\boxplus(c\boxplus a) = c\boxplus c \WAVY
          \o0.
\end{align*}
Thus $d(a,\o0)\,<\,2K$, in contradiction to our specification of
$a$. This contradiction completes the proof that
$\lambda_{\Reals}(\Sigma)\geq\text{radius}(\Reals)/2$.

\subsubsection{$A=\Reals^2$ ; $\Sigma=$ groups of exponent 2.}
   \label{subsub-square-not-exp2}
In \S\ref{subsub-square-not-exp2} we take $\Gamma_2$ to be a set
of axioms for groups of exponent $2$. (Or, one might wish to
consider the weaker equations of \S\ref{subsub-not-exp2}. One may
also extend the problem to $\Reals^n$ for all $n\geq2$.)
\vspace{0.1in}

\hspace*{-\parindent}%
{\bf Problem}\hspace{0.1in}
$\lambda_{\Reals^2}(\Gamma_2)\EQ\infty\,$?\vspace{0.1in}

We will sketch a proof that $\Gamma_2$ is not compatible with
$\Reals^2$, and then comment on a (still unknown) strengthening of
that proof that might establish that
$\lambda_{\Reals^2}(\Gamma_2)>0$ (and hence $\EQ\infty$).

Suppose $(\Reals^2,\o F,\o e)$ is a group of exponent $2$, with
multiplication $\o F$ and unit element $\o e$. Choose
$a\in\Reals^2$ with $a\neq \o e$. The map $\phi:x\GOESTO F(x,a)$
is a continuous involution of $\Reals^2$, i.e. a map satisfying
$\phi(\phi(x))=x$. According to a 1934 theorem of P. A. Smith (see
J. Dugundji and A. Granas \cite[Theorem 5.3, page 79]{dug-gran}),
$\phi$ has a fixed point $b$. In other words, we have $F(b,a)=b$,
with $a\neq\o e$. This contradiction to the laws of group theory
establishes that there can be no such continuous group operation
$\o F$ of exponent $2$.

Now in order to extend this proof to approximate models of
$\Gamma_2$, we need some way of finding an approximate fixed point
for an approximate involution of $\Reals^2$, in other words, an
approximated version of Smith's theorem mentioned above. As far as
the author is aware, no such result is available.

\subsubsection{Two auxiliary theories.}
                   \label{subsub-aux-to-log}
In \S\ref{subsub-aux-to-log} we consider two infinite theories,
which we shall (in this section and the next) denote $\Sigma_1$
and $\Sigma_2$. These are the only infinite theories that we shall
consider as such. (The two estimates derived in
\S\ref{subsub-aux-to-log} both require infinitely many equations;
they---and the corollary estimate in \S\ref{subsub-l-groups}---are
the only such estimates in the paper.) For $i=1,2$, and for any
compact metric space $A$, we shall have that
$\lambda_A(\Sigma_i)\geq\text{diameter}(A)/4.$ The result for
$\Sigma_2$ will be applied in \S\ref{subsub-l-groups} to make a
similar estimate for the (more naturally occurring) variety of
lattice-ordered groups.

Our methods for estimating $\lambda_A(\Sigma_i)$ ($i=1,2$) are
similar, but different enough that we will present both in some
detail. In both cases we consider points $a,b\in A$ whose distance
is the diameter of $A$, and apply the triangle inequality to a
certain triangle $\Tri{a}{b}{e}$. In one case, $e$ will be the
value of a term-function $e=\o K(a,b)$; in the other case it will
be the limit of a sequence of such values.

 $\Sigma_1$ will be the following infinite set of
equations
\begin{align}                \label{eq-old-fix-p-a}
                 F(\phi^k(x),x,y)&\Wavy x \\
                     F(x,x,y)&\Wavy y, \label{eq-old-fix-p-b}
\end{align}
for $k\in\omega$, $k\geq 1$. We introduced this theory in
1986---see \cite[\S3.18, page 35]{wtaylor-cots}---and proved that
it is incompatible with every compact Hausdorff space. Here we
shall prove the stronger result that {\em if $A$ is a compact
metric space, then}
\begin{align}                    \label{eq-diameter-infinite-set}
           \lambda_A(\Sigma_1)\;\geq\; \text{diameter}(A)/4.
\end{align}

It is easiest to prove this inequality by contradiction. If
(\ref{eq-diameter-infinite-set}) fails, then there is a
topological algebra  $\mathbf A = \la A, \o\phi, \o F\ra$ with
$\lambda_{\mathbf A}(\Sigma_1)\;<\; \text{diameter}(A)/4$.
Therefore there is a positive real number $\varepsilon$ such that
\begin{align*}
          \lambda_{\mathbf A}(\Sigma_1)\;<\; \text{diameter}(A)/4 -
                      \varepsilon.
\end{align*}

Let $a$ and $b$ be points of $A$ with $d(a,b)$ equal to the
diameter of $A$. Consider the sequence $\o\phi^i(a)$; by
compactness it has a convergent subsequence:
\begin{align*}
       \lim_{i\longrightarrow\infty} \o\phi^{n(i)}(a)\EQ c\in A.
\end{align*}
By the triangle inequality, we have either
$d(a,c)\geq\text{diameter}(A)/2$ or
$d(b,c)\geq\text{diameter}(A)/2$. Without loss of generality, we
will assume that
\begin{align}             \label{eq-distance-b-c}
               d(b,c)\geq\text{diameter}(A)/2.
\end{align}

We next consider the infinite sequence in $A$,
\begin{align*}
         \alpha_i \EQ \o F(\o\phi^{n(i+1)}(a),\,\o\phi^{n(i)}(a),\,b).
\end{align*}
By our choice of the subsequence $\o\phi^{n(i)}(a)$, and by the
continuity of $\o F$, there exists $i_0$ such that
\begin{align}                        \label{eq-F-continuity}
      d(\alpha_{i_0},\,\o F(c,c,b))\;&<\:\varepsilon
                    \quad\text{and}\\
             d( \o\phi^{n(i_0)}(a),\,c)\;&<\; \varepsilon.
                       \label{eq-F-continuity-b}
\end{align}
Now by the approximate satisfaction of
(\ref{eq-old-fix-p-a}--\ref{eq-old-fix-p-b}) we have
\begin{align*}
            d (\o F( c, c,b),\, b) \;&<\; \text{diameter}(A)/4 -
                                   \varepsilon\\
            d (\alpha_{i_0},\,\o\phi^{n(i_0)}(a))
                         \;&<\; \text{diameter}(A)/4 -
                                   \varepsilon.
\end{align*}
Combining these two inequalities with
(\ref{eq-F-continuity}--\ref{eq-F-continuity-b}), via the triangle
inequality, yields $d(b,c)<\text{diameter}(A)/2$. This
contradiction to (\ref{eq-distance-b-c}) completes the proof of
(\ref{eq-diameter-infinite-set}).

Now let us take $\Sigma_2$ to be the following (doubly infinite)
set of equations:
\begin{gather}         \label{eq-old-fix-p-a-bis}
           G(\psi_{m+k}(x,y),\,\psi_m(x,y),\,x,\,y) \Wavy x\\
             K(x,y)\Wavy G(u,u,x,y)\Wavy K(y,x),
                   \label{eq-old-fix-p-b-bis}
\end{gather}
for $m,k\in\omega$, with $k\geq1$. We shall again establish that
if $A$ is compact metric space, then
\begin{align}                    \label{eq-diameter-infinite-set-bis}
           \lambda_A(\Sigma_2)\;\geq\; \text{diameter}(A)/4.
\end{align}
It is easiest to prove this inequality by contradiction. If
(\ref{eq-diameter-infinite-set-bis}) fails, then there is a
topological algebra  $\mathbf A = \la A, \o\psi, \o G, \o K\ra$
with $\lambda_{\mathbf A}(\Sigma_2)\;<\; \text{diameter}(A)/4$.
Therefore there is a positive real number $\varepsilon$ such that
\begin{align*}
          \lambda_{\mathbf A}(\Sigma_2)\;<\; \text{diameter}(A)/4 -
                      \varepsilon.
\end{align*}

Let $a$ and $b$ be points of $A$ with $d(a,b)$ equal to the
diameter of $A$. By the triangle inequality, we have either
$d(a,\o K(a,b) \geq \text{diameter}(A)/2$ or $d(b,\o K(a,b) \geq
\text{diameter}(A)/2$. Without loss of generality, we shall assume
that
\begin{align}                      \label{eq-distance-b-K}
             d(b,\o K(a,b)) \;\geq\; \text{diameter}(A)/2
\end{align}
Consider the sequence $\o\psi_i(b,a)$; by compactness it has a
convergent subsequence:
\begin{align*}
       \lim_{i\longrightarrow\infty} \o\psi_{n(i)}(b,a)\EQ c\in A.
\end{align*}
We next consider the infinite sequence in $A$
\begin{align*}
              \beta_i \EQ \o
                   G(\o\psi_{n(i+1)}(b,a),\,\o\psi_{n(i)}(b,a),\,b,\,a).
\end{align*}
By the continuity of $\o G$, there exists $i_0$ such that
\begin{align}                         \label{eq-G-continuity}
               d(\o G(c,c,b,a),\,\beta_{i_0}) \;<\; 2\varepsilon.
\end{align}

Now by the approximate satisfaction of
(\ref{eq-old-fix-p-a-bis}--\ref{eq-old-fix-p-b-bis}) we have
\begin{align*}
            d (\o G( c, c,b,a),\, b) \;&<\; \text{diameter}(A)/4 -
                                   \varepsilon;\\
            d (\beta_{i_0},\,\o K(a,b))   \;&<\; \text{diameter}(A)/4 -
                                   \varepsilon.
\end{align*}
Combining these two inequalities with (\ref{eq-G-continuity}), via
the triangle inequality, yields $d(b,\o
K(a,b))<\text{diameter}(A)/2$. This contradiction to
(\ref{eq-distance-b-K}) completes the proof of
(\ref{eq-diameter-infinite-set-bis}).

It may be noted that in the proof for $\Sigma_1$ we used only the
continuity of $\o F$, and in the  proof for $\Sigma_2$ we used
only the continuity of $\o G$. We never invoked the continuity of
$\o K$, $\o\phi$ or the $\o\psi_m$'s. Similar occurrences will be
noted below, in \S\ref{subsub-inf-dim} and in Theorem
\ref{th-low-est-squares} of \S\ref{subsub-approx-cubes}.

\subsubsection{Lattice-ordered groups.} We will let $\Lambda\Gamma$ stand for a certain
set of equations (see (\ref{eq-lg-c}--\ref{eq-lg-a}) below) that
is satisfied by the variety of lattice-ordered groups. That
variety has a well-known finite axiomatization, which we
paraphrase below, but do not state formally. Our $\Lambda\Gamma$
is an infinite set of
equations true in LO-groups.\footnote{%
$\Lambda\Gamma$ is obviously weaker than the theory of LO-groups,
since it does not mention $\vee$.  See Equations
(\ref{eq-lg-c}--\ref{eq-lg-a}) below.} %

Like a Boolean algebra, a lattice-ordered group has both lattice
operations $\wedge,\vee$ and group operations $+,-,0$. Therefore
any incompatibilities for these two strong theories (see
\S\ref{sub-groups} and \S\ref{sub-lattices} {\em passim}) will be
inherited by $\Lambda\Sigma$. But LO-groups indeed form a theory
that is stronger than the join of groups and lattices. In addition
to the usual axioms of group theory and lattice theory, it is
assumed that addition distributes over meet and join: $x+(y\wedge
z)\wavy(x+y)\wedge(x+z)$, and dually. The resulting theory is
quite strong. As was proved by M. Ja.\ Antonovski\u\i\ and A. V.
Mironov \cite{anton-mir} in 1967 {\em no compact Hausdorff space
is compatible with lattice-ordered groups}. (Thus, for example,
the Cantor set is a space that is compatible with Boolean algebra,
but not with LO-groups.) Here we will prove, for any compact
metric space $A$, the stronger result that
\begin{align}                 \label{eq-estim-lo-gps}
             \lambda_A(\Lambda\Gamma) \;\geq\; \frac{1}{4}
                        \;\text{diameter}(A).
\end{align}

Our method for proving (\ref{eq-estim-lo-gps}) is very simple: in
\S\ref{subsub-aux-to-log} we proved the corresponding estimate
(\ref{eq-diameter-infinite-set-bis}) for the special theory
$\Sigma_2$. To prove (\ref{eq-estim-lo-gps}), we need only give a
careful definition of $\Lambda\Gamma$, prove that $\Sigma_2$ is
interpretable (\S\ref{sub-interpretable}) in $\Lambda\Gamma$, and
then invoke Theorem \ref{th-interp-lesseq} of
\S\ref{sub-interpretable} (monotonicity of $\lambda_A$ with
respect to the interpretability relation).

$\Lambda\Gamma$ is defined to be the following (doubly infinite)
set of equations:
\begin{align}
   x&\Wavy x\wedge[(z_{m+k} - z_m) \,+\, (x\wedge y)]
                                         \label{eq-lg-c}\\
        x\wedge y &\Wavy   x\wedge[(u\,-\,u) + (x\wedge y)]
                   \Wavy y\wedge x,
                                            \label{eq-lg-a}
\end{align}
where $z_n$ ($n\in\omega$) are terms defined recursively as
follows:
\begin{align*}        
        z_0 \EQ 0;\quad\quad\quad  z_{n+1}\EQ (z_n\,+\,(x-(x\wedge y))).
\end{align*}
Since $x\geq x\wedge y$ in LO-groups, it is easy to see that for
$r\leq s$, LO-groups satisfy $z_r\leq z_s$. To see that LO-groups
satisfy (\ref{eq-lg-c}), we calculate
\begin{align*}
   x &\Wavy [z_m \,+\, (x \,-\, x\wedge y)] 
            \,-\, z_m \,+\, x\wedge y\\
     &\Wavy    z_{m+1} \,-\,z_m\,+\, x\wedge y
     \;\;\leq\;\; z_{m+k} \,-\,z_m\,+\, x\wedge y,  
\end{align*}
for $m,k\in\omega$ with $k\geq1$. The validity of (\ref{eq-lg-a})
in LO-groups is evident. Thus $\Lambda\Gamma$ is a subset of the
equations holding in LO-groups.

Finally, we need to prove that $\Sigma_2\,\leq\,\Lambda\Gamma$,
where $\Sigma_2$ is as defined in  \S\ref{subsub-aux-to-log} above
(Equations (\ref{eq-old-fix-p-a-bis}--\ref{eq-old-fix-p-b-bis})).
In the theory of LO-groups we define the following terms:
\begin{align}          \label{eq-interp-log-a}
        G(u,v,x,y)  &\EQ x\wedge[(u-v)+(x\wedge y)] \\
           K(x,y)   &\EQ  x\wedge y \label{eq-interp-log-b} \\
        \psi_m(x,y) &\EQ z_m. \label{eq-interp-log-c}
\end{align}
These equations define $G(u,v,x,y)$ (resp.\ $K(x,y)$, resp.\
$\psi_m(x,y)$) as the interpreting term for the operation $G$
(resp.\ $K$, resp.\ $\psi_m$). Clearly Equation
(\ref{eq-old-fix-p-a-bis}) interprets as (\ref{eq-lg-c}) and
Equations (\ref{eq-old-fix-p-b-bis}) interpret as (\ref{eq-lg-a}).
Since Equations  (\ref{eq-lg-c}--\ref{eq-lg-a}) define
$\Lambda\Gamma$, the interpretation is valid; we do have
$\Sigma_2\,\leq\,\Lambda\Gamma$.

\subsection{$\Sigma$ related to lattice theory}
                \label{sub-lattices}

\subsubsection{$A=S^1$ (circle); $\Sigma=$ semilattice theory.}
    \label{subsub-lambda-homotopy}
 \label{subsub-l-groups} As
in \S\ref{sub-spheres-simple}, we let $S^n$ stand for the
$n$-sphere of diameter $1$, embedded as a metric subspace of
$\Reals^{n+1}$, in the usual way. In \S\ref{sub-sub-S1} we proved
that $\lambda_{S^1}(\Sigma)=1$ if $\Sigma$ is not Abelian. It is
not hard to see that the theory of semilattices is not Abelian. We
will establish this fact for the weaker theory of an {\em
idempotent commutative binary operation}:
\begin{align}                \label{eq-com-id}
          F(x,y)\Wavy F(y,x);\quad\quad F(x,x)\Wavy x.
\end{align}
Suppose we had $(\Integers,\o F)$ modeling Equations
(\ref{eq-com-id}), where $\o F(x,y) = ax + by$ for some integers
$a$ and $b$. The first equation of (\ref{eq-com-id}) forces $a=b$,
and the second forces $a+b=1$. These two equations in $a$ and $b$
have no solution in $\Integers$; hence Equations (\ref{eq-com-id})
do not form an Abelian set. Thus, in particular,
\begin{align*}
            \lambda_{S^1}(\text{idempotent commutative})\EQ 1.
\end{align*}

The reader may take a similar path to discover that
$\lambda_{S^1}(\Sigma)\EQ 1$, where $\Sigma$ is taken to be the
equations of {\em multiplication with zero and one}:
\begin{align}                \label{eq-zero-one}
          F(0,x)\Wavy 0;\quad\quad F(1,x)\Wavy x.
\end{align}
In \cite[\S3.6, page 28]{wtaylor-cots}, we proved that if $A$ is
arcwise-connected and not contractible, then $A$ is not compatible
with the equations (\ref{eq-zero-one}). Nevertheless, for these
equations we have a calculation of a non-zero value of $\lambda$
only for the one space $S^1$.

The reader may also wish to consider the equations of {\em median
algebras}:
\begin{gather}
     m(x,y,z) \WAVY m(x,z,y) \WAVY m(y,z,x) \label{eq-median-a} \\
      m(x,x,z) \WAVY x \label{eq-median-b} \\
      m(m(x,y,z),u,v) \WAVY m(x,m(y,u,v),m(z,u,v)). \label{eq-median-c}
\end{gather}
(See e.g.\ Bandelt and Hedl\'{\i}kov\'a \cite{band-hed}.) The
reader may check that the proof given here in
\S\ref{subsub-lambda-homotopy} for an idempotent commutative
algebra may be applied also for Equations
(\ref{eq-median-a}--\ref{eq-median-c}). (In fact Equations
(\ref{eq-median-a}--\ref{eq-median-b}) suffice for a proof.)
Therefore this theory also has $\lambda_{S^1}$ equal to $1$.

\subsubsection{$S^1$ with a different metric.}
                 \label{subsub-S1-different}
In \S\ref{subsub-S1-different} we will let $\Sigma$ be any theory
compatible with $[0,1]$---such as, for example, the theories
mentioned in \S\ref{subsub-lambda-homotopy} (semilattices,
multiplication with zero and one, median algebra). Here we shall
see that, for each $\varepsilon>0$, there is a metric $d$ on $S^1$
such that $\lambda_{(S^1,d)}(\Sigma^{\star})<\varepsilon$. The
method is like that of \S\ref{sub-sub-S2-epsilon}.

We begin with the realization of $S^1$ as the ellipse that is the
locus of
\begin{align*}
      4\varepsilon^2 x^2 \,+\, 4y^2  \EQ
                        \varepsilon^2
\end{align*}
in ordinary 2-space. Then $d$ is defined to be the ordinary
Euclidean distance of $\Reals^2$, restricted to the ellipsoid. One
easily checks that $(S^1,d)$ has diameter $1$. If we define
$\psi\FROM S^1\TO S^1$ via $\psi(x,y) = (x,|y|)$, then we clearly
have the hypotheses of Theorem \ref{thm-inequal-retract} of
\S\ref{sub-inequal-retract}, with the image of $\psi$ homeomorphic
to $[0,1]$. Therefore
\begin{align*}
           \lambda_{(S^1,d)}(\Sigma^{\star})
           \;\leq\;\lambda_{[0,1]}(\Sigma^{\star})
               \, +\,\varepsilon \EQ \varepsilon.
\end{align*}

\subsubsection{An extension of
             \S\ref{subsub-lambda-homotopy}}
                            \label{subsub-semil-homotop}
Here we consider an arbitrary path-connected and locally
path-connected metric space $A$ whose fundamental group is
isomorphic to $\Integers$ (under ordinary addition), and which has
the following property: {\em there exists $K>0$ such that if\/
$\gamma(t),\delta(t)$ are loops beginning and ending at the same
point\/ $a\in A$, and if\/ $d(\,\gamma(t),\,\delta(t)\,)\,<\,K$
for all\/ $t$, then $\gamma$ and\/ $\delta$ are homotopic as
loops.}

Let $\Sigma$ be an idempotent set of equations that is not Abelian
(as defined in \S\ref{sub-sub-S1}). We shall in fact
require\footnote%
{This requirement would be used in the proof of
(\ref{eq-def-sig-star}), which we skipped. It is used to get each
loop to start and stop at $a_0$, as we do just below in the
definition of $F^{\star}$.} %
that each term appearing in $\Sigma$ shall be declared as
idempotent by an equation in $\Sigma$. (For example, if $\Sigma$
contains $x+(y+z)\wavy(x+y)+z$, then $\Sigma$ must also contain
$x+(x+x)\wavy x$ and $(x+x)+x\wavy x$.) We shall prove that
$\lambda_A(\Sigma)\geq K$. (In \S\ref{subsub-lambda-homotopy} we
had these assumptions true for $K=1$, and the conclusion was
$\lambda_A(\Sigma)\,=\, 1$.)

For a contradiction, let us suppose that $\lambda_A(\Sigma)< K$.
This means that $A\models_{\varepsilon}\Sigma$ for some
$\varepsilon<K$. (See (\ref{eq-models-varep}) of
\S\ref{sub-metric-approx-compat}.) Thus there are operations $\o
F_t$ corresponding to the symbols $F_t$ appearing in $\Sigma$,
such that for each equation $\sigma\wavy\tau$ in $\Sigma$, the
terms $\sigma$ and $\tau$ evaluate to functions $\o\sigma$ and
$\o\tau$ that are within $\varepsilon$ of each other. Our
contradiction will be that $\Sigma$ is Abelian; for this we need
to satisfy $\Sigma$ in an algebra $(\Integers,F^{\star}_t)_{t\in
T}$, where each operation $F^{\star}_t(x_1,x_2,\ldots)$ has the
form $\sum m_ix_i$ in $\Integers$.

We pick a base point $a_0\in A$ and consider homotopy classes
$[\alpha]$ of loops in $A$ at $a_0$. We define
$[\beta]\,=\,F_t^{\star}([\alpha_1],[\alpha_2], \ldots)$ to be the
homotopy class of the loop defined by $\o
F_t(\alpha_1(t),\alpha_2(t),\ldots)$. Since $\o F_t$ is idempotent
within $\varepsilon$, this loop starts and stops at a point $a'$
near $a_0$. Any two paths that connect $a_0$ to $a'$ that stay
within $\varepsilon$ of $a_0$ are homotopic. Hence we may (up to
homotopy) unambiguously append such a path to the beginning, and
its reverse to the end, of $[\beta]$; the modified  $[\beta]$ is
what we take for $[F_t^{\star}([\alpha_1],[\alpha_2], \ldots)$.
Now this loop may be continuously deformed, without moving the
endpoints, to one that begins with $\o
F_t(\alpha_1(t),a_0,a_0,\ldots)$, followed by $\o F_t(a_0,
\alpha_2(t), a_0,\ldots)$, and so on. From this we see that
$F^{\star}([\alpha_1],[\alpha_2],\ldots)$ has the required linear
form $\sum m_i[\alpha_1]$.

For the satisfaction of $\Sigma$ in $(\Integers,F^{\star}_t)_{t\in
T}$, we consider an equation $\sigma\wavy\tau$ of $\Sigma$. We
skip the proof (inductive over complexity of $\sigma$) that the
derived operation $\sigma^{\star}$ corresponding to $\sigma$ in
$(\Integers,F^{\star}_t)_{t\in T}$ satisfies
\begin{align}              \label{eq-def-sig-star}
\sigma^{\star}([\alpha_1],[\alpha_2],\ldots) \Wavy [\sigma']\,,
\end{align}
where
\begin{align}                    \label{eq-def-sig-prime}
           \sigma'(t) \EQ \o\sigma
                        (\alpha_1(t),\alpha_2(t),\ldots).
\end{align}

Now, if $\sigma\wavy\tau$ is in $\Sigma$, then by assumption the
equation $\o\sigma(x_1,x_2,\ldots) = \o\tau(x_1,x_2,\ldots)$ holds
within $\varepsilon$ on $A$. By (\ref{eq-def-sig-prime}) we have
$\sigma'(t)=\tau'(t)$, within $\varepsilon$, for all $t$. By our
assumptions on the metric space $A$, we have $\sigma'$ homotopic
to $\tau'$; in other words that $[\sigma']=[\tau']$. It is
immediate from (\ref{eq-def-sig-star}) that
$\sigma^{\star}\,=\,\tau^{\star}$, as operations on $\Integers$.

One special case of the result of \S\ref{subsub-semil-homotop}
occurs when $\Sigma$ is the theory of an idempotent commutative
operation, as presented in (\ref{eq-com-id}), and $A$ has the form
$S^1\times Y$, where $Y$ has trivial fundamental group. This
configuration will be of interest in \S\ref{sub-prod-on-prod} and
\S\ref{sec-filters} below.

\subsubsection{$A=\Reals$; $\Sigma=$ join-semilattice theory with zero.}
                        \label{subsub-depend-metric}
We let $\Sigma$ consist of any axioms for semilattice theory,
expressed in terms of a join operation $\vee$, together with a
zero operation for $\vee$. More precisely, the equations that we
shall require in $\Sigma$ are these somewhat weaker equations:
\begin{gather*}
           (x\vee y)\vee y \Wavy x\vee y, \quad
                    y\vee(y\vee x) \Wavy y\vee x,\\
           0\vee x \Wavy  x\Wavy x\vee 0 , \\
                    x\vee x \Wavy x.
\end{gather*}

Here we shall show that each of the equations
\begin{align}   \label{eq-reals-with0-a}
                \lambda_{(A,d)}(\Sigma) &\EQ 0 \\
                \lambda_{(A,d)}(\Sigma) &\EQ \infty
                                  \label{eq-reals-with0-b}
\end{align}
holds for an appropriate choice of a metric $d$ inducing the usual
topology on a space $A$ homeomorphic to $\Reals$. (Moreover,
(\ref{eq-reals-with0-a}) provides us with a further example of $A$
not compatible with $\Sigma$ but for which $\lambda_A(\Sigma)=0$.
It differs from the example in \S\ref{subsub-lim-zero} in that
this $\Sigma$ is consistent. For an example with $A$ compact, see
\S\ref{subsub-lim-zero-consis} below.)

For  (\ref{eq-reals-with0-a}), we take $A=(0,1)$ (homeomorphic to
$\Reals$), with $d$ the usual metric on $(0,1)$. To establish
(\ref{eq-reals-with0-a}), it will suffice to find continuous
operations that satisfy $\Sigma$ within $\varepsilon$, for every
$\varepsilon>0$. This is easily accomplished: we let $\jj$ be the
usual semilattice operation, and take $\o 0$ (the element of our
algebra that is denoted by the constant symbol $0$) to be
$\varepsilon$. The detailed verification is left to the reader.

For  (\ref{eq-reals-with0-b}), we take $A=\Reals$, with $d$ the
usual metric. For a contradiction, suppose that $
\lambda_{(A,d)}(\Sigma) \EQ K < \infty$. By definition, there are
a constant $\o 0$ and a continuous operation $\o\vee$ on $\Reals$
satisfying $\Sigma$ within $K$. Take $a<b\in\Reals$, with $\o0$
between $a$ and $b$ and with $d(a,\o0)=d(b,\o0)=3K$.

Let us suppose, without loss of generality, that $a\jj
b\,\geq\,0$. Now consider the continuous real-valued function
\begin{align*}
         x\GOESTO a\jj x.
\end{align*}
 It follows
from $\Sigma$ that $a\jj a<\o 0$. Thus our function maps $a$ and
$b$ to $a\jj a$ and $a\jj b$, which lie on opposite sides of
$\o0$. By the IVT, there exists $e\in\Reals$ such that $a\jj e \EQ
\o0$. Using $\Sigma$, we have
\begin{align*}
       \o0 \EQ a\jj e \WAVY a\jj(a\jj e)\EQ a\jj\o0.
\end{align*}
Therefore $d(\o0,a\jj\o0)\,\leq\,K$; again invoking $\Sigma$, we
have $d(a\jj\o0,a)\,\leq\,K$ We now use the triangle inequality to
compute
\begin{align*}
        d(\o0,a) \;&\leq\; d(\o0,\,a\jj\o0) \,+\, d(a\jj\o0,\,a)\\
                 \;&\leq\; K\,+\, K \EQ 2K.
\end{align*}
 This contradiction completes the proof of
(\ref{eq-reals-with0-b}).

\subsubsection{A lemma on approximate satisfaction of lattice
equations.}                        \label{subsub-approx-lattice}

In Lemma \ref{lem-approx-lat} which follows, $Y$ is the space of a
connected acyclic one-dimen\-sion\-al simplicial complex, of which
$U=[u_0,u_1]$ is a designated edge. Each edge of $Y$ has a known
length, and then distance $d(A,B)$ is measured by accumulating
these lengths along the unique injective path from $A$ to $B$
(with distance pro-rated along an incomplete edge).

If we arbitrarily choose one ordering of the vertices of $U$, say
$u_0<u_1$, then all of $Y$ can be ordered, as follows. If there is
no injective path containing $A$, $B$, $u_0$ and $u_1$, then $A$
and $B$ are not $<$-related. Otherwise a unique such path exists,
and we can list the four vertices in the order they appear along
the path, making sure to list $u_0$ before $u_1$. For instance, we
might have $u_0ABu_1$, or $Bu_0u_1A$, or a number of other
possibilities. Then we simply take $A\leq B$ (resp. $A\geq B$) to
mean that $A$ appears before (resp. after) $B$ in this listing.
The notation $[A,B]$ denotes, in the usual way, a closed interval
under this ordering.

It is easy to see that {\em if\/ $A\in U$ and\/ $B\in Y$, then
trichotomy holds: $A<B$ or $A=B$ or $A>B$.}

\begin{lemma}           \label{lem-special-acyclic}
If\/ $K$ is a connected subset of\/ $Y$, if\/ $r,s,t\in Y$ with\/
$r\leq s\leq t$, and if\/ $r,t\in K$ and $s\in U$, then $s\in K$.
\end{lemma} \begin{Proof}
By contradiction. If $s\not\in K$, then the two sets
\begin{align*}
       K_0 \EQ \{\, s\in K \,:\, x\,<\,s\,\}\quad\text{and}\quad
          K_1 \EQ \{\, s\in K \,:\, x\,>\,s\,\}
\end{align*}
are disjoint non-empty open sets in $K$ with union $K$.
\end{Proof}

In the context of operation symbols $\wedge$ and $\vee$, the {\em
dual} of any lexical object is formed by interchanging these two
symbols.

{\em Commutative rearrangement\/} is the smallest equivalence relation on
($\meet,\join$)-terms that satisfies the following condition: if $\tau_i$
is a commutative rearrangement of $\sigma_i$  ($i=1,2$), then
$\tau_2\join\tau_1$
is a commutative rearrangement of $\sigma_1\join\sigma_2$ (and similarly
for $\meet$). A {\em commutative rearrangement} of an equation $\sigma_1
\wavy\sigma_2$ is any equation $\tau_1\wavy\tau_2$ where $\tau_i$
is a commutative rearrangement of $\sigma_i$ ($i=1,2$).

In the following lemma, $Y$ and $U$ continue to be as specified at
the start of \S\ref{subsub-approx-lattice}.

\begin{lemma}              \label{lem-approx-lat}
Let $\varepsilon > 0$. Suppose that a topological algebra $\mathbf
Y=(Y,\mm,\jj)$ is given, and that the following equation-set, when
closed under
duals and commutative rearrangements, holds within $\varepsilon$ for
all $x_0,x_1\in U$:
\begin{align}
             x_0\vee x_0 &\;\WAVY\; x_0 \nonumber \\
             x_0\vee (x_0\vee x_1) &\;\WAVY\; x_0\vee x_1
                 \nonumber 
       \\        x_0\vee(x_1\wedge x_0) &\;\WAVY\; x_0.\nonumber
\end{align}
Suppose that\/ $a,b\in U$, with $d(a,b)>\varepsilon$. Then
\begin{itemize}
  \item[(i)] Either $a\jj b\in[a,b]$, or $d(a\jj b,a)\leq\varepsilon$,
         or $d(a\jj b,b)\leq\varepsilon$. The same disjunction
         holds for $a\mm b$.
  \item[(ii)] Either $a\jj b$ is within $2\varepsilon$ of\/ $b$, or
                 $a\mm b$ is within $2\varepsilon$ of\/ $b$.
\end{itemize}
\end{lemma}\begin{Proof}
Proof of (i).
Without loss of generality, we will assume that $a<b$. We need
prove (i) only for $a\jj b$. (The dual of our proof for $a\jj b$
is a proof for $a\mm b$.) If the first alternative holds, we are
done; hence we will assume that it fails, i.e.\ that $a\jj
b\not\in[a,b]$. By trichotomy, this means that $a\jj b > b$ or
$a\jj b < a$. We consider the first of these possibilities; the
second is similar (subject to some commutative rearrangement),
and will be left to the reader.

We are assuming that $d(a,a\jj a)<\varepsilon$, but
$d(a,b)>\varepsilon$. Therefore $a\jj a<b$. On the other hand we
are in the case where $a\jj b > b$. Let us take $K$ to be the
image of the function $x\GOESTO a\jj x$ over the interval $[a,b]$.
Clearly $K$ is connected and contains $a\jj a$ and $a\jj b$ with
$a\jj a \,<\,b\,<\,a\jj b$. We shall apply Lemma\footnote%
{Here Lemma \ref{lem-special-acyclic} replaces the Intermediate
Value Theorem. If we knew that $K\subseteq U$, or even that $K$
lies in some subset of $Y$ homeomorphic to $\Reals$, then we could
use the IVT directly.}%
 \ \ref{lem-special-acyclic} with $r=a\jj a$, $s=b$ and $t=a\jj b$.

By the Lemma, we have $b\in K$. In other words, there exists $e\in
U$ for which $a\jj e = b$. By the approximate satisfaction of
our equations, we have
\begin{align*}
               a\jj b \EQ a\jj(a\jj e) \WAVY a\jj e \EQ b,
\end{align*}
and hence $d(a\jj b,b)<\varepsilon$. This completes the proof of
Part (i).

For Part (ii), we shall assume, by way of contradiction, that
\begin{align}            \label{eq-sortout-lattice}
           d(a\jj b, b)\,\geq\,2\varepsilon;\quad\quad
           d(a\mm b, b)\,\geq\,2\varepsilon.
\end{align}
Therefore, for both $a\jj b$ and $a\mm b$, the first or second
clause of (i) must hold. In particular, $a\jj b$ is either in
$[a,b]$ or $<a$, and the same may be said of $a\mm b$. These
alternatives divide the rest of the proof into three cases.
\vspace{0.2in}

\hspace*{-\parindent}%
{\bf Case 1.} $a\jj b < a$. We know that $b\jj b\wavy b$; which
implies that $a< b\jj b$ (since we are given that
$d(a,b)>\varepsilon$). Therefore, we may apply Lemma
\ref{lem-special-acyclic} to see the existence of $e$ with $e\jj b
= a$. By our equations, we now have
\begin{align*}
         b \Wavy b\mm(e\jj b) = b\mm a,
\end{align*}
which is one alternative of the desired conclusion of Part
(ii).\vspace{0.2in}

\hspace*{-\parindent}%
{\bf Case 2.} $a\mm b < a$. The proof is dual to that of Case 1.
\vspace{0.2in}

\hspace*{-\parindent}%
{\bf Case 3.} $a\jj b$ and $a\mm b$ both lie in $[a,b]$. In this
case, we may apply trichotomy and observe that without loss of
generality, we have
\begin{align}           \label{eq-sortout-lattice-b}
                a\mm b \;\leq\; a\jj b \;<\; b.
\end{align}
By one of our given approximate identities, we have $d(b\mm b,
b)<\varepsilon$; combining this with (\ref{eq-sortout-lattice})
and (\ref{eq-sortout-lattice-b}), we see that
\begin{align*}
               a\mm b \;\leq\; a\jj b \;\leq\; b\mm b.
\end{align*}
Now we have $a\jj b$ between two values of the continuous function
$x\GOESTO x\mm b$. Moreover, $[a,b]\subseteq U$, and hence $a\jj
b\in U$. Another application of Lemma \ref{lem-special-acyclic}
yields $e\in U$ such that $e\mm b\,=\,a\jj b$.

Using the approximate satisfaction of the given equations, we
calculate
\begin{align*}
   b \WAVY b\jj(e\mm b) \EQ b\jj(a\jj b) \WAVY a\jj b.
\end{align*}
Thus $d(b,a\jj b)\,<\,2\varepsilon$, and Part (ii) is proved.
\end{Proof}

The following corollary was first proved by A. D. Wallace
in the nineteen-fifties. See
\cite{faucett, koch-wall-a, koch-wall-b, wallace}.

\begin{corollary}  \label{cor-wallace-etal}
If $([0,1],\mm,\jj)$ is a topological lattice, then either $\mm$
and $\jj$ are ordinary meet and join of real numbers (i.e.\ $x\mm
y$ is the smaller of $x$ and $y$, and $x\jj y$ is the larger), or
dually.
\end{corollary}\begin{Proof}
We apply Lemma \ref{lem-approx-lat} in the special case where $Y$
has only one simplex of dimension $1$; in other words, in the case
where $Y=U=[0,1]$. Considering arbitrary $a,b\in[0,1]$, Part (ii)
of the lemma tells us that  $a\jj b$ is a limit point of the set
$\{a,b\}$. Hence $a\jj b\in\{a,b\}$. Similarly $a\mm b\in\{a,b\}$.

We may assume, without loss of generality, that $0\jj 1=1$.
We now consider two subsets of $(0,1]$, namely
\begin{align*}
         K_0 \EQ \{\, x\in (0,1]\,:\, 0\jj x = 0 \,\}; \quad
         K_1 \EQ \{\, x\in (0,1]\,:\, 0\jj x = x \,\}.
\end{align*}
It follows readily from the continuity of $0\jj x$ that $K_i$ is
open in $(0,1]$ ($i=0,1$), and from what we have already proved
that $(0,1]$ is the disjoint union of these two sets. By
connectedness, one of the two sets is empty. Clearly $1\in K_1$,
so $K_0$ is empty. We have now established that $0\jj x=x$ for all
$x$.

We next consider $b>0$, and define the two sets
\begin{align*}
         K_0 \EQ \{\, x\in [0,b)\,:\, b\jj x = b \,\}; \quad
         K_1 \EQ \{\, x\in [0,b)\,:\, b\jj x = x \,\}.
\end{align*}
As before, the two sets form an open partition of $[0,b)$. By what
we have done before, $0\in K_0$, so $K_0=[0,b)$. In other words we
have proved that if $a<b$, then $a\jj b =b$. Finally, for $a\leq
b$, we use lattice theory to calculate that $a\mm b = a\mm (a\jj
b) =a$. Thus $\mm$ and $\jj$ agree with the usual lattice
operations on $[0,1]$.
\end{Proof}

\subsubsection{$A=Y$, the  triode; $\Sigma=$ lattice theory.}
                                       \label{subsub-Y}

 Let $B,C,D,E$ be four non-collinear points in
the Euclidean plane, with $E$ in the interior of $\Tri{B}{C}{D}$.
Our space $Y$ is defined to be the union of the three (closed)
segments $BE$, $CE$ and $DE$, called {\em legs}, with the topology
inherited from the plane. In fact, in order to give $Y$ a definite
metric $d$, we will further require that $\Tri{B}{C}{D}$ be
equilateral with $E$ at its center, and that each leg have unit
length. We then let $d$ be the metric of the plane, as inherited
by $Y$.

We let $\Sigma$ consist of any axioms for lattice theory
(expressed in terms of $\wedge$ and $\vee$), which includes the
equations of Lemma \ref{lem-approx-lat} on page
\pageref{lem-approx-lat} of \S\ref{subsub-approx-lattice}, and
also the commutative law $x_0\wedge x_1\wavy x_1\wedge x_0$, and
its dual. It was proved by A. D. Wallace in the mid-1950's (see
\cite[{\em Alphabet Theorem}]{wallace} for a statement of
the result) that $Y$ is not compatible with lattice theory. Here
we shall prove the sharper result that
$\lambda_Y(\Sigma)\,\geq\,0.25$.

For a contradiction, suppose that $ \lambda_Y(\Sigma)\,<\,0.25$.
By definition, there are continuous operations $\o\wedge$ and
$\o\vee$ on $Y$ satisfying $\Sigma$ within $\lambda_Y(\Sigma)$,
which is $\,<\,0.25$.

We shall now apply Lemma \ref{lem-approx-lat} of
\S\ref{subsub-approx-lattice} to this $Y$, and to $U$ taken to be,
first $BE$, then  $CE$, then  $DE$.
Conclusion (ii) of that Lemma yields, respectively,
\begin{align*}
             \bigl(\,d(B,B\mm
             E)\,<\,0.5\,\bigr)\quad&\text{or}\quad
                  \bigl(\,d(B,B\jj E)\,<\,0.5\,\bigr);\\[0.05in]
             \bigl(\,d(C,C\mm
             E)\,<\,0.5\,\bigr)\quad&\text{or}\quad
                  \bigl(\,d(C,C\jj E)\,<\,0.5\,\bigr);\\[0.05in]
             \bigl(\,d(D,D\mm
             E)\,<\,0.5\,\bigr)\quad&\text{or}\quad
                  \bigl(\,d(D,D\jj E)\,<\,0.5\,\bigr).
\end{align*}
Clearly either the $\mm$-alternative holds for at least two of these
three propositions, or the $\jj$-alternative holds at least twice.
We will assume, without loss of generality, that the
$\mm$-alternative holds in at least two of the propositions, and in
particular that it holds for legs $BE$ and $CE$. Thus we have
\begin{align}               \label{eq-two-ends-near}
             \bigl(\,d(B,B\mm
             E)\,<\,0.5\,\bigr)\quad&\text{and}\quad
              \bigl(\,d(C,C\mm E)\,<\,0.5\,\bigr).
\end{align}
The proof now divides into two cases, depending on the location of
the point $B\mm C$. Clearly it must lie in one of the three
legs.\vspace{0.05in}

\hspace*{-\parindent}%
{\bf Case 1. $B\mm C\in DE$ or $B\mm C\in BE$.} Consider the
continuous map
\begin{align*}
                   z\GOESTO z\mm C,
\end{align*}
defined for $z\in Y$. Clearly the image of $C$ is $C\mm C$, which
lies within 0.25 of $C$ (by the approximate satisfaction of
$\Sigma$), and hence is a point of leg $CE$. On the other hand, the
image of $B$ is $B\mm C$, which we have assumed to lie either on leg
$DE$ or on leg $BE$. Thus the image of this map is a connected
subset of $Y$ that contains points on at least two legs. Therefore,
the image of our map must contain $E$, which is to say that
\begin{align*}
                   P\mm C \EQ E
\end{align*}
for some $P\in Y$.

Using the equations of $\Sigma$, we calculate
\begin{align*}
              E \EQ P\mm C  \WAVY (P\mm C)\mm C \EQ E\mm C
                 \WAVY C\mm E,
\end{align*}
which implies that
\begin{align*}
           d(E,C\mm E) \;<\;0.5.
\end{align*}
Combined with (\ref{eq-two-ends-near}), this last inequality yields
$d(E,C)\,<\,1$, in contradiction to our specification of $CE$ as a
segment of length $1$. This contradiction completes the proof for
Case 1.
 \vspace{0.05in}

\hspace*{-\parindent}%
{\bf Case 2. $B\mm C\in CE$.} This situation is symmetric
to $B\mm C\in BE$, which was considered in Case 1. (Interchanging
the letters $B$ and $C$ in that proof, while sending $z$
to $B\mm z$, yields a proof here.) This
completes the proof that $\lambda_Y(\Sigma)\,\geq\,0.25$.

\subsubsection{$A=Y$ with a new metric.}
\label{subsub-Y-new-metric} Let $\Sigma$ be the axiom-set
(expressing part of lattice theory) that was presented in
\S\ref{subsub-Y}, and let $Y$ be the space defined there (a union
of three segments). We proved that for a certain metric of
diameter $\sqrt{3}$ the space $Y$ satisfies $\lambda_Y(\Sigma)\geq
0.25$. In other words, $\lambda_Y(\Sigma)\geq
0.25/\sqrt{3}=0.14433\ldots$ for a certain natural metric of
diameter $1$.

Here we shall show that for each real $\varepsilon>0$, there
exists a diameter-$1$ metric $d$ for the topology of $Y$ with
$\lambda_{(Y,d)}(\Sigma^{\star})<\varepsilon$.

The metric $d$ is obtained very simply by expressing $Y$---as
before---as the union of unit-length segments $EB$, $EC$ and $ED$
in the plane; this time with $B$, $C$ and $D$ chosen to lie within
$\varepsilon$ of one another, and with no two of them collinear with
$E$. The metric $d$ is then simply the
restriction to $Y$ of the metric in the plane. Clearly, for
$\varepsilon<1$, the metric space $(Y,d)$ has diameter $1$, and so
it will be enough to prove that $\lambda_{(Y,d)}<\varepsilon$.

We first define
\begin{align*}
                \psi\FROM Y \TO Y
\end{align*}
as follows: for any point $X\in Y$, $\psi(X)$ is the unique point
on segment $EB$ that satisfies $d(E,\psi(X))=d(E,X)$. Clearly
$\psi$ is continuous and $\psi$ is the identity function on $EB$.

The reader may check that we now have the hypotheses of Theorem
\ref{thm-inequal-retract} of \S\ref{sub-inequal-retract} (with $Y$
for $A$, $EB$ for $B$, and $\varepsilon$ for $K$). In this
context, the conclusion of the theorem is
\begin{align*}
        \lambda_{(Y,d)}(\Lambda)\;\leq\;\lambda_{(EB,d)}(\Lambda)
                     \,+\,\varepsilon.
\end{align*}
Moreover $EB$ is topologically a closed segment, which is
compatible with lattice theory; hence
$\lambda_{(EB,d)}(\Lambda)=0$. We now have the desired conclusion
that $\lambda_{(Y,d)}(\Sigma)<\varepsilon$.

 \subsubsection{$A=[0,1]$, $\Sigma_c$ an extension of lattice theory}
 \label{subsub-lim-zero-consis}
 Let $\Sigma_c$ be the following set of equations:
\begin{align*}
         &\text{axioms of lattice theory in $\wedge$, $\vee$}\\
         &0\wedge x \WAVY 0\\
         &1\vee x \WAVY 1\\
         &F(0)\WAVY 0\quad\quad\quad\quad   F(1)\WAVY 1\\
         &G(F(x))\WAVY 1\quad\quad\quad\quad  G(a)\WAVY 0.
\end{align*}

We will see that $\Sigma_c$ is consistent, that $\Sigma_c$ is not
compatible with $[0,1]$, and that $\lambda_{[0,1]}(\Sigma_c)=0$.

Let us first see that $\Sigma_c$ is not compatible with the
compact interval $[0,1]$. Looking for a contradiction, we will
suppose that $([0,1],\o\wedge,\o\vee,\o0,\o1,\o F,\o G,\o
a)\models\Sigma_c$, with $\o\wedge$, $\o\vee$, $\o F$ and $\o G$
continuous. Applying Corollary \ref{cor-wallace-etal}, we may
assume, without loss of generality, that $\mm$ and $\jj$ are
ordinary meet and join on the interval $[0,1]$. Obviously then,
$\o0=0$ and $\o1=1$. By the continuity of $\o F$. its range is
connected, and clearly the range contains $0$ and $1$. Therefore,
$\o F$ maps $[0,1]$ onto itself. Thus the penultimate equation,
$G(F(x))\wavy 1$, tells us that $\o G$ is a constant function with
value $1$. The final equation, however, yields $\o G(\o a)=0$;
hence $\o G$ is not a constant function. This contradiction
concludes the proof that $[0,1]$ is not compatible with
$\Sigma_c$.

We next prove that
\begin{align}      \label{eq-arb-d-zero}
\lambda_{([0,1],d)}(\Sigma_c)=0 ,
\end{align}
where $d$ is the usual metric on $[0,1]$. In order to prove
(\ref{eq-arb-d-zero}), for each $\varepsilon>0$ we shall construct
continuous $[0,1]$-operations modeling $\Sigma_c$ to within
$\varepsilon$.

We let $\jj$ and $\mm$ be ordinary join and meet on $[0,1]$. These
are well known to be continuous, and certainly they satisfy the
equations of lattice theory. We define $\o 0$ to be a point of
$[0,1]$ other than $0$ (the usual zero of $\Reals$), such that
$d(\o 0,0)<\varepsilon$. The definition of $\o1$ is dual to that
of $\o0$. Let us check the equation $0\wedge
x\wavy 0$. Its $\varepsilon$-interpretation in this model is that
$d(\o0\mm x,\o0)<\varepsilon$ for all $x\in[0,1]$. For
$x\in[\o0,1]$, we in fact have $\o0\mm x=\o0$, so we need only
consider $x\in[0,\o0]$. Since $d(\o 0,0)<\varepsilon$, we clearly
have
\[
      d(\o0\mm x,\o0) \EQ d(x,\o0) \;\leq\;
                \varepsilon,
\]
as desired. A similar calculation yields $d(\o1\jj x,\o1)<
\varepsilon$.

We next select a continuous function $\o F\FROM[0,1] \TO[\o0,\o1]$,
such that $\o F(\o0)=\o0$ and $\o F(\o1)=\o1$. Now all equations are
true within $\varepsilon$, except those in the last line, involving
$G$. We define $\o a$ to be $0$ (i.e.\ the real zero of $\Reals$),
and define $\o G\FROM[0,1]\TO[0,1]$ to be a continuous function
(piecewise linear will do) such that $\o G$ is constantly $\o1$ on
$[\o0,\o1]$, and $\o G(0)=0$. The equations in the last line now
hold exactly.

This completes the proof of (\ref{eq-arb-d-zero}). In fact, it is
now evident from Lemma \ref{lemma-comp-two-metrics} that
(\ref{eq-arb-d-zero}) holds for any metric defining the usual
topology on $[0,1]$. (We revisit this fact in
\S\ref{subsub-incom-delta+} below.)

To see that $\Sigma_c$ is consistent, take any bounded lattice $L$
of three or more elements, and an $\o F\FROM L\TO L$ that is not
onto and satisfies $\o F(\o0)=\o0$ and  $\o F(\o1)=\o1$. Select
$\o a$ from the complement of the range of $\o F$, and then define
$\o G$ in the obvious way. Thus we have a model of $\Sigma_c$ with
more than one element.

\subsection{The theory of an injective binary operation.}

\subsubsection{$A=[0,1]$; $\Sigma=$ ``injective
binary.''}
 \label{subsub-inf-dim}
Consider $\Sigma$ consisting of the two equations
\begin{align*}          
                F_0(G(x_0,x_1)) &\WAVY x_0 \\
                F_1(G(x_0,x_1)) &\WAVY x_1. 
\end{align*}
They imply, among other things, that in any topological model $\o
G$ must be a one-one continuous binary operation. Euclidean spaces
and Euclidean-like spaces (such as simplicial complexes)
of non-zero finite dimension do not have such operations, hence
are not compatible with $\Sigma$. In \S\ref{subsub-inf-dim} we
will be concerned with $A=[0,1]$. Spaces $A$ of higher dimension
will be considered in \S\ref{subsub-square-no-inj} and
\S\ref{sub-delta-square-embed-square}.

For $A=[0,1]$ we will show that $\lambda_{(A,d)}(\Sigma)\geq
0.5\cdot\text{diam}(A)$, for any metric $d$ that defines the usual
topology on $A$. (And of course, then
$\lambda_{(A,d)}(\Sigma^{\star})\geq 0.5\cdot\text{diam}(A)$.) It
thus follows from \S\ref{subsub-lambda} that if $\rho$ is the
usual Euclidean metric, then $\lambda_{(A,\rho)}(\Sigma)=
\lambda_{(A,\rho)}(\Sigma^{\star})= 0.5\cdot\text{diam}(A)$.

For the proof, we let $d$ be an appropriate metric, and we
consider continuous operations $\o F_0$, $\o F_1$ and $\o G$ that
satisfy $\Sigma$ within $K$. We shall prove that
$K\geq0.5\cdot\text{diam}(A)$.

Since $[0,1]$ is compact, there exist $a_0,a_1\in A$ with
$d(a_0,a_1)=\text{diam}(A)$. For flexibility of notation, we
take two such pairs: $d(a_0,a_1)=d(b_0,b_1)=\text{diam}(A)$.
Considering the four real numbers
\begin{align*}
     \o G(a_0,b_0),\;\;\o G(a_1,b_0),\;\;\o G(a_0,b_1),\;\;
                 \o G(a_1,b_1),
\end{align*}
we may assume, without loss of generality, that the smallest
among them is $\o G(a_0,b_0)$. Again without loss of
generality, we may assume that $\o G(a_1,b_0))\,\leq\,\o
G(a_0,b_1).$ In other words, we have
\begin{align*}
             \o G(a_0,b_0)\;\leq\; \o G(a_1,b_0)\;\leq\;\o
             G(a_0,b_1).
\end{align*}
Thus, along the segment $\overline {(a_0,b_0)(a_0,b_1)}$ in the
square $[0,1]^2$, the continuous function $\o G$ takes values that
are above and below the value $\o G(a_1,b_0)$. By the IVT, there
exists $e\in [0,1]$ such that
\begin{align*}
                  \o G(a_0,e)\EQ \o G(a_1,b_0).
\end{align*}
From $\Sigma$ we now calculate
\begin{align*}
            a_0\WAVY \o F_0(\o G(a_0,e)) \EQ
             \o F_0(\o G(a_1,b_0)) \WAVY a_1.
\end{align*}
From this we immediately see the desired conclusion that
$\text{diam}(A)\leq 2K$.

One may notice, incidentally, that this argument does not require that the
 operations $\o F_i$ be continuous.

\subsubsection{$A=[0,1]^2$; $\Sigma=$ ``injective
binary.''}
  \label{subsub-square-no-inj}
 So, we consider the same $\o F_0$, $\o F_1$ and
$\o G$, as they might be realized on the square $B=[0,1]^2$ (with the
usual Euclidean metric).

Let $\Sigma$ be as in \S\ref{subsub-inf-dim}. We will show that
$\lambda_B(\Sigma) = 0.5$. We first show that it is $\leq0.5$.
To this end, consider the following definitions of operations
on $B=[0,1]^2$:
\begin{align*}
            \o G((a_0,a_1),(b_0,b_1)) &\EQ (a_0,b_1)\\
               \o F_0((a_0,a_1)) &\EQ (a_0,0.5)\\
               \o F_1((a_0,a_1)) &\EQ (0.5,a_1).
\end{align*}
Now $\o F_0(\o G((a_0,a_1),(b_0,b_1))) = (a_0,0.5)$, so its
distance from $(a_0,a_1)$ is $d(a_1, 0.5) \allowbreak\leq 0.5$.
A similar calculation applies to $\o F_1(\o
G((a_0,a_1),(b_0,b_1)))$, and so both equations of $\Sigma$ are
seen to hold within $0.5$.

We next show that $\lambda_B(\Sigma) \geq 0.5$. To this end,
suppose that its true value is $K$, and now consider any
operations $\o G, \o F_0,  \o F_1$ that satisfy $\Sigma$ within
$K$ on $B$, and with $\o G$ continuous. Consider the action of $\o
G$ on the boundary of $B^2$, which is homeomorphic to $[0,1]^4$.
The boundary of this space is a three-sphere $S^3$. By the
Borsuk-Ulam Theorem\footnote{%
Proved by K. Borsuk in 1933---see \cite{borsuk}. See Steinlein
\cite{steinlein} or Matou\v sek \cite{matousek} for a thorough
discussion
of this theorem and its applications in modern mathematics.}, %
$\o G$ takes on the same value at two antipodal points. Without
loss of generality, two antipodal points have the form
$((0,x_1)(y_0,y_1))$ and $((1,u_1)(v_0,v_1))$. So now, by the
triangle inequality,
\begin{align*}
      1 &\EQ  d((0,0),(1,0)) \,\leq\, d((0,x_1),(1,u_1))\\
        &\;\leq\; d((0,x_1),\o F_0 \o G((0,x_1),(y_0,y_1)))
             \,+\,  d((1,u_1),\o F_0 \o G((1,u_1)(v_0,v_1)))\\
        &\;\leq\; K\,+\, K.
\end{align*}

We remark that if the metric is scaled to make the diameter equal
to $1$, then $\lambda_B(\Sigma) = 0.5/\sqrt{2} \approx .3536.$

\subsubsection{$A=[0,1]^2$, with a new metric.}
\label{sub-delta-square-embed-square}
         Let $\Sigma$ be defined as in \S\ref{subsub-inf-dim},
and take $B=[0,1]^2$, in its usual topology. Here we will show how
to define a unit-diameter metric $d$ for the usual topology on $B$
such that $\lambda_{(B,d)}(\Sigma)$ is arbitrarily small. More
precisely, given real $\varepsilon>0$, we construct a
unit-diameter metric $d$ for the topology of $B$, such that
$\lambda_{(B,d)}(\Sigma) \leq\varepsilon$.

Let us replace $B$ by the set $[0,\sqrt{1-\varepsilon^2}]
\times [0,\varepsilon]$, while taking $d$ to be the Euclidean
metric, restricted to this set. Clearly this $B$ is
homeomorphic to our original $B$, and moreover $B$ has unit
diameter. For an upper estimate on $\lambda_{(B,d)}(\Sigma)$,
we define these three operations on $B$:
\begin{align*}
        \o G((a_0,a_1),(b_0,b_1)) &\EQ (a_0,b_0)\\
        \o F_0(a_0,a_1) &\EQ (a_0,0)\\
        \o F_1(a_0,a_1) &\EQ (a_1,0).
\end{align*}
We now calculate
\begin{align*}
           d((a_0,a_1),\o F_0(\o G((a_0,a_1),(b_0,b_1)))) &\EQ
       d((a_0,a_1),(a_0,0))\\
          &\;\leq\; \varepsilon.
\end{align*}
and
\begin{align*}
           d((b_0,b_1),\o F_1(\o G((a_0,a_1),(b_0,b_1)))) &\EQ
       d((b_0,b_1),(b_0,0))\\
          &\;\leq\; \varepsilon.
\end{align*}
Thus the equations of $\Sigma$ hold within $\varepsilon$, and
our estimate on $\lambda_{(B,d)}(\Sigma)$ is established.

\subsubsection{Generalization of \S\ref{subsub-inf-dim}
and \S\ref{sub-delta-square-embed-square}}
          \label{sub-genzations-delta}
For $m>n$, let $\Sigma_{m,n}$ denote this set of equational axioms:
\begin{align}   \label{eq-first-of-biject}
     F_1(G_1(x_1,\cdots,x_m),\cdots,G_n(x_1,\cdots,x_m))
                                        &\Wavy x_1 \\
      \vdots\quad\quad\quad\quad\quad\quad & \nonumber \\
     F_m(G_1(x_1,\cdots,x_m),\cdots,G_n(x_1,\cdots,x_m))
                                        &\Wavy x_m. \nonumber
\end{align}
If we let $F$ and $G$ stand for the appropriate tuples of
function symbols, and let $x$ stand for a tuple of variables
$x_i$, then the equations may be symbolically abbreviated as
\begin{align*}
                    F(G(x)) \Wavy x.
\end{align*}

To model $\Sigma_{m,n}$ on a space $A$ is to find an $n$-tuple
of continuous functions $\o G_i\FROM A^m\TO A$, in other words
a continuous function $\o G\FROM A^m\TO A^n$, and also a
continuous function $\o F\FROM A^n\TO A^m$, such that the
composite mapping
\begin{align} \label{eq-inj-power}
             A^m \stackrel{\o G}{\TO} A^n \stackrel{\o F}{\TO} A^m
\end{align}
is the identity map on $A^m$. This is of course generally
impossible for $m>n$ and for a space $A$ of dimension $d$ with
$0<d<\infty$; in other words, for such spaces $A$, the
compatibility relation $A\models\Sigma_{m,n}$ does not hold.
Nevertheless, the associated metrical invariant $\lambda$ may or
may not be close to zero, as we see in Theorem \ref{th-m=nk} just
below. The theorem was already proved for $m=2$, $n=1$ and
$k=1,2$, in \S\ref{subsub-inf-dim}, \S\ref{subsub-square-no-inj}
and \S\ref{sub-delta-square-embed-square} above.

\begin{theorem} \label{th-m=nk}
Let $B=[0,1]^k$, with $k \in\{ 1, 2, 3, \ldots\}$. Then
\begin{itemize}
\item[(i)] $\lambda_{(B,d)}(\Sigma_{m,n})
            \leq0.5$ for $d$ the Euclidean metric (scaled to diameter $1$).
\item[(ii)] If $m\geq nk+1$, then $\lambda_{(B,d)}(\Sigma_{m,n})
            \geq0.5$ for $d$ any metric of diameter $1$.
\item[(iii)] If $m\leq nk$, then for every $\varepsilon>0$ there is
            a diameter-$1$ metric $d$ on $B$ such that
            $\lambda_{(B,d)}(\Sigma_{m,n})\leq\varepsilon$.
\end{itemize}
\end{theorem}\begin{Proof}
We first observe that, in the Euclidean metric,
$\lambda_B(\Sigma_{m,n})\leq 0.5$, for any $m,n,k\geq 1$. As was
pointed out in \S\ref{subsub-lambda}, we may define all operations
to be constant at the centroid, thereby obtaining the radius (in
this case half the diameter) as an upper estimate for $\lambda_B$.
This observation takes care of point (i).

For (ii), we suppose that $d$ is a metric for $B$ with diameter
$1$, and that $\o F$ and $\o G$ are continuous functions as in
(\ref{eq-inj-power}) (with $B$ in place of $A$). There exist
$a,b\in B$, with $d(a,b)=1$. By arc connectedness, there is an arc
$I$ in $B$ with endpoints $a$ and $b$. Let $S^{m-1}\subseteq B^m$
be the subset of $B^m$ that is the boundary of $I^m$. It is a
sphere of dimension $m-1$. We now have a restriction of $\o G$
which maps as follows:
\begin{align*}
           S^{m-1}\;\subseteq \;B^m \;\stackrel{\o G}{\TO}
                 \;B^n\; \ISO\; [0,1]^{nk}.
\end{align*}
Now since $m-1\geq nk$, the Borsuk-Ulam Theorem applies, and so
there exist antipodal points $x,y\in S^{m-1}$ with $\o G(x)=\o
G(y)$. Without loss of generality, $x = (a,x_2, \ldots)$ and $y =
(b,x_2, \ldots)$. Thus
\begin{align*}
    1\EQ d(a,b) &\;\leq\;  d(a,\o F_1(\o G(x))) \,+\,
                      d(b,\o F_1(\o G(y))).
\end{align*}
Thus one of these two summands must be $\geq0.5$; without loss
of generality,
\begin{align*}
        d(x_1,\o F_1(\o G(x))) \EQ
          d(a,\o F_1(\o G(x))) \;\geq\; 0.5.
\end{align*}
In other words, the two sides of Equation
(\ref{eq-first-of-biject}) must be distant by $0.5$. This
establishes our claim that $\lambda_{(B,d)}(\Sigma_{m,n})
            \geq0.5$.

Finally, we consider assertion (iii). Given arbitrary
$\varepsilon>0$, we give $B=[0,1]^k$ the metric
\begin{align*}
         d(\mathbf x,\mathbf y) \EQ \sup(|x_1-y_1|,
      \,\varepsilon|x_2-y_2|,
         \,  \cdots,\,\varepsilon |x_{k}-y_k|).
\end{align*}
Clearly $B$ has diameter $1$ in this metric, and retains its
original topology.

Further, define $\o F$ and $\o G$ (as in (\ref{eq-inj-power})) via
\begin{align}
   \o G(a_1,\cdots,a_m) &\EQ   \label{eq-zero-biject-G}
   ((a^1_1,\cdots,a^1_k),(a^1_{k+1},\cdots,a^1_{2k}),\;\;\cdots\;)
     \;\in\;([0,1]^k)^n\\
      \o F(b_1,\cdots,b_n) &\EQ ((b_1^1,0,\cdots,0),
      (b_1^2,0,\cdots,0),\;\;\cdots\;)\;\in\;([0,1]^k)^m
                  \label{eq-zero-biject-F}
\end{align}
(Our notation is that $b^v_u$ is the $v$-th component of an
element $b_u$ of $B^m$ (or $B^n$, as needed).) For Equation
(\ref{eq-zero-biject-G}), if $m=nk$ the right-hand is well
defined. If $m<nk$, the formula calls for $a^1_j$ with $m<j\leq
nk$, which is not defined. In this case, we simply use $a^1_j=0$.

As for Equation (\ref{eq-zero-biject-F}), the intent is to
cycle through all co-ordinates of all the $b_i$. In other
words,
\begin{align*}
   \o F(b_1,\cdots,b_n) &\EQ (c_1,\cdots,c_m)\;\in\; B^m,
   \;\;\;\text{where}\\
         \;\;\;c_{k(u-1)+v} &\EQ
         (b_u^v,0,\cdots,0)\in[0,1]^k\EQ B,
\end{align*}
for $u=1,\ldots n$ and $v=1,\ldots,k$.

The continuity of the operations $\o G$ and $\o F$ is evident.
To check the approximate satisfaction of the equations
$\Sigma_{m,n}$, let us evaluate $\o F_j(\o G(a_1,\cdots,a_m))$,
where $j=k(u-1)+v$, for $u=1,\ldots n$ and $v=1,\ldots,k$.
(Since $m\leq nk$, every $j$ is representable in this form.)
Now
\begin{align*}
   \o F_j(\o G(a_1,\cdots,a_m)) &\EQ \bigl(
            \bigl[\o G(a_1,\cdots,a_m)\bigr]^v_u
            ,0,\cdots,0 \bigr),
\end{align*}
where the subscript $u$ indicates the $u^{\text{\tiny th}}$
vector in the vector of vectors that appears on the right-hand
side of Equation (\ref{eq-zero-biject-G}), namely
$(a_{k(u-1)+1}^1,a_{k(u-1)+2}^1,\allowbreak\cdots)$. Moreover
the superscript $v$ refers to the $v^{\text{\tiny th}}$
component of that vector, namely $a_{k(u-1)+2}^1$. In other
words,
\begin{align*}
   \o F_j(\o G(a_1,\cdots,a_m)) &\EQ \bigl(
              a_{k(u-1)+v}^1
            ,0,\cdots,0 \bigr)  \EQ
            \bigl(
              a_j^1
            ,0,\cdots,0 \bigr).
\end{align*}
As for approximate satisfaction, if we let $\varepsilon_i$ denote
$1$ for $i=1$ and $\varepsilon$ for $i>1$, we now have
\begin{align*}
        d(a_j,\o F_j(\o G(a_1,\cdots,a_m))) &\EQ
                \sup_{i=1\cdots m} \;\varepsilon_i\,\bigl|a^i_j
                  -\o F_j(\o G(a_1,\cdots,a_m))^i\bigr|\\
                  &\EQ \sup\,\bigl(|a^1_j-a^1_j|,\,\varepsilon|a^2_j-0|,
                            \cdots,\,\varepsilon|a^k_j-0|\bigr)\\
                  &\EQ \sup\,\bigl(0,\,\varepsilon|a^2_j|,
                           \cdots,\,\varepsilon|a^k_j|\bigr)\\
                  &\;\leq\;\varepsilon.
\end{align*}
\end{Proof}

In case $k=1$, the above proof, while of course correct, supplies
more. In that case, the operations $\o F$ and $\o G$ defined in
(\ref{eq-zero-biject-G}) and (\ref{eq-zero-biject-F}) model
$\Sigma_{m,n}$ exactly on $B=I$.

\subsection{$[n]$-th power varieties.} \label{sub-nth-power}
 \subsubsection{The definition of $\Sigma^{[n]}$.}
               \label{subsub-nth-def}
 In
\S\ref{sub-nth-power} we begin with a similarity type (i.e.\
list of operation symbols) that does not contain $d$ or $g$.
Thus an equation-set $\Sigma$ under consideration does not
mention $d$ or $g$, and its models have no operations $\o d$ or
$\o g$. The construction of $\Sigma^{[n]}$ presented below
enlarges $\Sigma$ by adding some new equations involving the
given operation symbols and new operations symbols $d$ and $g$.
Throughout \S\ref{sub-nth-power} we adopt the convention that
$d$ and $g$ appear only in the guise of having been expressly
added in the formation of $\Sigma^{[n]}$. Obviously, this
convention imposes no essential restriction on the generality
of what is written here.

So let $\Sigma$ be a set of such equations (where the type is
implicitly determined by the operation symbols appearing in
$\Sigma$). By the {\em $[n]$-th power equations of $\Sigma$} we
mean the set $\Sigma^{[n]}$ that is formed by adjoining the
following equations to $\Sigma$:
\begin{align}
            g^n(x) &\Wavy x  \label{eq-[n]-a}\\
           d(x,\cdots,x) &\Wavy x \\
        d(g(x_1),\cdots,g(x_n)) &\Wavy
                       g(d(x_2,\cdots,x_n,x_1))\label{eq-[n]-c} \\
   d(d(x_{11},\cdots,x_{1n}),\cdots,d&(x_{n1},\cdots,x_{nn}))
                          \nonumber\\
      &\Wavy d(x_{11},\cdots,x_{nn})\label{eq-[n]-d}  \\
      F(g(x_1),g(x_2),\cdots) &\Wavy g(F(x_1,x_2,\cdots))
                               \label{eq-[n]-e}   \\
      F(d(x_{11},\cdots,x_{1n}),
                          d(x_{21},&\cdots,x_{2n}), \cdots)
                                   \nonumber\\
                      &\Wavy
         d(F(x_{11},x_{21},\cdots), \cdots,
         F(x_{1n},x_{2n},\cdots)),  \label{eq-scheme-nth}
\end{align}
where (\ref{eq-[n]-e}) and (\ref{eq-scheme-nth}) are really
schemes of equations, one for each operation $F$ of $\Sigma$.
The construction of $\Sigma^{[n]}$ first appears in McKenzie
\cite{mckenzie-cubes}. Further exposition and application of
the theory appears in Taylor \cite{wtaylor-fsv,wtaylor-sae} and
Garcia and Taylor \cite{ogwt-mem}.

The essential fact about $\Sigma^{[n]}$ [{\em op.\ cit.}] is that,
within isomorphism, eacb of its (topological) models is the $n$-th
power of a (topological) model of $\Sigma$, where the new
operations $\o d$ and $\o g$ operate on the $n$-th power by the
shuffling of co-ordinates:
\begin{align}
                     \label{eq-d-bar}
 \o d(\la\alpha_{11},\cdots,\alpha_{1n}\ra, \cdots
                         ,\la\alpha_{n1},\cdots,\alpha_{nn}\ra)
             &\EQ \la\alpha_{11},\cdots,\alpha_{nn}\ra\\
      \o g(\la\alpha_1, \cdots,\alpha_n\ra) &\EQ \la\alpha_n,
                                   \alpha_1,\cdots,\alpha_{n-1}\ra.
                                     \label{eq-g-bar}
\end{align}
If the original $\Sigma$-algebra is denoted $\mathbf A$, then
the $n$-th power model so described is denoted $\mathbf
A^{[n]}$.

\subsubsection{$\mbox{Sets}^{[n]}$.} Here we consider Equations
(\ref{eq-[n]-a}--\ref{eq-scheme-nth}) in the special case where
$\Sigma$ has no operation symbols and no equations (the so-called
variety of sets). This essentially means that Equations
(\ref{eq-[n]-e}--\ref{eq-scheme-nth}) do not appear, and so we
will be considering only Equations
(\ref{eq-[n]-a}--\ref{eq-[n]-d}). In this context, a (topological)
model of $\Sigma^{[n]}$ is completely determined by Equations
(\ref{eq-d-bar}--\ref{eq-g-bar}). It follows from
\S\ref{subsub-nth-def} that
\begin{theorem}
$A\models\mbox{Sets}^n$ iff there exists a space $B$ such that
$A$ is homeomorphic to $B^n$.~\rule[-2.2mm]{1.5mm}{3.7mm}
\end{theorem}

We skip the proof of Corollary \ref{cor-sat-n|m}, which
involves a study of when $[0,1]^m$ is homeomorphic to an $n$-th
power. In \S\ref{subsub-approx-cubes} immediately below, we
shall implicitly include a proof for the case $m=n+1$.

\begin{corollary} \label{cor-sat-n|m}
$[0,1]^m\models{Sets}^{[n]}$ iff
$n|m$.~\rule[-2.2mm]{1.5mm}{3.7mm}
\end{corollary}

\subsubsection{$A=[0,1]^m; \Sigma=\mbox{Sets}^{[n]}$.}
          \label{subsub-approx-cubes}

The proof that follows is valid for all $k\geq n$, but the
statement of the theorem ignores a sharper conclusion that can be
drawn when $k=n$, or indeed when $k$ is any multiple of $n$. (By
Corollary \ref{cor-sat-n|m}, $B\models\Sigma$ for any such $k$.)

\begin{theorem}               \label{th-nthpow-on-inttok}
If\/ $B=[0,1]^k$ and $\Sigma=\mbox{Sets}^{[n]}$, with\/ $k\geq n$,
then for each\/ $\varepsilon>0$ there exists a diameter-$1$ metric
$d$ for $B$ with\/
$\lambda_{(B,d)}(\Sigma^{\star})\leq\varepsilon$.
\end{theorem}\begin{Proof}
We may assume that $\varepsilon\leq 1$. Let $d$ be defined by
\begin{align*}
         d(\mathbf x,\mathbf y) \EQ \sup(|x_1-y_1|,\cdots,|x_n-y_n|,
         \,\varepsilon |x_{n+1}-y_{n+1}|,
                      \cdots,\,\varepsilon |x_{k}-y_k|).
\end{align*}
Now define $\psi\FROM B\TO B$ via
\begin{align*}
       \psi(\mathbf x) \EQ (x_1,\cdots,x_n,0.\cdots,0).
\end{align*}
We omit the easy proofs that the image of $\psi$ is a subset $E$
homeomorphic to $[0,1]^n$, that $\psi\upharpoonright
E=\text{identity}$, and that $d(\mathbf x,\psi(\mathbf
x))\leq\varepsilon$ for all $\mathbf x$.

We are now in a position to apply Theorem
\ref{thm-inequal-retract} of \S\ref{sub-inequal-retract} (with $B$
for $A$, $E$ for $B$, and $\varepsilon$ for $K$). The conclusion
is that
\begin{align*}
             \lambda_{(B,d)}(\Sigma^{\star})\;\leq\;
             \lambda_{(E,d)}(\Sigma^{\star})
                     \,+\,\varepsilon.
\end{align*}
But $(E,d)\models\Sigma^{\star}$, by Corollary \ref{cor-sat-n|m},
and hence we have the desired result that
$\lambda_{(B,d)}(\Sigma^{\star})\leq\varepsilon$.
 \end{Proof}

To state a result in the opposite direction (a non-zero lower
estimate on certain values of $\lambda$ for $n=2$ and $k=1$), we
need to modify $\Sigma$ slightly, to include those consequences of
Equations (\ref{eq-[n]-a}--\ref{eq-[n]-d}) (for $n=2$) that we
will be using in the proof of Theorem \ref{th-low-est-squares}.
(Recall from \S\ref{sub-dep-deduct} that $\lambda$ is not
independent of equational deductions. Here we know the result only
for these consequences of $\Sigma$, and not for the original
$\Sigma$.) We consider the equations
\begin{align}    \label{eq-square-deriv-a}
    x &\WAVY g( g(x,u), g(v,g(w,x)))\\
         x &\WAVY g(g(x,u), d(g(g(d(x),w),v)) ).
         \label{eq-square-deriv-b}
\end{align}
If $\sigma$ is a term in a single binary operation $g$, then its
{\em opposite} is the term obtained by recursively replacing
$g(x,y)$ by $g(y,x)$. The opposite of an equation
$\sigma\wavy\tau$ simply pairs the opposite of $\sigma$ with the
opposite of $\tau$.

\begin{theorem}            \label{th-low-est-squares}
   Let $B=[0,1]$ and let $\Sigma=$ be defined by Equations
   (\ref{eq-square-deriv-a}--\ref{eq-square-deriv-b}) and
their opposite equations. ($\Sigma$ is thus a subset of the theory
of $\mbox{Sets}^{[2]}$.) In the Euclidean metric (scaled to
diameter 1) $\lambda_B(\Sigma)=0.5$, and moreover
$\lambda_B(\Sigma)\geq0.5$ in every diameter-$1$ metric.
\end{theorem}\begin{Proof}
The inequality $\lambda_B(\Sigma^+)\leq0.5$ in the Euclidean case
follows as at the start of the proof of Theorem \ref{th-m=nk} on page
\pageref{th-m=nk}. For
$\lambda_B(\Sigma^+)\geq0.5$, we suppose that $\rho$ is a metric
$\rho$ for $B$ with diameter $1$, and that $\o g$ and $\o d$ are
continuous operations on $B$, binary and unary, respectively. For
the desired inequality, we will assume that $(B,\o g,\o d)$
satisfies Equations
(\ref{eq-square-deriv-a}--\ref{eq-square-deriv-b}) and their
opposites within $K\in\Reals$, and then prove that $K\geq 0.5$.

By compactness, there exist $a,b\in B$ with $\rho(a,b)=1$. Let
$I$ be an arc in $B$ joining $a$ to $b$. Let $S^1$ denote the
boundary of $I^2\subseteq B^2$. By the Borsuk-Ulam Theorem, $\o
g$ takes on equal values at antipodal points of $S^1$. Without
loss of generality, $\o g(a,c)=\o g(b,d)$ for some $c,d\in B$.

We now wish to establish the existence of $e,f\in B$ such that
either
\begin{align}         \label{eq-second-BU-alts}
           \o g(e,\o g(\o d(a),a)) \EQ \o g(f,\o g(\o d(b),b))
           \;\;\;\mbox{or}\;\;\;
           \o g(\o g(\o d(a),a),e)) \EQ \o g(\o g(\o d(b),b),f).
\end{align}
Clearly such $e,f$ exist if $\o g(\o d(a),a)=\o g(\o d(b),b)$;
thus we may assume that these are distinct elements of $B$.
Thus there is an arc $J$ joining $\o g(\o d(a),a)$ and $\o g(\o
d(b),b)$.  As before, two antipodal boundary points of $J^2$
take equal values under $\o g$; one easily sees that this
observation is tantamount to (\ref{eq-second-BU-alts}). To
complete the proof of Theorem \ref{th-low-est-squares}, we
consider two cases, which correspond to the two alternatives of
(\ref{eq-second-BU-alts}):\vspace{0.1in}

\hspace*{-\parindent}%
{\bf Case 1.} $\o g(a,c)=\o g(b,d)$ and $\o g(e,\o g(\o
d(a),a)) = \o g(f,\o g(\o d(b),b))$. In $B$, we may define
\begin{align*}
         r &\EQ \o g(\o g(a,c),\o g(e,\o g(\o d(a),a))) \\
           &\EQ \o g(\o g(b,d),\o g(f,\o g(\o d(b),b)))
\end{align*}
(where the second equality clearly comes from the assumptions
of Case 1). If now we substitute into Equation
(\ref{eq-square-deriv-a})---taking $a$ for $x$, $c$ for $u$,
$e$ for $v$, and $d(a)$ for $w$---clearly the right-hand side
of (\ref{eq-square-deriv-a}) takes the value $r$. Approximate
satisfaction now yields $\rho(a,r)\;<\;K$. A similar
calculation---substituting $b$ for $x$, $d$ for $u$, $f$ for
$v$, and $d(b)$ for $w$---yields $\rho(b,r)\;<\;K$. Finally, by
the triangle inequality, we have
\begin{align*}
          1\EQ \rho(a,b) \;\leq\; \rho(a,r) + \rho(r,b)
               \;<\; K + K \EQ 2K,
\end{align*}
and so $K\geq0.5$, as desired.
 \vspace{0.1in}

\hspace*{-\parindent}%
{\bf Case 2.} $\o g(a,c)=\o g(b,d)$ and $\o g(\o g(\o
d(a),a),e)) = \o g(\o g(\o d(b),b),f)$. In $B$, we may define
\begin{align*}
         s &\EQ \o g(\o g(a,c),\o d(\o g(\o g(\o d(a),a),e))) \\
           &\EQ \o g(\o g(b,d),\o d(\o g(\o g(\o d(b),b),f)))
\end{align*}
(where the second equality clearly comes from the assumptions
of Case 2). Now if we substitute into Equation
(\ref{eq-square-deriv-b})---taking $a$ for $x$, $c$ for $u$,
$e$ for $v$ and $a$ for $w$---the right-hand side of
(\ref{eq-square-deriv-b}) takes the value $a$. Approximate
satisfaction yields $\rho(a,s)<K$. Continuing as in Case 1, we
again obtain $K\geq 0.5$.
\end{Proof}

\hspace*{-\parindent}%
{\bf Remark on the proof of Theorem \ref{th-low-est-squares}.} The
proof never invoked the continuity of $\o d$. Thus, even if we
relax the notion of compatibility to allow non-continuous $\o d$
(along with continuous $\o g$), we still cannot satisfy
(\ref{eq-square-deriv-a}--\ref{eq-square-deriv-b}) any closer than
half the diameter of $[0,1]$. (Similar remarks hold in
\S\ref{subsub-aux-to-log} above.)

\section{Two topological invariants.}   \label{sec-top-inv}
As we have seen, in \S\ref{subsub-depend-metric},
\S\ref{subsub-Y-new-metric}, \S\ref{sub-delta-square-embed-square}
and elsewhere, $\lambda_A$ depends on the choice of the metric $d$
(among those metrics $d$ that yield the given topology on $A$),
and hence is not a topological invariant. For a space $A$ of
finite diameter, one may narrow the choice of $d$ by insisting
that $\text{diam}(A,d)=1$, but $\lambda_A$ is still  not
invariant.

 On the other hand, we can obviously obtain a topological
invariant by considering the extreme values that occur for
$\lambda_{(A,d)}(\Sigma)$ when $d$ is allowed to range over all
appropriate metrics.

\begin{align*}
       \delta^-_A(\Sigma) &\EQ \inf\; \bigl\{\, \lambda_{(A,d)}(\Sigma)\,:\,
                \text{diam}(A,d) \geq 1 \bigr\} \\
       \delta^+_A(\Sigma) &\EQ \sup \bigl\{\, \lambda_{(A,d)}(\Sigma)\,:\,
                \text{diam}(A,d)\leq 1 \bigr\} \\
\end{align*}
where
\begin{align*}
          \text{diam}(A,d) &\EQ \sup \bigl\{d(x,y) \,:\, x,y \in
          A\bigr\}.
\end{align*}
Obviously, the values $ \delta^-_A(\Sigma)$ and $
\delta^+_A(\Sigma)$ are equal, respectively, to the $\inf$ and
$\sup$ of the single set
$$ \bigl\{\, \lambda_{(A,d)}(\Sigma)\,:\,
                \text{diam}(A,d)\EQ 1 \bigr\}.$$

{\bf Short versions:}
\begin{align*}
           \delta^-_A(\Sigma) &\EQ \inf_{ d} \;\inf_{
                            \mathbf A}
                 \; \sup_{\stackrel{(\sigma,\tau)}{\mathbf a}}
                      \;d(\sigma^{\mathbf A}(\mathbf a),
                      \tau^{\mathbf A}(\mathbf a)) \\
\leq\; \delta^+_A(\Sigma) &\EQ \sup_{ d} \;\inf_{ \mathbf A}
                  \;\sup_{\stackrel{(\sigma,\tau)}{\mathbf a}}
                      \;d(\sigma^{\mathbf A}(\mathbf a),
                      \tau^{\mathbf A}(\mathbf a)).
\end{align*}

\subsection{Connection with radius and diameter.}
                 \label{sub-conn-rad-diam}
If $\Sigma$ contains no equation $x_i\wavy x_j$ for $i\neq j$,
then by \S\ref{subsub-lambda}, $\lambda_A(\Sigma)\leq
\text{radius}(A)$. It follows immediately that
\begin{align*}
            \delta^-_A(\Sigma) \;\leq\;
                  \frac{\text{radius}(A)}{\text{diameter}(A)}\;.
\end{align*}
In particular, if $A$ can be metrized with
$\text{diameter}(A)=2\cdot\text{radius}(A)$---as can any simplex
of finite dimension---then the above inequality yields
\begin{align*}
            \delta^-_A(\Sigma) \;\leq\; 0.5.
\end{align*}
The value $0.5$ is realized for $A$ a simplex in
\S\ref{subsub-delt-injbin-in} and \S\ref{subsub--m=nk-delta}
below.

\subsection{Some estimates of $\delta^+_A(\Sigma)$.}
             \label{sub-est-delt-plu}
\begin{theorem}\quad      \label{th-delt-plus}
         $\delta^+_A(\Sigma) \EQ 0\;\mbox{or}\;1$. If\/ $A$ is
         compact, then $\delta^+_A(\Sigma)=0$ iff\/
         $\lambda_{(A,d)}(\Sigma)=0$ for all metrics $d$ defining
         the topology of\/ $A$.  If\/ $A$ is
         compact, then $\delta^+_A(\Sigma)=1$ iff\/
         $\lambda_{(A,d)}(\Sigma)>0$ for all metrics $d$ defining
         the topology of\/ $A$.
\end{theorem}\begin{Proof}
If $\delta$ is a metric for the topology of $A$, then so also is
$f\Compos\delta$, where $f(x) = 1\wedge(2x)$. Thus
$\delta^+_A(\Sigma)$ is the sup of an $f$-closed subset of
$[0,1]$. The last two sentences are immediate from Lemma
\ref{lemma-comp-two-metrics}.
\end{Proof}

\begin{corollary}
If\/ $\delta^-_A(\Sigma)>0$, then $\delta^+_A(\Sigma)=1$.
\end{corollary}

\subsubsection{An inconsistent example where $\delta^+=0$.}
      \label{subsub-delta+0}
Let $A$ be any non-discrete metric space where Tietze's theorem is
applicable. Let $\Gamma$ be the (inconsistent) equations
 \begin{align*}
 F(a,x_0,x_1) &\WAVY x_0\\
  F(b,x_0,x_1) &\WAVY x_1\\
       a &\WAVY b
\end{align*}
(which were first mentioned in \S\ref{sub-sub-no-wavy}). It
follows from \S\ref{sub-sub-no-wavy} that $\delta^+_{A}(\Gamma) =
0$.

\subsubsection{A consistent incompatible example where $\delta^+=0$.}
           \label{subsub-incom-delta+}
In \S\ref{subsub-lim-zero-consis} we considered a consistent
expansion $\Sigma_c$ of lattice theory that is not compatible with
$[0,1]$. We saw there that $\lambda_{([0,1],d)}(\Sigma_c)=0$ for
any metric $d$ that generates the usual topology on $[0,1]$.
Therefore
\begin{align*}
            \delta^+_{[0,1]}(\Sigma_c)\EQ
            \delta^-_{[0,1]}(\Sigma_c)\EQ0.
\end{align*}

Notice that $\delta^-\,=\,0$ does not distinguish  the situation
of this $\Sigma_c$---that all appropriate $\lambda$-values are
zero---from the situation that is seen in
\S\ref{subsub-delta-H-space}, \S\ref{subsub-triode-delta},
\S\ref{subsub-delt-injbin-sq}, \S\ref{subsub--m=nk-delta} and
\S\ref{subsub-nth-delta}, etc., below. In these latter situations,
$\lambda$ is positive but approaches zero (by suitable choice of
metrics); hence $\delta^-$ is again zero. Nevertheless $\delta^+$
does distinguish these two situations.


\subsection{Some estimates of $\delta^-_A(\Sigma)$.}
             \label{sub-estim-delta-minus}
For an initial exploration of the possible values of
$\delta^-_A(\Sigma)$, we revisit the examples of
\S\ref{sub-num-ex}. All these examples will have
$\delta^+_{A}(\Sigma)\EQ1$. In some cases we will be able to state
$\delta^-_A(\Sigma^{\star})=0$; obviously this assertion entails
$\delta^-_A(\Sigma)=0$ (which we will not mention).

\subsubsection{$A=S^1$; $\Sigma=$ semilattice theory.}
              \label{sub-circle-revisit}

From \S\ref{subsub-S1-different}, we see that $
          \delta^-_{S^1}(\Sigma^{\star})\EQ0
$. This method also applies to any theory $\Sigma$ that is
compatible with $[0,1]$.

More generally, let us suppose that $\Sigma$ is not compatible
with $[0,1]$, but that $\delta^-_{[0,1]}(\Sigma)=0$. We shall see
that we still have $\delta^-_{S^1}(\Sigma)\EQ0$. Let us be given
$\varepsilon>0$. By our hypothesis, there is a metric $d$ on
$[0,1]$ such that $\lambda_{([0,1],d)}\leq\varepsilon$. We give
$[0,1]\times\Reals$ the metric
\begin{align*}
           d'((x_0,y_0),(x_1,y_1))\EQ
                     d(x_0,y_0) \,+\,|x_1-y_1|.
\end{align*}

We define $\psi$ as \S\ref{subsub-S1-different}; clearly things go
as before, and one may apply Theorem \ref{thm-inequal-retract} of
\S\ref{sub-inequal-retract}. Much as before, we have
\begin{align*}
           \lambda_{(S^1,d')}(\Sigma)
           \;\leq\;\lambda_{([0,1],d)}(\Sigma)
               \, +\,\varepsilon \;\leq\; \varepsilon
               \,+\,\varepsilon \EQ 2\varepsilon.
\end{align*}

In fact, our proof easily extends to
\begin{align*}
         \delta^-_{S^n}(\Sigma) \;\leq\;
               \delta^-_{[0,1]}(\Sigma).
\end{align*}

\subsubsection{$A=S^2$; $\Sigma=$ H-spaces.}
              \label{subsub-delta-H-space}
From \S\ref{sub-sub-S2-epsilon}, we see that $
          \delta^-_{S^2}(\Sigma^{\star})\EQ0
$. This method also applies to any theory $\Sigma$ that is
compatible with $[0,1]$.

If $\Sigma$ is not compatible with $[0,1]$, but
$\delta^-_{[0,1]}(\Sigma)=0$, then we still have
$\delta^-_{S^2}(\Sigma)\EQ0$. (The proof is like that in
\S\ref{sub-circle-revisit}.)

In the same way, some of the positive $\lambda$-values of
\S\ref{sub-spheres-simple} can be seen to correspond to
zero-values for $\delta^-$.

\subsubsection{Group theory on spaces with the
fixed-point property.} \label{subsub=group-f.p.-delta}

It is immediate from \S\ref{subsub=group-f.p.} that
$\delta^-_A(\Sigma^{\star})\geq0.5$. For a number of spaces, such
as $[0,1]$ and its direct powers (cubes),
\S\ref{sub-conn-rad-diam} yields $\delta^-_A(\Sigma^{\star})=0.5$

\subsubsection{Groups of exponent $2$ on $\Reals$.} It is immediate
from \S\ref{subsub-not-exp2} and \S\ref{subsub-lambda} that (if
$\Sigma$ is as described in \S\ref{subsub-not-exp2}, then)
\begin{align*}
        \delta^-_{\Reals}(\Sigma)\;\geq\;
           \frac{\text{radius}(\Reals)}{2\cdot
           \text{diameter}(\Reals)}\;\geq\;0.25.
\end{align*}

\subsubsection{$A=\Reals$; $\Sigma=$ lattices with zero.}
It is immediate from \S\ref{subsub-depend-metric} that
$\delta^-_{\Reals}(\Sigma^{\star})=0$.

\subsubsection{The triode $Y$ and $\Lambda =$ lattice theory, revisited.}
          \label{subsub-triode-delta}
From \S\ref{subsub-Y-new-metric} we have $\delta^-_Y(\Lambda) =
0$.

\subsubsection{No injective binary on interval, revisited.}
 \label{subsub-delt-injbin-in}
        Let $\Sigma$ be defined as above  in \S\ref{subsub-inf-dim},
and take $A=[0,1]$, in its usual topology. In
\S\ref{subsub-inf-dim} it was proved that
$\lambda_{(A,d)}(\Sigma^{\star})\geq\lambda_{(A,d)}(\Sigma)\geq
0.5$ for any metric $d$ of diameter $1$, and moreover that for one
such metric (namely the usual metric), we have
$\lambda_{(A,d)}(\Sigma^{\star})=\lambda_{(A,d)}(\Sigma)= 0.5$.
Therefore, clearly
\begin{align*}
           \delta^-_{[0,1]}(\Sigma)\EQ \delta^-_{[0,1]}(\Sigma^{\star}) \EQ 0.5.
\end{align*}

\subsubsection{No injective binary on a square, revisited.}
     \label{subsub-delt-injbin-sq}
From \S\ref{sub-delta-square-embed-square},
$\delta^-_B(\Sigma)=0$.

\subsubsection{Generalizations of \S\ref{subsub-delt-injbin-in}
          and  \S\ref{subsub-delt-injbin-sq}.}
               \label{subsub--m=nk-delta}
\S\ref{sub-genzations-delta}   revisited. $\Sigma_{m,n}$ as
before.
\begin{theorem} 
Let $B=[0,1]^k$, with $k \in\{ 1, 2, 3, \ldots\}$. Then
\begin{align*}
     \delta^-_B(\Sigma_{m,n}) \EQ
          \begin{cases}
                0 \quad &\text{if $m\leq nk$} \\
                0.5  \quad &\text{if $m\geq nk + 1$}.
          \end{cases}
\end{align*}
\end{theorem}

\subsubsection{$n$-th power varieties.}
               \label{subsub-nth-delta}
(Refer to \S\ref{sub-nth-power} for definitions.)  By Theorem
\ref{th-nthpow-on-inttok} of \S\ref{subsub-approx-cubes}, if\/
$B=[0,1]^k$ and $\Sigma=\mbox{Sets}^{[n]}$, with\/ $k\geq n$, then
$\delta^-_{B}(\Sigma^{\star})=0$.

\section{Product varieties} \label{sec-prod-vars}
For any sets  $\Gamma$ and $\Delta$ of equations, (finite or
infinite, deductively closed or not), one may construct a new
theory $\Gamma\times\Delta$, which has the following properties
(among others):
\begin{itemize}
 \item[(i)] if $\Gamma$ and $\Delta$ are finite, then
               $\Gamma\times\Delta$ is also finite;
 \item[(ii)]  $\Gamma^{\star}\times\Delta^{\star}\subseteq
                    (\Gamma\times\Delta)^{\star}$;
 \item[(iii)] if $[\Gamma]$ denotes the class of $\Gamma$ in the lattice
       $\mathcal L$ of varieties ordered by interpretability,
       then  $[\Gamma\times\Delta]$ is the meet of $[\Gamma]$ and
       $[\Delta]$ in $\mathcal L$---i.e. $[\Gamma\times\Delta]
        =[\Gamma]\wedge[\Delta]$;
 \item[(iv)] if $A$ is a topological space, then $A\models\Gamma\times
       \Delta$ iff there are spaces $C$ and $D$ such that $A$ is
       homeomorphic to $C\times D$, $C\models\Gamma$ and $D\models
       \Delta$;
 \item[(v)] consequently, if $A$ is a space that cannot be factored
        non-trivially, and if $A\models\Gamma\times\Delta$, then
        either $A\models\Gamma$ or $A\models\Delta$.
\end{itemize}

For an explication of the lattice $\mathcal L$, and for a detailed
explication of these definitions and results, the reader is
referred to \cite{ogwt-mem}. We include a definition of
$\Gamma\times\Delta$ immediately below. In \S\ref{sec-prod-vars}
we shall investigate to what extent there are analogs of the last
two points for approximate satisfaction.

 In defining the product
$\Gamma\times\Delta$ of finite equation-sets $\Gamma$ and
$\Delta$, we will assume that the operations of $\Gamma$ are $G_i$
($i\in I$), and the operations of $\Delta$ are $D_j$ ($j\in J$),
for sets $I$ and $J$ with $I\cap J = \emptyset$. We augment
$\Gamma$ with a new operation $D_j$ and a new equation
\begin{align}
           D_j(x_1,x_2,\cdots) \WAVY x_1,
\end{align}
for each $j\in J$. Similarly, we augment $\Delta$ with a new
operation $G_i$ and a new equation
\begin{align}
           G_i(x_1,x_2,\cdots) \WAVY x_1,
\end{align}
for each $i\in I$. Clearly the augmented $\Gamma$ defines a
variety that is definitionally equivalent to the one defined by
the original $\Gamma$, and similarly for $\Delta$. Now when we
refer to $\Gamma$ and $\Delta$, we mean the augmented
equation-sets.

The operations of $\Gamma\times\Delta$ will be those of $\Gamma$
and $\Delta$, together with a new binary operation denoted $p$. To
define the equation-set $\Gamma\times\Delta$, we first need one
more piece of notation. Let $y$ be a variable not in the set
$\{x_1,x_2,\cdots\}$. For any operation symbol $F\in \Gamma \cup
\Delta\cup\{p\}$, we define $F^R$ and $F^L$ to be the terms
$p(F(x_1,x_2,\cdots),y)$ and $p(y,F(x_1,x_2,\cdots))$,
respectively.  Let $\tau$ be any term whose operations symbols lie
in $\Gamma\cup\Delta\cup\{p\}$, and whose variables are among
$x_1,x_2, \ldots$~. We now recursively define $\tau^R$ as follows.
\begin{itemize}
      \item If $\tau=x_i$, then $\tau^R = p(x_i,y)$ and $\tau^L= p(y,x_i)$.
      \item If $\tau=F(\tau_1,\cdots,\tau_n)$, then
\begin{align}       \label{eq-rec-def-tau-R}
             \tau^R &\EQ p(F(\tau_1^R,\tau_2^R,\cdots),y)\\
             \tau^L &\EQ
             p(y,F(\tau_1^L,\tau_2^L,\cdots)).\nonumber
\end{align}
\end{itemize}

The equations of $\Gamma\times\Delta$ may now be declared as
Equations (\ref{eq-product-a}--\ref{eq-product-g}) that follow.
One easily sees that if $\Gamma$ and $\Delta$ are both finite,
then $\Gamma\times\Delta$ is finite.
\begin{align}            \label{eq-product-a}
               &p(x,x)\WAVY x\; \\
     &p(p(x,y),p(u,v)) \WAVY p(x,v)
     \label{eq-product-b}\\
       &\tau^R \WAVY  p(\tau,y);\quad\quad \tau^L \WAVY\label{eq-product-d}
       p(y,\tau)\\
     &\tau(p(x_1,y_1),p(x_2,y_2),\cdots)\WAVY p(\tau(x_1,x_2,\label{eq-product-e}
     \cdots),\tau(y_1,y_2,      \cdots))\\ &\quad\quad\quad  \quad\quad\quad
                \quad\quad\quad\quad  \quad\quad\quad\quad \nonumber
                 \text{($\tau$ any $\Gamma$-term or any
                 $\Delta$-term)}\\\label{eq-product-f}
                &p(\sigma ,x)\WAVY p(\tau, x)
                \,\quad\quad\quad\quad\quad\quad\quad\quad\quad\quad\quad
                (\sigma\wavy\tau\in \Gamma)\\\label{eq-product-g}
                &p(x,\alpha) \WAVY p(x,\beta)  \quad\quad\quad\quad\quad\quad
                \quad\quad\quad\quad\quad
                (\alpha\wavy\beta\in \Delta).
\end{align}
In equational logic, i.e.\ the logic of exact equality, Equations
(\ref{eq-product-d}) may be deduced from
(\ref{eq-product-a}--\ref{eq-product-b}) and(\ref{eq-product-e}),
and thus are redundant from that point of view. We shall need
them, however, as approximate identities, and in that context they
cannot be deduced.

\subsection{Product identities on product-indecomposable spaces}
               \label{sub-prod-prodindec}
Here we consider a possible analog, for approximate satisfaction,
of the result, mentioned above, that if $A$ is
product-indecomposable and if $A\models\Gamma\times\Delta$, then
$A\models\Gamma$ or $A\models\Delta$. We know of no analog that
holds generally for product-indecomposable spaces, but rather we
must find examples case by case. See Theorems
\ref{th-I-approx-indec} 
and \ref{th-Y-approx-indec}. The first lemma is easy:

\begin{lemma} \label{lem-ifgam-thengamx} If $A$ is a metric space, and if
$\lambda_A(\Gamma)<\varepsilon$, then
$\lambda_A(\Gamma\times\Delta)<\varepsilon$.
\end{lemma}\begin{Proof}
Suppose we have operations that satisfy $\Gamma$ within
$\varepsilon$. We obtain ($\Gamma\times\Delta$)-operations as
follows. The $\Gamma$-operations are retained unchanged. The
$\Delta$-operations, as well as $\o p$, are taken to be first
co-ordinate projection. It is straightforward to verify that
Equations (\ref{eq-product-a}--\ref{eq-product-g}) hold within
$\varepsilon$.
\end{Proof}

It is obvious that the conclusion of Lemma
\ref{lem-ifgam-thengamx} in fact holds if we are given either
$\lambda_A(\Gamma)<\varepsilon$ or $\lambda_A(\Delta)
<\varepsilon$. Therefore an appropriate converse would, when true,
deduce this disjunction (or something slightly weaker) from
$\lambda_A(\Gamma\times\Delta)<\varepsilon$. Such a converse
generally requires a re-adjustment of $\varepsilon.$ Even so, the
converse tends to be false for decomposable spaces $A$---a
specific example will be given as Theorem
\ref{th-prod-smaller-factors} in \S\ref{sub-prod-on-prod} below.
We now turn to two product-indecomposable spaces ($[0,1]$ and
$Y$), where some version of the converse holds---see Theorems
\ref{th-I-approx-indec} and \ref{th-Y-approx-indec} below. First,
a lemma that applies to both of these spaces:

\begin{lemma} \label{lem-Ra-Gamma}
Suppose that $\Gamma\times\Delta$ is satisfiable within
$\varepsilon$ on a metric space $A$, by continuous operations that
include $\o p$ (realizing the $p$ that appears in
{\em(\ref{eq-product-a}--\ref{eq-product-g})}). For $a\in A$, let
$R_a$ and $L_a$ be the images of right and left\/ $\o
p$-multiplication by $a$:
\begin{align*}
        R_a \EQ \{\,\o p(x,a)\,:\,x\in A\};   \quad\quad \quad
        L_a \EQ \{\,\o p(a,x)\,:\,x\in A\}.
\end{align*}
Let us endow $R_a$ and $L_a$ each with the metric obtained by
restricting the given metric on $A$. Then $\Gamma$ is satisfiable
within $3\varepsilon$ by continuous operations on $R_a$, and
$\Delta$ is  is satisfiable within $3\varepsilon$ by continuous
operations on $L_a$.
\end{lemma} \begin{Proof}
It will be enough to prove the part about $\Gamma$ and $R_a$. So
let us suppose that continuous operations $\o G_i$ ($i\in I$) and
$\o D_j$ ($j\in J$) on $A$ form, along with $\o p$, a topological
algebra $\mathcal A$ that satisfies $\Gamma\times\Delta$ to within
$\varepsilon$.

We define a topological algebra $\mathcal R_a$  based on the
subspace $R_a$ as follows. For each operation symbol $G_i$ of
$\Gamma$, we define
\[
           \o G_i^a(x_1,x_2,\cdots) \EQ
                     \o p(\o G_i(x_1,x_2,\cdots),a).
\]
It is clear that each $\o G_i^a$ is continuous and maps into
$R_a$.
We thus have a topological algebra $\mathcal R_a = (R_a,\o
G_i^a)_{i\in I}$. Our task is to prove that $\mathcal R_a$ satisfies
$\Gamma$ within $3\varepsilon$.

So let us first consider an arbitrary term $\tau$ in the language of
$\Gamma$. Let $\o\tau$ denote the realization of $\tau$ in the
original algebra $\mathcal A$ (it is a continuous function
$A^{\omega}\TO A$, which depends on only finitely many variables).
Further, let $\o\tau^a$ denote the realization of $\tau$ in the
algebra $\mathcal R_a$. We propose to prove, by induction, that
\begin{align}               \label{eq-restrict-tau}
      \o\tau^a(a_1,a_2,\cdots)\EQ\o{\tau^R}(a_1, a_2,\cdots,a)
\end{align}
wherever both sides are defined, which is to say for
$a_1,a_2,\cdots\in R_a$, where $\tau^R$ is as defined in
(\ref{eq-rec-def-tau-R}) (and where the $a$ has been substituted for
the variable $y$ of $\tau^R$). Now (\ref{eq-restrict-tau}) obviously
holds for $\tau$ a variable. For $\tau$ a composite term, we may
write $\tau=G_i(\tau_1,\tau_2,\cdots)$. By induction, we have
\begin{align*}
 \o\tau^a(a_1,a_2,\cdots) &\EQ \o G^a_i(\o\tau_1^a(a_1,a_2,\cdots),
     \o\tau_2^a(a_1,a_2,\cdots),\cdots)\\
     &\EQ \o
     G^a_i(\o{\tau_1^R}(a_1,a_2,\cdots,a),\o{\tau_2^R}(a_1,a_2,\cdots,a),
        \cdots)\\
     &\EQ \o p(\o
     G_i(\o{\tau_1^R}(a_1,a_2,\cdots,a),\o{\tau_2^R}(a_1,a_2,\cdots,a),
        \cdots),a)\\
     &\EQ \o{\tau^R}(a_1,a_2,\cdots,a).
\end{align*}

Finally, if $\sigma\wavy\tau$ is an equation of $\Gamma$, then we
may use equations (\ref{eq-product-d}), (\ref{eq-product-f}) and
(\ref{eq-restrict-tau}) to see
\begin{align*}
        \o\sigma^a(a_1,a_2,\cdots) &\EQ \o{\sigma^R}(a_1,
        a_2,\cdots,a)
        \WAVY \o p(\o{\sigma}(a_1,a_2,\cdots),a)\\
        &\WAVY \o p(\o{\tau}(a_1,a_2,\cdots),a) \WAVY \o{\tau^R}(a_1,
        a_2,\cdots,a)\\ &\EQ \o\tau^a(a_1,a_2,\cdots).
\end{align*}
In other words, the realization of $\sigma$ in $\mathcal R_a$ is
equal within $3\varepsilon$ to the realization of $\tau$ in
$\mathcal R_a$.\end{Proof}

The following theorem will be applied in Theorem
\ref{th-01-filter} of \S\ref{sec-filters}, to show that the
condition $\lambda_{[0,1]}(\Sigma)>0$ defines a filter of
theories.
\begin{theorem}  \label{th-I-approx-indec}
$[0,1]$ denotes the unit interval with
its usual metric. Suppose $0<\varepsilon < 1/16$. If
$\lambda_{[0,1]}(\Gamma\times\Delta) <\varepsilon$, then either
$\lambda_{[0,1]}(\Gamma) <4\varepsilon$ or
$\lambda_{[0,1]}(\Delta) <4\varepsilon$.
\end{theorem} \begin{Proof}
Suppose that we have continuous operations $\o p$, $\o G_i$ ($i\in
I$) and $\o D_j$ ($j\in J$) on $[0,1]$, forming a topological
algebra $\mathcal A$, and satisfying equations
$\Gamma\times\Delta$ within $\varepsilon$.

Let us take $R_0$, $R_1$, $L_0$ and $L_1$ as defined in the
statement of Lemma \ref{lem-Ra-Gamma}. By continuity of $\o p$,
each of the four sets is a compact interval $\subseteq[0,1]$. Now
consider $ R_0$ and $L_1$: their intersection contains the point
$\o p(1,0)$; hence their union is again an interval. Now obviously
this interval contains $\o p(0,0)\wavy 0$ and $\o p(1,1)\wavy 1$,
and hence
\[
 R_0 \cup L_1 \; \;\supseteq\;\;
               [\o p(0,0), \o p(1,1)] \;\;\supseteq\;\;
               [\varepsilon,1-\varepsilon],
\]
where the second inclusion comes from (\ref{eq-product-a}).
Therefore
\begin{align}             \label{eq-measure01}
           \mu(R_0) \,+\, \mu(L_1) \;\geq\;
                      1-2\varepsilon,
\end{align}
where $\mu$ stands for Lebesgue measure (in this case, the length
of the interval).
 Similar considerations yield
\begin{align}              \label{eq-measure10}
           \mu(R_1) \,+\, \mu(L_0) \;\geq\;
                      1-2\varepsilon.
\end{align}
It is an easy consequence of (\ref{eq-measure01}) and
(\ref{eq-measure10}) that either
\begin{align}                   \label{eq-measure-0side}
           \mu(R_0) \,+\, \mu(L_0) \;\geq\;
                      1-2\varepsilon
\end{align}
or
\begin{align*}
           \mu(R_1) \,+\, \mu(L_1)\;\geq\;
                      1-2\varepsilon.
\end{align*}
Without loss of generality, we shall assume that
(\ref{eq-measure-0side}) holds.

The next part of the proof (through (\ref{eq-R0-big})) is borrowed
from the example in \S\ref{subsub-inf-dim}. Let us suppose that
the diameter of $L_0$ is attained by the two points $\o p(0,a)$
and $\o p(0,b)$ (in other words, they are its endpoints, in
unspecified order); and further that  $\o p(c,0)$ and $\o p(d,0)$
play the same role for $R_0$. Consider the four real numbers $\o
p(c,a)$, $\o p(d,a)$, $\o p(c,b)$,  $\o p(d,b)$. We may assume,
without loss of generality, that the minimum of these four number
is $\o p(c,a)$, and then that of the two numbers $\o p(d,a)$ and
$\o p(c,b)$, the latter is smaller (or equal). In other words, we
may suppose that
\[
           \o p(c,a) \;\leq\; \o p(c,b) \;\leq\; \o p(d,a).
\]
Applying the Intermediate Value Theorem to the continuous function
$\o p(\cdot,a)$, we infer the existence of $e\in[0,1]$ such that
\[
          \o  p(e,a) \EQ \o p(c,b).
\]

Now, by the approximate satisfaction of (\ref{eq-product-b}), we
have
\[
\o p(0,a) \WAVY \o p(\o p(0,0),\o p(e,a)) \EQ
       \o p(\o p(0,0),\o p(c,b)) \WAVY \o p(0,b).
\]
Thus the distance between the two endpoints of $L_0$ is
$<2\varepsilon$; in other words $\mu(L_0)<2\varepsilon$. It now
follows from (\ref{eq-measure-0side}) that
\begin{align}             \label{eq-R0-big}
               \mu(R_0) \;\geq\; 1-4\varepsilon.
\end{align}

By Lemma \ref{lem-Ra-Gamma}, there are continuous operations on
$R_0$ that satisfy $\Gamma$ within $3\varepsilon$. If we expand
this algebra by a factor of $1/\mu(R_0)$ (the reciprocal of its
length), we now have a topological algebra based on a unit
interval. In this expanded algebra, the laws of $\Gamma$ are
satisfied within
\[
    \frac{3\varepsilon}{\mu(R_0)} \;\leq\;
    \frac{3\varepsilon}{1-4\varepsilon} \;\leq\; 4\varepsilon,
\]
where the final inequality is easily derived from our assumption
that $0<\varepsilon<1/16$.
\end{Proof}

\hspace*{-\parindent}%
{\bf Remark.} It is immediate from (\ref{eq-measure01}) that
either $\mu(R_0)\geq (1-2\varepsilon)/2$ or $\mu(L_1)\geq
(1-2\varepsilon)/2$; without loss of generality, $\mu(R_0)\geq
(1-2\varepsilon)/2$. For small $\varepsilon$, this estimate
differs from (\ref{eq-R0-big}) approximately by a factor of $2$,
yielding approximate satisfaction to within about $8\varepsilon$.
Thus the steps between (\ref{eq-measure01}) and (\ref{eq-R0-big})
are unnecessary, unless one really cares about the $4$ versus the
$8$.

\begin{theorem} \label{th-Y-approx-indec} Let  $Y$ denote the figure-Y
space, as defined and metrized in \S{\em\ref{subsub-Y}}. Suppose
$0<\varepsilon < 1/28$. If $\lambda_{Y}(\Gamma\times\Delta)
<\varepsilon$, then either $\lambda_{Y}(\Gamma) <4\varepsilon$ or
$\lambda_{Y}(\Delta) <4\varepsilon$.
\end{theorem} \begin{Proof}
Suppose that we have continuous operations $\o p$, $\o G_i$ ($i\in
I$) and $\o D_j$ ($j\in J$) on $Y$, forming a topological algebra
$\mathcal Y$, and satisfying equations $\Gamma\times\Delta$ within
$\varepsilon$.

Retaining the notation of \S{\ref{subsub-Y}}, we let $E$, $B$,
$C$, $D$ denote the center point of $Y$ and its three endpoints
(i.e., non-cutpoints). We refer to the segments $BE$, $CE$ and
$DE$ as the {\em legs} of $Y$. Let us consider the sets $R_B$,
$L_B$, $R_C$, $L_C$, $R_D$, $L_D$, as defined in the statement of
Lemma \ref{lem-Ra-Gamma}.

We first prove that either all $R_x$ have measure $<2\varepsilon$
or all $L_x$ have measure $<2\varepsilon$ ($x\in\{A,B,C\}$).
Suppose, to the contrary, that e.g.\ $\mu(L_A)>2\varepsilon$ and
$\mu(R_B)>2\varepsilon$. Thus there exist $P,Q,R,S\in Y$ such that
\begin{align}            \label{largeR-sets}
        d(\o p(A,P),\o p(A,Q)) \,>\, 2\varepsilon;\quad\quad\quad
            d(\o p(R,B),\o p(S,B)) \,>\, 2\varepsilon.
\end{align}
We now arrange the four pairs $(R,P)$, $(R,Q)$, $(S,P)$, $(S,Q)$
in a graph, as follows
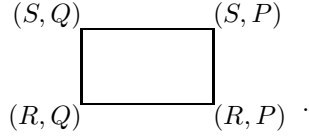
\[\begin{xy}
      0;<5mm,0mm>:
      (3.5,0)="C",  (3.5,2)="D",
      (0,0)="A",  (0,2)="B",
      "C"*!UL(1.0)[o]+{(R,P)}; "D"*!DL(1.0)[o]+{(S,P)};
      "A"*!UR(1.0)[o]+{(R,Q)}; "B"*!DR(1.0)[o]+{(S,Q)};  "A";"C"**@{-};
          "A";"B"**@{-};  "D";"C"**@{-}; "D";"B"**@{-};
\end{xy}.\]
Of these four vertices, two must have $\o p$-values that lie in a
single leg of $Y$. Suppose e.g.\ that $V_1$ and $V_2$ are two of
the four vertices, with $\o p(V_1)$ and $\o p(V_2)$ on a single
leg, say $EX$ where $X\in\{B,C,D\}$, and with
 \[
 d(\o
p(V_1),X)\,\leq\,d(\o p(V_2),X) \,\leq\,\text{the smaller
of}\;\{d(\o p(V_3),X),
             d(\o p(V_4),X)\}
\]
for $V_3$ and $V_4 $ the other two vertices (labeled in any
manner). Now either $V_3$ or $V_4$ is related to $V_1$ in the
above graph; by choice of notation, we will assume that $V_1$ is
related to $V_3$. From our choice of the relative positions of $\o
p(V_i)$ ($i=1,2,3$), one easily sees that every connected subset
of $Y$ that contains $\o p(V_1)$ and $\o p(V_3)$ also contains $\o
p(V_2)$.

From here the proof depends slightly on whether the edge $V_1V_3$
is vertical or horizontal. We present only the vertical case; the
horizontal case is similar. Without loss of generality, we may in
fact assume these values:
\[
       V_1\EQ (S,Q);\quad\quad  V_2\EQ (S,P);\quad\quad
        V_3\EQ (R,Q).
\]
By continuity, the image of $\o p(\cdot,Q)$ is a connected subset
of $Y$ that contains  $\o p(V_1)$ and $\o p(V_3)$; by our remarks
above, it also contains $\o p(V_2)$. Thus there exists $T\in Y$
such that
\[
            \o p(T,Q)\EQ \o p(S,P).
\]
Finally
\[
           \o p(A,P) \WAVY \o p(\o p(A,A),\o p(S,P)) \EQ
                  \o p(\o p(A,A),\o p(T,Q)) \WAVY \o p(A,Q),
\]
and so $d(\o p(A,P),\o p(A,Q))<2\varepsilon$, in contradiction to
(\ref{largeR-sets}). This contradiction completes the proof of our
claim that  either all $R_x$ have measure $<2\varepsilon$ or all
$L_x$ have measure $<2\varepsilon$ ($x\in\{A,B,C\}$). Without loss
of generality, we may from now on suppose that all $L_x$ have
measure $<2\varepsilon$.

Now let us consider $R_B$, $L_C$ and $L_D$. Each is a connected
subset of $Y$, the point $\o p(C,B)$ lies in $R_B\cap L_C$, and
the point $\o p(D,B)$ lies in $R_B\cap L_D$. Therefore $R_B\cup
L_C\cup L_D$ is connected. Moreover it contains the three points
\[
            \o p(B,B),\quad \o p(C,C), \quad  \o p(D,D),
\]
which are within $\varepsilon$ of $B$, $C$ and $D$, respectively,
by (\ref{eq-product-a}) The only such set is all of $Y$, minus
three intervals, each of size $<\,\varepsilon$. Thus
\[
              \mu(R_B\cup L_C\cup L_D) \;\geq\; 3-3\varepsilon.
\]
So finally
\begin{align*}
         \mu(R_B) \;&\geq\;  \mu(R_B\cup L_C\cup L_D) - \mu(L_C) -
         \mu(L_D)\\
         &\geq (3-3\varepsilon) - 2\varepsilon - 2\varepsilon \EQ
              3 - 7\varepsilon.
\end{align*}

By Lemma \ref{lem-Ra-Gamma}, there are continuous operations
making the space $R_B$ into a topological algebra modeling
$\Gamma$ to within $3\varepsilon$ . Clearly $R_B$ is homeomorphic
to $Y$, and will become isometric to $Y$ upon rescaling each leg.
The maximum rescaling factor occurs if all the error is in a
single leg, namely $1/1-7\varepsilon$. Therefore, after rescaling,
the maximum error in the equations of $\Gamma$ will be
\[
          \frac{3\varepsilon}{1-7\varepsilon} \;\leq\;
          4\varepsilon,
\]
where the final inequality is easily derived from our assumption
that $0<\varepsilon<1/28$.
\end{Proof}

There are results like Theorems \ref{th-I-approx-indec} and
\ref{th-Y-approx-indec}, for example for $\lambda_{S^1}$ and
$\lambda_{S^n}$ ($n\neq 1,3,7$), that hold simply because in these
special cases, approximate satisfaction of $\Sigma$ implies a
strict structural property of $\Sigma$. Consider the case of
$\lambda_{S^n}$ for the indicated values of $n$. If
$\lambda_{S^n}(\Gamma\times\Delta)<\varepsilon$ with
$\varepsilon\leq 1$, then by \S\ref{sub-sub-Sn},
$\Gamma\times\Delta$ is undemanding. From there, it is not hard to
see that either $\Gamma$ or $\Delta$ must be undemanding. Hence
$\lambda_{S^n}(\Gamma)=0$ or $\lambda_{S^n}(\Delta)=0$.

The case of $\lambda_{S^1}$ is similar. If If
$\lambda_{S^1}(\Gamma\times\Delta)<\varepsilon$ with
$\varepsilon\leq 1$, then by \S\ref{sub-sub-S1}
$\Gamma\times\Delta$ is Abelian. From there, it is not hard to see
that either $\Gamma$ or $\Delta$ must be Abelian. Hence
$\lambda_{S^n}(\Gamma)=0$ or $\lambda_{S^n}(\Delta)=0$.

\subsection{Product identities on a product of two spaces}
\label{sub-prod-on-prod} In comparing $\lambda_A$ and $\lambda_B$
with $\lambda_{A\times B}$, one needs to have a metric on $A\times
B$ that is related in some way to the metrics on $A$ and $B$. For
example, if $A=(A,d)$ and $B=(B,e)$ are metric spaces, one can use
the Pythagorean metric on $A\times B$, namely $\rho$, where
\[
             \rho((a,b),(a',b'))\EQ \sqrt{d(a,a')^2+e(b,b')^2}.
\]
Two further possible metrics for $A\times B$ are these:
\begin{align*}
          \sigma((a,b),(a',b')) &\EQ \max(d(a,a'),e(b,b'))\\
          \tau((a,b),(a',b')) &\EQ d(a,a')+e(b,b').
\end{align*}
For  \S\ref{sub-prod-on-prod} only, we will retain the notations
$\rho$, $\sigma$ and $\tau$ for these three metrics. In any case,
the three are related by the obvious inequalities
\[
          \sigma(\alpha,\alpha') \;\leq\;  \rho(\alpha,\alpha') \;\leq\;
           \tau(\alpha,\alpha') \;\leq\;  2\sigma(\alpha,\alpha').
\]
Thus in any case, switching between $\rho$, $\sigma$ and $\tau$
will change $\lambda_{A\times B}$ at most by a factor of $2$. In
Theorem \ref{th-prod-smaller-factors} below, we will use XXX as
our metric on $A\times B$, simply because the calculations are
easier with this metric.

The following theorem is easy; its proof will be omitted. Versions
obviously exist for other metrics on $A\times B$, as well.
 \begin{theorem} \label{th-prod-on-prod}
If $\lambda_{(A,d)(\Gamma)}(\Gamma)<\delta$ and
$\lambda_{(B,e)}(\Delta)<\varepsilon$, then
\begin{align*}
    \lambda_{(A\times B,\rho)}(\Gamma\times\Delta) &\;<\;
                                    \sqrt{\delta^2+\varepsilon^2}\\
    \lambda_{(A\times B,\sigma)}(\Gamma\times\Delta)
                             &\;<\; \max(\delta,\varepsilon)\\
    \lambda_{(A\times B,\tau)}(\Gamma\times\Delta)
                               &\;<\;      \delta+\varepsilon.\ENDPROOF
\end{align*}
\end{theorem}

Theorem \ref{th-prod-on-prod} makes it rather easy for a product
of metric spaces to satisfy the product equations
$\Gamma\times\Delta$ approximately. This realization will
facilitate our seeing that the approximate-primality results of
\S\ref{sub-prod-prodindec} (Theorems \ref{th-I-approx-indec} and
\ref{th-Y-approx-indec}) do not extend to product-decomposable
spaces. A specific counter-example is provided by Theorem
\ref{th-prod-smaller-factors} immediately below.

 In Theorem \ref{th-prod-smaller-factors},
$(S^1,d)$ is the one-dimensional sphere with the diameter-1 metric
$d$ that is proportional to arc length along shortest paths, and
$(Y,e)$  is the triode described in \S\ref{subsub-Y}. $S^1\times
Y$ will be given the metric $\tau$ described above (addition of
metrics in the two components). $\Gamma$ is the theory of groups
(a weaker version could be used instead---see
\S\ref{subsub-group-cyl}), and $\Sigma$ is semilattice theory. (Or
one may use $\Sigma=$ commutative idempotent algebras of type
$\langle2\rangle$.) Topologically, it has long been known that
$S^1\times Y\models\Gamma\times\Sigma$, that $S^1\times
Y\not\models\Gamma$, and that $S^1\times Y\not\models\Sigma$. The
following metric theorem is sharper.
\begin{theorem} \label{th-prod-smaller-factors}
$S^1\times Y\models\Gamma\times\Sigma$---and thus
$\lambda_{S^1\times Y}(\Gamma\times\Sigma)=0$---while
$\lambda_{S^1\times Y}(\Gamma)\geq 0.1$ and $\lambda_{S^1\times
Y}(\Sigma)\geq 1.0$.
\end{theorem} \begin{Proof}
Since $S^1\models\Gamma$ and $Y\models\Sigma$, the first assertion
is immediate from Theorem \ref{th-prod-on-prod} (or from Point
(iv) at the start of \S\ref{sec-prod-vars}). The estimate for
$\lambda_{S^1\times Y}(\Gamma)$ comes from
\S\ref{subsub-group-cyl}, and the estimate for $\lambda_{S^1\times
Y}(\Sigma)\geq 0.4$ comes from \S\ref{subsub-semil-homotop}.
\end{Proof}

\section{Approximate satisfaction by piecewise linear (simplicial) operations.}
\label{sec-simplicial}

In \S\ref{sec-simplicial} we will deal with finite simplicial
complexes, corresponding to compact spaces that can be
triangulated. (Most of the compact spaces in this paper fall into
this category.) Our objective will be to prove that continuous
satisfaction within $\varepsilon$ is equivalent to piecewise
linear (i.e. simplicial) satisfaction within $\varepsilon$, in a
manner that can be recursively enumerated (see results of
\S\ref{sub-main-simp-alg}).

\subsection{Simplicial complexes and maps} \label{sub-simpmaps}

For definiteness, we paraphrase the definitions and notation of
Spanier \cite[pages 108--128]{spanier}. A {\em (simplicial) complex}
is a family $K$ of nonempty finite sets such that
\begin{itemize}
   \item[(a)] if $v\in\bigcup K$, then $\{v\}\in K$;
   \item[(b)] if $\emptyset\neq s\subseteq t\in K$, then $s\in K$.
\end{itemize}
The elements of $K$ are sometimes called the (abstract) {\em
simplices} of $K$. The elements of $V=\bigcup K$ are called {\em
vertices} of $K$. (Condition (a) provides a one-one correspondence
between vertices of $K$ and  one-element simplices of $K$.) If $K$
and $L$ are complexes and $L\subseteq K$, then $L$ is called a {\em
subcomplex} of $K$. For each simplex $s$ of $K$, we let $\o s$
denote the set of all non-empty subsets of $s$, and $\dot s$ denote
$\o s\smallsetminus\{s\}$. One easily checks that $\o s$ and $\dot s
$ are both subcomplexes of $K$.

The connection with metric topology is this. The {\em standard
geometric realization} of $K$ is a metric space $|K|$ whose
underlying set  consists of all vectors $\alpha\in[0,1]^V$ such that
\begin{itemize}
\item[(c)] The {\em carrier} of $\alpha$,
       namely the set $\{v\in V:\alpha_v >0\}$, is a simplex of
       $K$. (In particular, it is a finite subset of $K$.)
\item[(d)] $\sum_{v\in V}\alpha_v = 1$.
\end{itemize}

For finite $K$---for which $V$ is also finite---the metric on $|K|$
is the restriction of the usual Euclidean metric on $[0,1]^V$. In
other words,
\begin{align}   \label{eq-geom-real-dist}
    d(\alpha,\beta) \EQ \sqrt{\smash[b]
                      {\sum_{v\in V} \,(\alpha_v-\beta_v)^2}}.
\end{align}
In fact, even for infinite $K$, the sum appearing in Equation
(\ref{eq-geom-real-dist}) is finite, and hence
(\ref{eq-geom-real-dist}) may be used as a distance formula for
all $K$, finite or infinite. See Spanier \cite[{\em loc.\
cit.}]{spanier} for the fact that $|K|$ is compact if and only if
$K$ is finite. (In fact [{\em loc.\ cit.}] another topology is
often used in the case of infinite $K$.) In \S\ref{sec-simplicial}
our interest is in compact spaces; hence we will mainly work with
finite $K$.

For $w\in V$, define $w^{\star}\in|K|$ {\em via}
\[
            w^{\star}(v) \EQ \begin{cases}
                             \;1 &\text{if $v=w$};\\
                             \;0 &\text{otherwise}.
                             \end{cases}
\]
In this notation, the $\alpha$ introduced in clause (c) above may
also be denoted $\sum_{v\in s}\alpha_v v^{\star}$. It is then not
hard to see that $w\GOESTO w^{\star}$ maps $V$ injectively to $|K|$,
and that for each simplex $s$, the convex hull of
$\{w^{\star}\,:\,w\in s\}$ in $[0,1]^V$ is the set of points in
$|K|$ with carrier $\subseteq s$. This convex hull is denoted $|s|$.
Since $s\subseteq V$, there is a natural projection $[0,1]^V \TO
[0,1]^s$; this projection restricts to a homeomorphism $|s|\TO |\o
s|$. The space $|\o s|$---the geometric realization of an abstract
simplex---is a {\em geometric $n$-simplex} for some $n$, i.e., the
set of points $(\alpha_0,\cdots,\alpha_n)\in[0,1]^{n+1}$ satisfying
$\sum \alpha_i = 1$, with the usual Euclidean metric.

For $s\in V$ we call $|s|$ the {\em closed geometric simplex of\/}
$s$. The {\em open geometric simplex of\/} $s$ is
\[
\la s\ra \EQ |s|\,\smallsetminus\,
                \bigcup\,\{\,|t|\,:\,t\subseteq s, \;t\neq s\}.
\]
$\la s\ra$ is an open subset of $|s|$, homeomorphic to an open
$n$-disk (for appropriate $n$). Each vector $\alpha\in|K|$ lies in
$\la s\ra$, where $s$ is the carrier of $\alpha$, and in no other
set of the form $\la t \ra$. Thus the sets $\la s\ra$ partition the
space $|K|$.

A {\em simplicial map} $\phi\FROM K_1\TO K_2$ is  a map between the
corresponding vertex-sets, $\phi\FROM V_1\TO V_2$, such that if
$s\in K_1$, then $\phi[s]\in K_2$. The {\em geometric realization
of\/} $\phi$ is the map $|\phi|\FROM|K_1|\TO|K_2|$ that is defined
as follows: if
\[
         \alpha\EQ \sum_{v\in s} \,\alpha_v\, v^{\star},
\]
where each $\alpha_v>0$ and $\sum\alpha_v = 1$, then
\[
         |\phi|\,(\alpha) \EQ \sum_{v\in s} \,\alpha_v \,
           (\phi(v))^{\star}.
\]
It is not hard to check that the right-hand-side here lies in $|K_2$
and is well-defined, and that the resulting map between metric
spaces is continuous.

\subsection{The product complex.} \label{sub-prod}
Suppose that $K_1$ and $K_2$ are complexes, and that $<$ is a total
order on $V_1\cup V_2$. We will define a complex $K_1\times_{<}K_2$
with vertex-set $V_1\times V_2$.

Suppose that $s^1_m=\{v_0^1,\cdots,v_m^1\}$ is an $m$-simplex of
$K_1$, with $v_0^1<\cdots<v_m^1$, and that
$s^2_n=\{v_0^2,\cdots,v_n^2\}$ is an $n$-simplex of $K_2$, with
$v_0^2<\cdots<v_n^2$. Now let
\begin{align}         \label{eq-prod-simplex}
   \lambda \EQ \bigl(
  (r_0,s_0),\;  (r_1,s_1),\;  \cdots  \;(r_{m+n},s_{m+n})
         \bigr)
\end{align}
be a finite sequence of pairs with $r_0=s_0=0$, and such that, for
each $i$, either $r_{i+1}=r_i +1$ and $s_{i+1} = s_i$, or
$r_{i+1}=r_i$ and $s_{i+1} = s_i +1$ and such that the first
alternative happens $m$ times and the second alternative happens
$n$ times. (It of course follows that $r_{m+n}=m$ and $s_{m+n}=n$.
There are clearly ${m+n}\choose{n}$ sequences
(\ref{eq-prod-simplex}).) We then define
\[
       s^{\lambda}_{m+n} \EQ
       \bigl\{\,(v_{r_0}^0,v_{s_0}^1), \cdots\;, (v_{r_m}^0,v_{s_n}^1)
                    \,\bigr\}
\]
to be an $(m+n)$-simplex of $K_1\times_{<}K_2$. We omit the
(standard) proof that the set of all such $s_{m+n}^{\lambda}$ (for
all appropriate $\lambda$, for all $m$ and $n$, and for $s_m^1$ and
$s_n^2$ ranging over all simplices in $K_1$ and $K_2$,
respectively), form a complex, which we will denote
$K_1\times_<K_2$.

It is not hard to check (using formula (\ref{eq-prod-simplex})) that
each of the two coordinate projections $\pi_i\FROM V_1\times V_2\TO
V_i$ is a simplicial map. Therefore we have continuous maps
\[
           |\pi_i|\FROM |K_1\times_<K_2| \TO |K_i|
\]
for $i=1,2$. We therefore have the continuous map
\begin{align}            \label{eq-triangulate-prod}
           |\pi_1|\!\times\!|\pi_2|\,\FROM \;|K_1\times_<K_2|
           \;\TO \;|K_1|
                    \times|K_2|.
\end{align}
We omit the (standard) proof that the map appearing in
(\ref{eq-triangulate-prod}) is a homeomorphism; in other words,
$K_1\times_< K_2$ triangulates the product space $|K_1|\times|K_2|$.

In fact, there is a standard way to triangulate any finite product
of geometric realizations, $|K_1|\times\cdots\times|K_n|$. One may
give an analogous construction, or simply iterate the binary product
that we already have. We omit the details.

\subsection{Barycentric subdivision } \label{sub-bary} Let $K$ be a complex, with
vertex-set $V$, and let $|K|$ be the standard geometric realization
of $K$ (\S\ref{sub-simpmaps}). For each $n$, and each $n$-simplex
$s\in K$, we define the {\em barycenter} of $s=\{w_0,\cdots,w_n\}$
to be the point in $|K|$ defined by
\begin{align}
              b(s)\EQ \frac{1}{n+1}\;\sum_{w\in s}\,\, w^{\star}\EQ
                 \frac{1}{n+1}\;\sum_{j=0}^n \, w_j^{\star}\,.
\end{align}
In other words,
\begin{align*}
           [b(s)]\,(v) \EQ \begin{cases}
                  \frac{1}{n+1} &\text{if  } v\in s;\\
                    0           &\text{otherwise.}
                    \end{cases}
\end{align*}
It is not hard to see that $b(s)$ is in the open simplex $\la
s\ra$ of $|K|$; therefore there is a one-one correspondence
between $K$ and the set of all barycenters.

We define a complex $K'$, called the {\em barycentric subdivision}
of $K$, as follows. The vertex set $V'$ for $K'$ is the set of
barycenters $b(s)$ for $s\in K$, and an $n$-simplex of $K'$ is a
finite set
\begin{align}    \label{eq-K-prime-simplex}
    \Gamma \EQ   \bigl\{ \,b(s_0),\,\cdots\,,\,b(s_n)\, \bigr\},
\end{align}
where
\begin{align}         \label{eq-prime-simp-condition}
         s_0\;\subset\; s_1 \;\subset\;\cdots
                              \;\subset \;s_n \;\in \;K.
\end{align}
Clearly $K$ and $K'$ have the same dimension. Conditions (a) and
(b) of \S\ref{sub-simpmaps} are easily seen to hold for $K'$, and
hence $K'$ is a complex. We will see that in fact $|K'|$ is
homeomorphic to $|K|$.

Since $V'$ is a set of elements of $|K|$, we have the identity
inclusion $\iota\FROM V'\TO|K|$. Let us suppose we have a
$K'$-simplex $\Gamma$ as in (\ref{eq-K-prime-simplex}).  Then
Condition (\ref{eq-prime-simp-condition}) assures us that each
vertex $b(s_j)$ of $\Gamma$ in fact lies in the realization $|s_n|$
of $s_n$ (a geometric simplex $\subseteq |K_n|$). Thus $\iota$ has
an affine-linear extension $|\iota|$ to the geometric realization
$|K'|$ of $K'$. In particular, on the simplex $\Gamma$ defined by
(\ref{eq-K-prime-simplex}) and (\ref{eq-prime-simp-condition}), we
have
\[
          |\iota|\;\bigl(\,\sum_{j=0}^n \alpha_j\,
          [b(s_j)]^{\star}\bigr)
          \EQ
          \sum_{j=0}^n\alpha_j\,b(s_j) \,\in\,|K|
\]
for all $\alpha_0,\cdots,\alpha_n$ with $\sum\,\alpha_j=0$. For
the proof of the following well-known result see \cite[{\em loc.\
cit.}]{spanier}.

\begin{lemma}          \label{lem-sub-homeo}
$|\iota|$ maps $|K'|$ homeomorphically to $|K|$.\ENDPROOF
\end{lemma}

Recall that we have metrized $|K|$ according to formula
(\ref{eq-geom-real-dist}). We could of course metrize the geometric
realization $|K'|$ according to (\ref{eq-geom-real-dist}), but in
that case the homeomorphism $|\iota|$ of Lemma \ref{lem-sub-homeo}
would not be an isometry. We prefer to keep a single metric; in
other words, we shall metrize  $|K'|$  in the unique way that makes
$|\iota|$ an isometry.

We can now iterate the process of subdivision, yielding a tower of
isometries\footnote{We call them all $|\iota|$, even though,
strictly speaking, they are distinct maps.}
\begin{align*}
       \;\cdots\;\Mapnamed{|\iota|}\;|K^{(m)}|\;
            \Mapnamed{|\iota|}\;\cdots\;\Mapnamed{|\iota|}\;
       |K''|\;\Mapnamed{|\iota|}\;|K'|\;\Mapnamed{|\iota|}\;|K|.
\end{align*}
We saw in \S\ref{sub-simpmaps} that the open simplices $\la s\ra$ of
$|K|$ partition $|K|$. This result applies, of course, to each
iteration of the subdivision. Thus, for example, the open simplices
of $|K'|$ partition $|K'|$. These can be carried by $|\iota|$ to
sets in $|K|$, which of course partition the original $|K|$. It is
not hard to check that, if $\Gamma$ is a simplex of $K'$ (see
(\ref{eq-K-prime-simplex})) with largest vertex $s_n\in K$, then
$|\iota|\,[\la \Gamma\ra]\subseteq \la s_n\ra$. In other words,

\begin{lemma} \label{lem-partition}  The partition of $|K|$ by
sets of the form $|\iota|\,[\la \Gamma\ra]$ refines the partition by
sets of the form $\la s\ra$.\ENDPROOF
\end{lemma}

The next lemma is almost obvious from what has come before. The
notation continues from above.
\begin{lemma}
$|\iota|\,[\la \Gamma\ra]$ is a subset of $\la s\ra$ that is defined
by linear inequalities.\ENDPROOF
\end{lemma}

The partition of $|K|$ by sets of the form $\la s\ra$ will be called
the {\em natural triangulation} of $|K|$. The partition of $|K|$
mentioned at the start of Lemma \ref{lem-partition} will be called
the {\em first subdivision} of the natural triangulation. Obviously
one also has the second, third, \ldots, $m^{\text{th}}$, \ldots
subdivisions of the natural triangulation. (The elements of the
$m^{\text{th}}$ subdivision are the $|\iota|^m$-images of sets in
the natural partition of $|K^{(m)}|$.)

If $A$ is a metric space and if $\mathcal Z$  is a family of subsets
of $A$, the {\em mesh} of $\mathcal Z$ is
\[
            \sup_{B\in \mathcal Z}\, \sup_{b,c\in B}\, d(b,c).
\]
Suppose that $K$ is a complex of dimension $n$.  It is clear from
the definitions in \S\ref{sub-simpmaps} that the natural
triangulation of $|K|$ has mesh $\leq\sqrt{2}$. The next lemma
follows from some elementary considerations of metric geometry;
see e.g.\ \cite[{\em loc.\ cit.}]{spanier}.

\begin{lemma}   \label{lem-subd-mesh}
Let $K$ be a complex of dimension $n$. The $m^{\text{th}}$
subdivision of the natural triangulation has mesh
\[
                 \sqrt{2}\;\left(\frac{n}{n+1}\right)^m.\;\ENDPROOF
\]
\end{lemma}
In particular, the iterated subdivisions have mesh approaching zero.

\subsection{Simplicial approximation of continuous operations}
\label{sub-simp-approx} At the end of \S\ref{sub-simpmaps}, we
defined the notion of simplicial maps and their associated geometric
realizations. Here we extend this idea as follows. Let  $K$ and $L$
be complexes. An $(M,N)$-{\em simplicial map} from $K$ to $L$ is a
simplicial map $\phi\FROM K^{(M)}\TO L^{(N)}$, where $K^{(M)}$ and
$L^{(N)}$ are as defined in \S\ref{sub-bary}. Its {\em ground-level
geometric realization} is the (piecewise linear, continuous)
composite map
\begin{align}                     \label{eq-phi-0-def}
  |\phi|_0 \;\EQ \;|K| \;\xrightarrow{\;\;\INV{{(|\iota|^M)}}\;\;} \;|K^{(M)}|
              \; \xrightarrow{\;\;|\phi|\;\;}\;|L^{(N)}|\;
              \xrightarrow{\;\;|\iota|^N}\;\;|L|\,.
\end{align}

The virtues of such maps are two: (1) they are amenable to an
algorithmic approach; (2) (allowing free choice of $M$ and $N$)
they can approximate any continuous function between spaces $|K|$
and $|L|$ (for $K$ and $L$ finite). As for (1), that will be the
topic of \S\ref{sub-main-simp-alg} below. As for (2), we have the
following version of the classical {\em Simplicial Approximation
Theorem}.

\begin{theorem} \label{th-SAT} For $K$ finite and $L$ an arbitrary
complex, given real $\varepsilon>0$ and continuous $f\FROM |K|\TO
|L|$, there exist positive integers $M$ and $N$ and a simplicial map
$\phi\FROM K^{(M)}\TO L^{(N)}$ such that $|\phi|_0$ approximates $f$
within $\varepsilon$ on $|K|$.
\end{theorem}\begin{Proof}
The traditional version of this theorem (see e.g.\ Spanier
\cite[Theorem 8, page 128]{spanier}) says this. Given $f$ and $N$,
there exist $M$ and a simplicial map $\phi\FROM K^{(M)}\TO L^{(N)}$
such that $|\phi|_0$ approximates $f$ in the following sense: if
$x\in |K|$ and $f(x)\in |\iota|^N[\la s\ra]$ for a simplex $s$ of
$L^{(N)}$, then $|\phi|_0(x)\in |\iota|^N[|s|]$. In order to obtain
the conclusion of Theorem \ref{th-SAT} as stated, we merely need to
choose $N$ here large enough that each geometric simplex $|s|$ of
$L^{(N)}$ has diameter $<\varepsilon$. That this is possible is
immediate from Lemma \ref{lem-subd-mesh}.
\end{Proof}

Now, for definiteness, let us suppose that the ground-level
complexes under discussion (i.e.\ $K$, $L$, etc., but not $K^{(M)}$,
etc.) all have vertices in a fixed set $V_1$, which has a fixed
strict total order $<$. If $K$ is a simplicial complex and
$n\in\Integers^+$, then by $K^n$ we mean the complex
\[
              (\cdots((K\times_<K)\times_<K)\times_<\dots\times_<K),
\]
having $n$ factors $K$, that is formed by an iteration of the
construction in \S\ref{sub-prod}. An iterated version of Equation
(\ref{eq-triangulate-prod}) provides a natural homeomorphism
\begin{align}                  \label{eq-pi-n}
       |\pi|^n\EQ(\cdots(|\pi_1|\times |\pi_2|)\times
                   \dots\times |\pi_n|)\;
         \FROM\;|K^n|\,\TO \,|K|^n.
\end{align}
Again, this map is a homeomorphism, but not an isometry according to
the way the metrics have been defined. We get around this by
redefining the metric on $|K^n|$ so that $|\pi|^n$ becomes an
isometry.

\begin{theorem} \label{th-SAT-op} For $K$ finite, $n\in\Integers^+$,
and $L$ an arbitrary complex, given real $\varepsilon>0$ and
continuous $F\FROM |K|^n\TO |L|$, there exist positive integers $M$
and $N$ and a simplicial map $\phi\FROM (K^n)^{(M)}\TO L^{(N)}$ such
that $|\phi|_0\Compos\INV{(|\pi|^n{})}$ approximates $F$ within
$\varepsilon$ on $|K|^n$.
\end{theorem}\begin{Proof}
$F\Compos|\pi|^n$ maps $|K^n|$ to $|L|$; hence by the Simplicial
Approximation Theorem (\ref{th-SAT}), $F\Compos|\pi|^n$ is
$\varepsilon$-approximated by $|\phi_0|$ for some simplicial map
$\phi\FROM (K^n)^{(M)}\TO L^{(N)}$. Therefore $F$ is
$\varepsilon$-approximated by $|\phi|_0\Compos\INV{(|\pi|^n{})}$.
\end{Proof}

\subsection{Approximation of term operations}
\label{sub-approx-term} The considerations of
\S\ref{sub-approx-term} apply to operations on any compact metric
space; triangulation is not required here.

We consider a finite similarity type $\bigl\la F_i\bigr\ra_{i<N}$,
with each $F_i$ of arity $n_i$. The {\em depth} of a term $\tau$ in
this language is defined recursively as follows:
\begin{itemize}
\item[(i)] If $\tau$ is a variable, then $\tau$ has depth $0$;
\item[(ii)] if $\tau=F_i(\tau_1,\cdots,\tau_{n_i})$, then
\[
      \text{depth}(\tau)\EQ 1\,+\,
                      \max\,\{\,\text{depth}(\tau_j)\,:\,
                                  1\leq j\leq n_i\}.
\]
\end{itemize}
Suppose now that $\mathbf A=(A,F_i^{\mathbf A})_{i<N}$ is a
topological algebra of this similarity type, with $A$ a compact
metric space. As is well known, the operations $F_i^{\mathbf A}$ are
in fact {\em uniformly continuous}.

The aim in \S\ref{sub-approx-term} is to see how appropriate
approximations to the operations $F_i^{\mathbf A}$ (say by piecewise
linear operations, or by differentiable operations) can lead to
approximations of term operations $\tau^{\mathbf A}$. Thus suppose
we are given a real $\varepsilon>0$, and we would like to be able to
approximate any
depth-$n$ term operation\footnote{%
The precise meaning of ``approximating a term-operation'' is
deferred until the statement of Lemma \ref{lem-approx-term}.}
$\tau^{\mathbf A}$ within $\varepsilon$. We recursively
define\footnote{%
In the forthcoming work \cite{wtaylor-dise} we shall say that the
operation $\o F$ is {\em constrained} by
$(\varepsilon_1,\varepsilon_0/2)$, by
$(\varepsilon_2,\varepsilon_1/2)$, and so on.} real numbers
$\varepsilon_k>0$ as follows:
\begin{itemize}
\item[(i)] $\varepsilon_0\,=\,\varepsilon$;
\item[(ii)] suppose that $\varepsilon_k$ has already been defined.
   By uniform continuity, there exists\footnote{%
In the unobstructed realm of all continuous operations, there is
no effective value that can be assigned to $\varepsilon^{\star}$.
This is why we need the enumerative approach in
\S\ref{sub-main-simp-alg} below.} real
   $\varepsilon^{\star}>0$
   such that, for all $i<N$, and for all $x_j,y_j\in A$ ($1\leq
   j\leq n_i$), we have:
   \begin{align} \nonumber
     \text{if}\;\; (d(x_j,y_j) < &\varepsilon^{\star} \;\;\text{for}
            \;\; 1\leq j\leq n_i,\;\;\;\text{then}\\
                 &d(F_i^{\mathbf A}(x_1,\ldots,x_{n_i}),
                 F_i^{\mathbf A}(y_1,\ldots,y_{n_i}))\;<\;\frac
                 {\varepsilon_k}{2}.   \label{eq-rec-eps-k}
   \end{align}
   We choose such an $\varepsilon^{\star}$ and define
   \begin{align}   \label{eq-rec-def-epsk}
            \varepsilon_{k+1}\EQ \min\,\{\,\varepsilon^{\star},\,
                     \frac{\varepsilon_k}{2}\,\}.
   \end{align}
\end{itemize}

The next lemma says, in paraphrase, that if we approximate the
operations $F_i^{\mathbf A}$ within $\varepsilon_M$, then for all
terms $\tau$ of depth $\leq M$, we will also approximate
$\tau^{\mathbf A}$ within $\varepsilon$. The proof is a simple
recursion based on the triangle inequality.

\begin{lemma} \label{lem-approx-term}
Let us be given a topological algebra $\mathbf A=(A,F_i^{\mathbf
A})_{i<N}$ (as above), an integer $M\geq1$, and a real
$\varepsilon>0$. Let $\varepsilon_M$ be as defined above in
{\em(\ref{eq-rec-eps-k}--\ref{eq-rec-def-epsk})}. If\/ $\mathbf
B=(A,F_i^{\mathbf B})_{i<N}$ is a similar algebra defined on $A$
(topological or not), and if the operations $F_i^{\mathbf B}$
satisfy
\begin{align}                \label{eq-ops-M-1-close}
              d( F_i^{\mathbf B}(\mathbf x), F_i^{\mathbf A}(\mathbf
              x)) \;< \;\varepsilon_{M-1}
\end{align}
for each $\mathbf x\in A^{n_i}$, then for each term $\tau$ of depth
$\leq M$, we have
\[
      d(\tau^{\mathbf B}(\mathbf x),\tau^{\mathbf A}(\mathbf x))
                   \;<\; \varepsilon
\]
for each $\mathbf x\in A^{\omega}$.
\end{lemma} \begin{Proof}
We will prove, by induction on $k$, that if $\sigma$ is a term (in
this language) of depth $\leq k$, then
\begin{align}   \label{eq-inner-ind-cont}
      d(\sigma^{\mathbf B}(\mathbf x),\sigma^{\mathbf A}(\mathbf x))
                   \;<\; \varepsilon_{M-k}
\end{align}
for each $\mathbf x\in A^{\omega}$. The lemma then follows from the
truth of (\ref{eq-inner-ind-cont}) for $k=M$.

For $k=1$, the inductive assertion (\ref{eq-inner-ind-cont}) is
immediate from (\ref{eq-ops-M-1-close}). We now suppose that the
assertion holds for $k\geq1$ and prove it for $k+1$. So suppose that
$\sigma$ is a term of depth $\leq k+1$. We may clearly assume that
$\sigma=F_i(\sigma_1,\cdots,\sigma_{n_i})$ for some $i<N$ and some
terms $\sigma_j$ of depth $\leq k$. By an inductive appeal to
(\ref{eq-inner-ind-cont}), for each $j$  we have
\begin{align*}
      d(\sigma_j^{\mathbf B}(\mathbf x),\sigma_j^{\mathbf A}(\mathbf x))
                   \;<\; \varepsilon_{M-k}
\end{align*}
for each $\mathbf x\in A^{\omega}$.  To calculate
(\ref{eq-inner-ind-cont}) for $\sigma$, we begin with the triangle
inequality, obtaining
\begin{align*}
      d(\sigma^{\mathbf B}(\mathbf x),\sigma^{\mathbf A}(\mathbf x))
        &\EQ d(F^{\mathbf B}(\sigma_1^{\mathbf B}(\mathbf x),\ldots),
        F^{\mathbf A}(\sigma_1^{\mathbf A}(\mathbf
        x),\ldots))\\[0.1in]
        &\;\leq\;d(F^{\mathbf B}(\sigma_1^{\mathbf B}(\mathbf x),\ldots),
        F^{\mathbf A}(\sigma_1^{\mathbf B}(\mathbf x),\ldots))\\
        &\quad\quad\quad\quad\quad\quad
        +\;
        d(F^{\mathbf A}(\sigma_1^{\mathbf B}(\mathbf x),\ldots),
        F^{\mathbf A}(\sigma_1^{\mathbf A}(\mathbf x),\ldots))\\[0.1in]
        &\;\;\leq \;\varepsilon_{M-1} \;+\;
                 \frac{\varepsilon_{M-k-1}}{2}\;\leq \;
                    \frac{\varepsilon_{M-k-1}}{2}\;+\;
                    \frac{\varepsilon_{M-k-1}}{2}
                    \EQ \varepsilon_{M-k-1}\;.
\end{align*}
(In the final line, the first term of the estimate (i.e.\
$\varepsilon_{M-1}$) comes from (\ref{eq-ops-M-1-close}), and the
second term (i.e.\ $\varepsilon_{M-k-1}$) comes from
(\ref{eq-rec-eps-k}) and (\ref{eq-inner-ind-cont}). The final
inequality comes from (\ref{eq-rec-def-epsk}).
\end{Proof}

\subsection{Approximate satisfaction by simplicial operations.}
Suppose that $A=|K|$, the geometric realization of a complex $K$. We
call an operation $G\FROM A^n\TO A$ {\em simplicial} if there exist
$M$, $N$ and an $(M,N)$-simplicial map $\phi\FROM (K^n)^{(M)}\TO
L^{(N)}$ such that $G\EQ |\phi|_0\Compos\INV{(|\pi|^n{})}$. (For
notation, see Theorem \ref{th-SAT-op} on page \pageref{th-SAT-op},
and material preceding Theorem \ref{th-SAT-op}.)

\begin{corollary}    \label{cor-tau-close}
Let $A=|K|$, with $K$ finite, and let $\mathbf A=(A,
F_i^{\mathbf A})_{i<N}$ be a topological algebra based on $A$. For
each $M\in\Integers^+$ and each real $\varepsilon>0$, there exist
simplicial operations $F_i^{\mathbf B}\FROM A^{n_i}\TO A$ such that
each term $\tau$ of depth $\leq M$ in the operation symbols $F_i$
satisfies
\[
      d(\tau^{\mathbf B}(\mathbf x),\tau^{\mathbf A}(\mathbf x))
                   \;<\; \varepsilon
\]
for each $\mathbf x\in A^{\omega}$.
\end{corollary}\begin{Proof}
Clearly $A$ is compact, and so we may appeal to Theorem
\ref{lem-approx-term}. Let $\varepsilon_{M-1}$ be as supplied by
Theorem \ref{lem-approx-term} (in other words, coming from Equations
(\ref{eq-rec-eps-k}--\ref{eq-rec-def-epsk})). By Theorem
\ref{th-SAT-op}, for each operation $F_i^{\mathbf A}$, there is a
simplicial map $\phi_i$ such that
$|\phi_i|_0\Compos\INV{(|\pi|^n{})}$ approximates $F_i^{\mathbf A}$
within $\varepsilon_{M-1}$. In other words, if we define
$F_i^{\mathbf B}$ to be the simplicial operation
$|\phi_i|_0\Compos\INV{(|\pi|^n{})}$, then we have Equation
(\ref{eq-ops-M-1-close}) holding for all $\mathbf x\in A^{n_i}$. The
desired conclusion is now immediate from Theorem
\ref{lem-approx-term}.
\end{Proof}

Recall from \S\ref{sub-metric} that
\begin{align*}
         \lambda_{\mathbf A}(\sigma,\tau) &\EQ
              \sup \bigl\{\,d(\sigma^{\mathbf A}({\mathbf a}),
              \tau^{\mathbf A}({\mathbf a}))
              \,:\, {\mathbf a}\in A^{\omega} \bigr\}
              \;\in\; \Reals^{\geq0}\cup\{\infty\};\\
    \lambda_{\mathbf A}(\Sigma) &\EQ \sup\bigl\{\lambda_{\mathbf A}
      (\sigma,\tau)\,:\, \sigma\wavy\tau\,\in \Sigma \bigr\}
      \;\in\; \Reals^{\geq0}\cup\{\infty\};\\
     \lambda_A(\Sigma) &\EQ \inf\bigl \{\lambda_{\mathbf
     A}(\Sigma) \,:\,{\mathbf A} =(A;\o F_t)_{t\in T},
       \;\text{$\o F_t$ any continuous operations} \bigr\}.
\end{align*}

\begin{lemma}\label{lem-approx-sig-tau}
If\/ $A$ is a metric space, and\/ $\mathbf A=(A, F_t^{\mathbf
A})_{t\in T}$ and $\mathbf B=(A, F_i^{\mathbf B})_{i<N}$ are two
similar algebras based on $A$, and if\/ $\sigma$ and\/ $\tau$ are
terms in the language of $\mathbf A$ satisfying
\[
      d(\sigma^{\mathbf B}(\mathbf x),\sigma^{\mathbf A}(\mathbf x))
                   \;<\; \varepsilon\,;\quad
      d(\tau^{\mathbf B}(\mathbf x),\tau^{\mathbf A}(\mathbf x))
                   \;<\; \varepsilon
\]
for all $\mathbf x\in A^{\omega}$, then $\lambda_{\mathbf
B}(\sigma,\tau)\leq\lambda_{\mathbf
A}(\sigma,\tau)\,+\,2\,\varepsilon$.\ENDPROOF
\end{lemma}

\begin{theorem}\label{th-lambda-approx-simp}
Suppose that $A=|F|$, the geometric realization of a finite complex,
and suppose that\/ $\Sigma$ is a finite set of equations of type
$\la F_t\,:\,t\in T\ra$. For each topological algebra $\mathbf A=(A,
F_t^{\mathbf A})_{t\in T}$ of this type based on $A$, and for each
real\/ $\varepsilon>0$, there exists an algebra $\mathbf B=(A,
F_t^{\mathbf B})_{t\in T}$, each of whose operations is a simplicial
operation, and such that
\[
        \lambda_{\mathbf B}(\Sigma) \;\leq\; \lambda_{\mathbf A}(\Sigma)
                     \,+\,\varepsilon.
\]
\end{theorem}\begin{Proof}
Let $M$ be the maximum depth of terms appearing in $\Sigma$. Take
$M$ and $\varepsilon/2$ to Corollary \ref{cor-tau-close}, and let
$F_i^{\mathbf B}\FROM A^{n_t}\TO A$ be the simplicial operations
that it yields. Thus Corollary \ref{cor-tau-close} yields
\[
      d(\sigma^{\mathbf B}(\mathbf x),\sigma^{\mathbf A}(\mathbf x))
                   \;<\; \varepsilon/2\,;\quad
      d(\tau^{\mathbf B}(\mathbf x),\tau^{\mathbf A}(\mathbf x))
                   \;<\; \varepsilon/2
\]
for $\mathbf x\in A^{\omega}$ and for each equation
$\sigma\wavy\tau$ of $\Sigma$. Then Lemma \ref{lem-approx-sig-tau}
yields $\lambda_{\mathbf B}(\sigma,\tau)\leq\lambda_{\mathbf
A}(\sigma,\tau)\,+\,\varepsilon$ for each equation $\sigma\wavy\tau$
of $\Sigma$. Collecting the individual equations into one $\Sigma$
yields the desired result.
\end{Proof}

In performing the infimum that occurs in the definition of
$\lambda_A$ for a metric space $A$, we might wish to restrict our
attention to simplicial operations. In other words, we may define
\begin{align*}
     \lambda_A^{\text{simp}}(\Sigma) \EQ \inf\bigl \{\lambda_{\mathbf
     A}(\Sigma) \,:\,{\mathbf A} =(A;\o F_t)_{t\in T},
       \;\text{$\o F_t$ any simplicial operations} \bigr\}.
\end{align*}
Clearly $\lambda_A^{\text{simp}}(\Sigma) \,\geq\, \lambda_A(\Sigma)$
for any $\Sigma$ and for $A$ the geometric realization of any
complex. In one special case we have equality:

\begin{corollary}    \label{cor-simp-lambda}
Suppose that $A=|F|$, the geometric realization of a finite complex,
and suppose that\/ $\Sigma$ is a finite set of equations. Then $$
\lambda_A^{\text{simp}}(\Sigma) \EQ \lambda_A(\Sigma).\ENDPROOF$$
\end{corollary}

\section{Algorithmic considerations.}  \label{sec-alg}
Recall from \cite{wtaylor-eri} that there is no algorithm to
determine whether a finite $\Sigma$ is compatible with $\Reals$. In
fact, the set of finite $\Sigma$ that are not $\Reals$-compatible
fails to be recursively enumerable \cite{wtaylor-eri}.

Our main tool will be an algorithm $\mathcal T$ that decides the
truth in $\la\Reals;+,\cdot,0,1\ra$ of all first-order sentences in
$+$, $\cdot$, 0, and 1. Such an algorithm was devised by Alfred
Tarski in 1931
--- see \cite{tarski-1931,tarski-1948,tarski-1951} --- with
improved versions devised later by G. E. Collins in 1975 (see
\cite{collins}), and by several others in later years (see e.g.\
\cite{caviness}). (Notice that $x\leq y$ can be replaced in this
context by $\exists z(y=x+z^2)$, so the algorithm $\mathcal T$ can
also work with inequalities.)

\subsection{The main simplicial algorithm} \label{sub-main-simp-alg}

We shall assume that we have available an input language that
accommodates the description of a finite sets of sets (e.g.\ a
simplicial complex), a finite similarity type $(n_1,\cdots,n_k)$,
and equations in operation symbols $F_1$ ($n_1$-ary),\ldots,$F_k$
($n_k$-ary).
\begin{theorem} \label{th-alg-root}
There exists an algorithm $\mathcal A$ with the following
behavior. $\mathcal A$ accepts as input positive integers $M$,
$N$, $r$ and $s$, a finite complex $K$, and a finite
equation-set\/ $\Sigma$. $\mathcal A$ yields an answer to the
question, do there exist $M,N$-simplicial maps on the appropriate
powers of\/ $|K|$ so that the equations of\/ $\Sigma$ are
satisfied within $\,<\,r/s$?
\end{theorem}\begin{Proof}
We consider only a single equation $\sigma\wavy\tau$; the extension
to a finite set offers no important complications.
\begin{itemize}
\item Adopt an algorithmic system for the handling of complexes,
their products (\S\ref{sub-prod}), their subdivision
(\S\ref{sub-bary}), and their geometric representations
(\S\ref{sub-simpmaps}). For use with the latter, we shall require
real-vector calculations within any specified tolerance. In
particular the system must be capable of a symbolic modeling of the
computations seen in \S\ref{sub-simp-approx}: the ground-level
geometric realization $|\phi|_0$ defined in Equation
(\ref{eq-phi-0-def}), the function $|\pi|^n$ described in Equation
(\ref{eq-pi-n}), and composites involving these, as used in Theorem
\ref{th-SAT-op}.
\item Establish a list --- $F_1$($n_1$-ary),\ldots,$F_k$
 ($n_k$-ary) --- of the operation symbols that appear either in
$\sigma$ or in $\tau$, and their arities.
\item Establish a representation in this system for the complex
$K^{(N)}$, and for each complex $(K^{n_j})^{(M)}$ (for $1\leq j\leq
k$).
\item Loop through all of the finitely many $k$-tuples
$(\phi_1,\cdots,\phi_k)$, where each $\phi_j$ maps the vertices of
$(K^{n_j})^{(M)}$ to the vertices of $K^{(N)}$. For each $k$-tuple,
loop through these instructions:
   \begin{itemize}
   \item Examine each $\phi_j$ to see whether it is a simplicial
   map. If all are simplicial, continue this loop; if one is not
   simplicial, jump to the next $k$-tuple of maps and return to
   this instruction.
   \item Represent each $|\phi_j|\,\Compos\,\INV{(|\pi|^{n_j})}$ ---
   as in Theorem \ref{th-SAT-op} --- as a piecewise affine map with
   unknown real coefficients $\alpha_u^v$. (With $u$ ranging over all
   simplices of $(K^{n_j})^{(M)}$, and $v$ ranging over $\{0,\cdots,
   \dim(u)\}$.)
   \item Recursively realize each subterm of $\sigma$ by a piecewise
   affine map with unknown coefficients. These unknown coefficients
   may be expressed as ring-theoretic combinations of the
   $\alpha_u^v$. The same is to be done for $\tau$. Now distances
   between $\sigma(\mathbf x)$ and $\tau(\mathbf x)$ may be computed
   on each simplex of the subdivision, in terms of the unknown
   coefficients $\alpha_u^v$.
   \item Finally we express a question of whether this distance
   can made less than $r/s$ over each simplex. This question may
   be expressed as the existential closure of a conjunction of
   ring-theoretic inequalities in the unknowns $\alpha_u^v$.
   \item Using Tarski's algorithm $\mathcal T$, we determine the truth
   in $\Reals$ of this existential sentence. If the answer is affirmative,
   we terminate the
   algorithm, answering yes to the question of satisfiability
   within $\,<\,r/s$.
   \end{itemize}
\item If we reach this point, having never answered yes, we
terminate the algorithm, answering no the question of satisfiability
within $\,<\,r/s$.
\end{itemize}
It is obvious that the algorithm terminates at some point. The only
way for it to terminate with a ``yes'' answer is if, in the
penultimate instruction, Tarski's algorithm tells us that a certain
family of piecewise-linear operations, namely $|\phi_j|\,\Compos\,
\INV{(|\pi|^{n_j})}$ ---    as in Theorem \ref{th-SAT-op} --- yields
a topological algebra that satisfies $\sigma\wavy\tau$ within
$\,<\,r/s$. On the other hand, if a simplicial map exists allowing
satisfaction within $r/s$, then it must have been considered, and
Tarksi algorithm must have answered ``yes.'' Hence a ``yes'' answer
must have been obtained by our algorithm. This shows that the
algorithm is correct.
\end{Proof}

For the next corollary, let us fix a list of operation symbols $F_i$
($i\in\omega$), which includes each arity infinitely often.
\begin{corollary}   \label{cor-enumerate-sext}
There exists an algorithm $\mathcal B$ that takes no input and whose
output is an infinite sequence of sextuples $\la K, M, N, r,s,
\Sigma\ra$. Each sextuple satisfies
\begin{enumerate}
\item $K$ is a finite simplicial complex.
\item $M$, $N$, $r$ and $s$ are positive integers.
\item $\Sigma$ is a finite set of equations, whose operation symbols are
    among the $F_i$.
\item There exist $M,N$-simplicial maps on the appropriate powers
    of\/ $|K|$ so that the equations of\/ $\Sigma$ are satisfied
    within $\,<\,r/s$.
\end{enumerate}
Moreover every sextuple satisfying (1)-(4) is in the output of the
algorithm $\mathcal B$.
\end{corollary} \begin{Proof}
It is well known that there exists an algorithm $\mathcal C$ that
lists all sextuples satisfying 1, 2 and 3. To obtain the desired
$\mathcal B$, we merely filter $\mathcal C$ using the algorithm
$\mathcal A$ of Theorem \ref{th-alg-root}. In other words
$\mathcal B$ successively takes each output of $\mathcal B$ and
passes it to $\mathcal A$, which returns an answer of whether
there exist $M,N$-simplicial maps on the appropriate powers of\/
$|K|$ so that the equations of\/ $\Sigma$ are satisfied within
$\,<\,r/s$ (Point 4). If the answer is ``no,'' $\mathcal B$ takes
no further action at that stage; if the answer is ``yes,'' then
$\mathcal B$ outputs the sextuple in question.
\end{Proof}

\subsection{Recursive enumerability of
$\lambda_{|K|}(\Sigma)<\alpha$}    \label{sub-re-ity}

For {\em computable} real numbers, the reader is referred to 
\cite{pour-richards}. All we shall require of computability is this:
if $\alpha$ is computable, then there is an algorithm to decide
$s\,\alpha>\,r$ for positive integers $r$ and $s$. Clearly every
rational is computable; hence the computable reals are dense in
$\Reals$.

As in Corollary \ref{cor-enumerate-sext}, we fix a list of operation
symbols $F_i$ ($i\in\omega$), which includes each arity infinitely
often. Clearly every finite set of equations is definitionally
equivalent to a finite set involving only the operation symbols
$F_i$.

\begin{corollary}              \label{cor-sigmas-re}
Let $K$ be a finite simplicial complex, with\/ $|K|$ its geometric
realization (as usual), and let\/ $\alpha>0$ be a computable real
number. There is an algorithm $\mathcal E_{K,\alpha}$ whose output
consists of those finite sets $\Sigma$ of equations in $F_i$
($i\in\omega$) for which $\lambda_{|K|}(\Sigma)\,<\,\alpha$.
\end{corollary} \begin{Proof}
We run the algorithm $\mathcal B$ of Corollary
\ref{cor-enumerate-sext}, and filter the output as follows. When
$\mathcal B$ outputs $\la K, M, N, r,s, \Sigma\ra$, the new
algorithm $\mathcal E_{K,\alpha}$ either outputs $\Sigma$ or rests.
If $s\,\alpha\,>\,r$, then $\mathcal E_{K,\alpha}$ outputs $\Sigma$;
if not, then $\mathcal E_{K,\alpha}$ has no output. It is immediate
from Corollary \ref{cor-enumerate-sext} that $\mathcal E_{K,\alpha}$
has the desired output.
\end{Proof}

\begin{corollary}   \label{cor-less-alpha-re}
Let $K$ be a finite complex and $\alpha>0$ a computable real number.
Restricting attention to equations in $F_i$ ($i\in\omega$), the set
of finite sets $\Sigma$ with $\lambda_{|K|}(\Sigma)\,<\,\alpha$ is
recursively enumerable.\ENDPROOF
\end{corollary}

\begin{corollary}   \label{cor-less-alpha-re-pairs}
Let $\alpha>0$ be a computable real number. Restricting attention
to equations in $F_i$ ($i\in\omega$), the set of pairs
$(K,\Sigma)$ where $K$ is a finite complex, $\Sigma$ is a finite
set of equations, and where $\lambda_{|K|}(\Sigma)\,<\,\alpha$, is
recursively enumerable.\ENDPROOF
\end{corollary}

%

\subsection{The arithmetic character of $\lambda_{|K|}(\Sigma)=0$.}
                          \label{sub-ar-char-zero}
As in Corollaries
\ref{cor-enumerate-sext}--\ref{cor-less-alpha-re}, we fix a list
of operation symbols $F_i$ ($i\in\omega$), which includes each
arity infinitely often. We also fix a syntax that describes
nothing but complexes and describes each complex at least once up
to isomorphism.
\begin{corollary}
There is an algorithm $\mathcal F$ that takes one finite simplicial
complex $K$ as input, and whose output is an infinite stream of
pairs $(\Sigma,s)$, with each\/ $\Sigma$ a finite set of equations
in the symbols $F_i$ ($i\in\omega$) and with each $s\in\Integers^+$,
such that the following condition holds. Such a finite set $\Sigma$
satisfies $\lambda_{|K|}(\Sigma)=0$ iff\/ $(\Sigma,s)$ occurs in the
output stream of\/ $\mathcal F_{K}$ for arbitrarily large $s$.
\end{corollary}\begin{Proof}
The algorithm $\mathcal F$ accepts $K$, and then runs the algorithm
$\mathcal B$ of Corollary \ref{cor-enumerate-sext}, while filtering
the output as follows. Suppose that one piece of output from
$\mathcal B$ is $\la K', M, N, r,s, \Sigma\ra$. If $K'\neq K$, or if
$r\neq 1$, this output is filtered out completely. On the other
hand, if $K'=K$ and $r=1$, then $\mathcal F$ outputs the pair
$(\Sigma,s)$.

To see the equivalence of $\lambda_{|K|}(\Sigma)=0$ with the
algorithmic condition, we suppose first that $\Sigma$ is a finite
set of equations for which $(\Sigma,s)$ appears with arbitrarily
large $s$. By Corollary \ref{cor-enumerate-sext}(4), for arbitrary
large $s$ there are simplicial operations on $|K|$ that cause the
equations $\Sigma$ to be satisfied within $1/s$. In other words
$\lambda_{|K|}(\Sigma)\leq 1/s$ for arbitrarily large $s$, i.e.\
$\lambda_{|K|}(\Sigma)=0.$

Conversely, if $\lambda_{|K|}(\Sigma)=0$, then by Corollary
\ref{cor-simp-lambda}, we have $\lambda_{|K|}^{\text{simp}}(\Sigma)
\allowbreak\EQ\allowbreak \lambda_{|K|}(\Sigma)=0$. Thus for
arbitrarily large $s$, there exist simplicial operations on $|K|$
that cause the equations $\Sigma$ to hold within $1/s$. Therefore
the output stream of $\mathcal B$ will contain $\la K, M, N, 1,s,
\Sigma\ra$ for some $M$ and $N$. Therefore the output stream of
$\mathcal F$, when started with $K$, will contain the pair
$(\Sigma,s)$.
\end{Proof}

\begin{corollary}     \label{cor-pi-2}
For a fixed finite simplicial complex $K$, the set of\/ $\Sigma$
with $\lambda_{|K|}(\Sigma)=0$ is a $\Pi_2$-set.\ENDPROOF
\end{corollary}


\section{Filters.}\label{sec-filters} For any metric space $Z$, we may define the
class of theories
\begin{align*}
            \mathcal L_Z \EQ \bigl\{\Sigma^{\ast}\,:\,
             \lambda_Z(\Sigma^{\ast}) >
            0\bigr\}.
\end{align*}
$\Sigma^{\ast}$ is in this class iff there exists real
$\varepsilon>0$ so that any continuous operations on $Z$ violate
$\Sigma^{\ast}$ by more than $\varepsilon$ at some point of
$Z^{\omega}$.  In this context, it is best to work exclusively
with deductively closed sets of equations; hence our use of the
notation $\Sigma^{\ast}$. By Lemma \ref{lemma-comp-two-metrics} of
\S\ref{sub-dep-metric}, if $Z$ is compact, then the class
$\mathcal L_Z$ is a topological invariant. Hence one may also
regard $\mathcal L_Z$ as well-defined for any compact metrizable
topological space $Z$.

From Theorem \ref{th-interp-lesseq} of \S\ref{sub-interpretable},
it is not hard to see that $\mathcal L_Z$ is an upward-closed
subclass of the class of equational theories, ordered according to
interpretability. That is, if $\Gamma^{\ast}$ is interpretable in
$\Sigma^{\ast}$, and if $\Gamma^{\ast}\in \mathcal L_Z$, then
$\Sigma^{\ast}\in \mathcal L_Z$.

For some special spaces $Z$, such as $Z=[0,1]$ we can also show
that $\mathcal L_Z$ is a filter (i.e.\ that it is also closed
under the meet operation for this ordering). As is well known, the
meet of theories $\Gamma$ and $\Delta$ in this context is the
product theory $\Gamma\times\Delta$ that was described at the
start of \S\ref{sec-prod-vars}. For some $Z$, $\mathcal L_Z$ is
not a filter, for example for $Z=S^1\times Y$, as was seen in
Theorem \ref{th-prod-smaller-factors} of \S\ref{sub-prod-on-prod}.

In keeping with the convention of \S\ref{sec-filters}, the
following theorem will be stated only for deductively closed
theories $\Gamma^{\star}$ and $\Delta^{\star}$. It is easily
proved (e.g.\ from Condition (ii) at the start of
\S\ref{sec-prod-vars}) that $(\Gamma\times\Delta)^{\star} =
(\Gamma^{\star}\times\Delta^{\star})^{\star}$

\begin{theorem}           \label{th-01-filter}
If\/ $\Gamma^{\star}, \Delta^{\star}\in\mathcal L_{[0,1]}$, then
$(\Gamma\times \Delta)^{\star}\in\mathcal L_{[0,1]}$. Therefore
$\mathcal L_{[0,1]}$ is a filter.
\end{theorem} \begin{Proof}
For a proof by contradiction, let us suppose that $(\Gamma\times
\Delta)^{\star}\not\in\mathcal L_{[0,1]}$. Then
$\lambda_{[0,1]}((\Gamma^{\star}\times\Delta^{\star})^{\star})
=0$, and hence
$\lambda_{[0,1]}(\Gamma^{\star}\times\Delta^{\star}) =0.$ By
Theorem \ref{th-I-approx-indec} of \S\ref{sub-prod-prodindec}, for
any real $\varepsilon$ with $0<\varepsilon<1/16$, we have
$\lambda_{[0,1]}(\Gamma^{\star}) <4\varepsilon$ or
$\lambda_{[0,1]}(\Delta^{\star}) <4\varepsilon$. This clearly
implies that $\lambda_{[0,1]}(\Gamma^{\star}) \,=\,0$ or
$\lambda_{[0,1]}(\Delta^{\star}) \,=\,0.$ This conclusion
contradicts our assumption that both
$\lambda_{[0,1]}(\Gamma^{\star})$ and
$\lambda_{[0,1]}(\Delta^{\star})$ are $>0$.
\end{Proof}

 For $Z$ product-indecomposable, there is a better-known filter,
that of all theories $\Sigma^{\ast}$ that are incompatible with
$Z$. Our filter $\mathcal L_Z$ is of course a subset of that
one---generally a proper subset---but no further relationship is
known at this time. Notice also that the argument for Theorem
\ref{th-01-filter} is quite general: for any space $Z$, if we have
a result like Theorem \ref{th-I-approx-indec} for $Z$, then we can
prove a result like Theorem \ref{th-01-filter} for $Z$. For
example Theorem \ref{th-Y-approx-indec} yields that $\mathcal L_Y$
is a filter, for $Y$ being the figure-Y space.

Further up-sets may be considered, for example using $\delta^-$ (of
\S\ref{sec-top-inv}) in place of $\lambda$. \vspace{0.1in}

\section{Approximate satisfaction by differentiable
operations.} \label{sub-approx-diff} Obviously one could make a
whole new list of definitions here; we refrain from this action
until it may be warranted by further developments. Instead, we
content ourselves with sketching one example.

Recall that we proved in \S8.1.1 of \cite{wtaylor-eri} that
semilattice theory is not $C^1$-compatible with $\Reals$. Here we
exhibit some $C^1$ approximants to the theory. For arbitrary real
$p>1$, and for arbitrary reals $a,b>1$, define
\[
       \o {\wedge_p}(a,b) \EQ \bigl( a^p \,+\,
       b^p\bigr)^{\frac{1}{p}}.
\]
It is not hard to verify that each $\o {\wedge_p}$ is $C^1$, and
\[
              \lim_{p\rightarrow\infty} \o {\wedge_p}(a,b)
                       \EQ a\wedge b,
\]
uniformly in $a$ and $b$. It is also not hard to show that for
every bounded interval [c,d], and for every real
$\varepsilon>0$ there exists $p$ such that
\[
          \bigl|\, \o {\wedge_p}(a,b) \,-\, a\wedge b\, \bigr|
                       < \varepsilon
\]
for all $a,b$ in $[c,d]$.


\newcommand{\AEQ}{{\em Aequationes Mathematicae}}
\newcommand{\AAM}{{\em Advances in Applied Mathematics}}
\newcommand{\AMS}{{American Mathematical Society}}
\newcommand{\AU}{{\em Algebra Universalis}}
\newcommand{\ANNALS}{{\em Annals of Mathematics}}
\newcommand{\AML}{{\em Annals of Mathematical Logic}}
\newcommand{\ANNALEN}{{\em Mathematische Annalen}}
\newcommand{\BAMS}{{\em Bulletin of the \AMS}}
\newcommand{\BPAN}{{\em Bulletin de l'Aca\-d\'{e}\-mie Po\-lo\-naise
    des Sciences, S\'{e}rie des sciences math., astr. et phys.}}
\newcommand{\CAMB}{{\em Proceedings of the Cambridge Philosophical Society}}
\newcommand{\CANAD}{{\em Canadian Journal of Mathematics}}
\newcommand{\COLLOQ}{{\em Col\-lo\-qui\-um Ma\-the\-ma\-ti\-cum}}
\newcommand{\COLLOQUIA}{{\em Col\-lo\-qui\-a Ma\-the\-ma\-ti\-ca
    Soc\-i\-e\-ta\-tis Bolyai J\'{a}nos}}
\newcommand{\DM}{{\em Discrete Mathematics}}
\newcommand{\EM}{{\em l'En\-seigne\-ment Ma\-th\'{e}\-matique}}
\newcommand{\FM}{{\em Fun\-da\-men\-ta Ma\-the\-ma\-ti\-cae}}
\newcommand{\HOUS}{{\em Hous\-ton Journal of Mathematics}}
\newcommand{\INDAG}{{\em In\-da\-ga\-ti\-o\-nes Ma\-the\-ma\-ti\-cae}}
\newcommand{\JALG}{{\em Journal of Algebra}}
\newcommand{\JAMS}{{\em Journal of the Australian Mathematical Society}}
\newcommand{\JPAA}{{\em Journal of Pure and Applied Algebra}}
\newcommand{\JPAL}{{\em Journal of Pure and Applied Logic}}
\newcommand{\JSL}{{\em J. Symbolic Logic}}
\newcommand{\LNM}{{\em Lecture Notes in Mathematics}}
\newcommand{\MONTHLY}{{\em American Mathematical Monthly}}
\newcommand{\MUSSR}{{\em Mathematics of the {\sc USSR} --- Sbor\-nik}}
\newcommand{\NAMS}{{\em Notices of the \AMS}}
\newcommand{\NORSK}{{\em Norske Vid. Selssk. Skr. I, Mat. Nat.
                Kl. Chris\-ti\-a\-nia}}
\newcommand{\ORD}{{\em Order}}
\newcommand{\PACIFIC}{{\em Pacific Journal of Mathematics}}
\newcommand{\PAMS}{{\em Proceedings of the \AMS}}
\newcommand{\REM}{{\em Research and Exposition in Mathematics}}
\newcommand{\SF}{{\em Semigroup Forum}}
\newcommand{\SCAND}{{\em Ma\-the\-ma\-ti\-ca Scan\-di\-na\-vi\-ca}}
\newcommand{\SIAMJC}{{\em {\sc Siam} Journal of Computing}}
\newcommand{\SZEGED}{{\em Acta Sci\-en\-ti\-a\-rum Ma\-the\-ma\-ti\-ca\-rum}
                (Sze\-ged)}
\newcommand{\TAMS}{{\em Transactions of the \AMS}}
\newcommand{\TCS}{{\em Theoretical Computer Science}}

\newcommand{\ALVIN}{{\em Algebras, Lattices, Varieties}}
\newcommand{\Wads}{{Wads\-worth and Brooks-Cole Publishing Company,
                    Monterey, CA}}
\newcommand{\NorthHolland}{{North-Holland Publishing Company, Amsterdam}}
\newcommand{\Birk}{{Birkh\"{a}user}}
\newcommand{\Garcia}{Garc\'{\i}a}

\vspace{\fill} \hspace*{-\parindent}%
\parbox[t]{3.0in}{ Walter Taylor \\ Mathematics Department\\ University
of Colorado\\
 Boulder, Colorado \ 80309--0395\\ USA\\ Email:
 {\tt walter.taylor@colorado.edu}}
\end{document}